\documentclass[a4paper,12pt]{amsart}

\usepackage{graphicx,amsmath,amsfonts,amsthm,amssymb,url,dsfont,bbm}

\usepackage[OT2,T1]{fontenc}
\DeclareSymbolFont{cyrletters}{OT2}{wncyr}{m}{n}
\DeclareMathSymbol{\Sha}{\mathalpha}{cyrletters}{"58}

\def\]{\textup{\mbox{]\hspace{-.15em}]}}}
\def\[{\textup{\mbox{[\hspace{-.15em}[}}}

\def\got{\mathfrak}

\newenvironment{pf}
{\medskip\noindent {\it Proof --- \ }}
{\hfill\nobreak $\Box$ \par\bigbreak}

\newcommand{\F}{{\mathbb F}}
\newcommand{\HH}{{\mathcal H}}

\newcommand{\isomo}{\overset{\sim}{\rightarrow}}
\newcommand{\cal}{\mathcal}

\newcommand{\WW}{\mathcal{W}}
\newcommand{\TT}{\mathcal{T}}
\newcommand{\Hom}{\mathrm{Hom}}

\newcommand{\Gal}{\mathrm{Gal}}
\newcommand{\G}{\mathrm{G}}

\newcommand{\pn}{\par \noindent}

\newcommand{\Qpb}{\overline{\mathbb{Q}}_p}

\newcommand{\GL}{\mathrm{GL}}
\newcommand{\SL}{\mathrm{SL}}
\renewcommand{\U}{\mathrm{U}}

\newcommand{\ps}{\par \smallskip}

\newcommand{\End}{\mathrm{End}}
\newcommand{\N}{\mathbb{N}}
\newcommand{\Z}{\mathbb{Z}}

\newcommand{\Fpb}{\overline{\mathbb{F}}_p}
\newcommand{\Q}{\mathbb{Q}}
\newcommand{\Qb}{\overline{\Q}}
\newcommand{\Qp}{\mathbb{Q}_p}
\newcommand{\R}{\mathbb{R}}
\newcommand{\C}{\mathbb{C}}

\newcommand{\rhob}{{\bar\rho}}

\newcommand{\Ker}{\mathrm{Ker}}
\newcommand{\disc}{\mathrm{disc}}

\newcommand{\AAA}{\mathbb{A}}
\newcommand{\OO}{\mathcal{O}}

\newcommand{\D}{\mathrm{D}}

\newcommand{\arb}{{\rm Ad'}(\rhob)}
\newcommand{\cB}{{\mathcal B}}

\newcommand{\Ref}{\cal{F}}
\newcommand{\fg}{(\varphi,\Gamma)}
\newcommand{\Ro}{\mathcal{R}}

\newcommand{\pr}{{\mathrm{pr}}}

\newcommand{\Nob}{\overline{N}_0}
\newcommand{\Cc}{\mathcal{C}}
\newcommand{\ATL}{\mathfrak{A}}
\newcommand{\Spf}{{\rm Spf}}
\newcommand{\rk}{{\mathrm{rk}}}
\newcommand{\gr}{{\rm gr}}
\newcommand{\cDc}{\mathcal{D}_{\rm crys}}
\newcommand{\cDS}{\cD_{\rm Sen}}

\newcommand{\PSen}{{{\rm P}_{\rm Sen}}}
\newcommand{\fga}{\widehat{\mathbb G}_a}
\newcommand{\pht}{\partial \PSen}
\newcommand{\tg}{{\rm t}}

\newtheorem{thm}[subsection]{Theorem}
\newtheorem{lemma}[subsection]{Lemma}
\newtheorem{remark}[subsection]{Remark}
\newtheorem{cor}[subsection]{Corollary}
\newtheorem{prop}[subsection]{Proposition}
\newtheorem{example}[subsection]{Example}
\newtheorem{question}[subsection]{Question}

\newtheorem{conj}[subsection]{Conjecture}

\newtheorem{definition}[subsection]{Definition}

\newcommand{\Dc}{D_{\rm crys}}
\newcommand{\cE}{{\mathcal{E}}}
\newcommand{\cF}{{\mathcal{F}}}
\newcommand{\cX}{{\mathfrak X}}
\newcommand{\cZ}{{\mathcal Z}}
\newcommand{\cXg}{\cX^{\rm g}}
\newcommand{\cXc}{\cX^{\rm c}}
\newcommand{\cXm}{\cX^{\rm mod}}
\newcommand{\Symm}{{\rm Symm}}

\newcommand{\cXem}{{\cX^{\rm emod}}}
\newcommand{\cXgem}{{\cX^{\rm gen}}}
\newcommand{\adr}{{\rm Ad}'}
\newcommand{\Sym}{{\rm Sym}} 

\title{On the infinite fern of Galois representations of unitary type}

\begin{document}

\author[G. Chenevier]{Ga\"etan Chenevier}
\email{chenevier@math.polytechnique.fr}
\address{Ga\"etan Chenevier\\ C.N.R.S, Centre de Math\'ematiques Laurent Schwartz, \'Ecole Polytechnique, 91128 Palaiseau Cedex\\ France}
\maketitle

\begin{flushright} {\it To Julia and Valeria.} \end{flushright}

\tableofcontents

\section*{Introduction}

Let\footnote{November 30, 2009. The author is supported by the C.N.R.S.} $E$ be a number field, $p$ a prime number, $S$ a finite set of places of
$E$ containing the places above $p$ and $\infty$, $\G_{E,S}$ the Galois
group of a maximal algebraic extension of $E$ unramified outside $S$, and let $n\geq 1$ be an integer.
We are interested in the set of $n$-dimensional, semisimple, continuous
representations $$\rho : \G_{E,S}\longrightarrow \GL_n(\Qpb)$$
taken up to isomorphism. This set turns out to be the $\Qpb$-points of a rigid analytic 
space $\cX$ (or $\cX_n$) over $\Qp$ in a natural way\footnote{We will not use this space in the sequel, but for its definition see \cite[Thm D]{chdet}.}. An interesting subset 
$$\cXg \subset \cX(\Qpb)$$ is the subset of representations which are {\it geometric}, in the sense they occur as a
subquotient of $H^i_{et}(X_{\overline{E}},\Qpb)(m)$ for some 
proper smooth algebraic variety $X$ over $E$ and some integers $i\geq 0$ and $m \in
\Z$. There are several basic open questions that we can ask about 
$\cX$ and its locus $\cXg$, here a some of them : 
\ps \medskip
\noindent {\bf Questions : } Does $\cXg$ possess some specific structure inside $\cX$ ? What can we say about its
various closures in $\cX$ ? (for example, for the Zariski or the $p$-adic topologies
)  What if we replace $\cXg$ by its subset $\cXc$ of $\rho$'s which are crystalline at the places of $E$ above 
$p$ ? 

\ps\medskip

Regarding the first question, a trivial observation is that $\cXg$ is countable, as so are
algebraic varieties over $E$, thus it certainly contains no analytic subset of $\cX$.
In the simplest case $E=\Q$ and $S=\{\infty,p\}$ then $\cX_1$ is 
the space of $p$-adic continuous
characters of $\Z_p^*$ (a finite union of $1$-dim open balls) 
and $\cXc$ is the subset of characters of the form $x \mapsto x^m$ for $m \in \Z$, 
which is Zariski-dense in $\cX_1$. For a general $E$
and $S$, we leave as an exercise to the reader to check that class field theory
and the theory of complex multiplication show that $\cXc$ is also
Zariski-dense in $\cX_1$ assuming Leopold's conjecture for $E$ at $p$. \ps

As a second and much more interesting example, let us recall the discovery of Gouv\^ea and Mazur \cite{gm}. They 
assume that $d=2$, $E=\Q$, and say $S=\{\infty,p\}$ to simplify. Let $q$ be a power of an odd prime $p$ and let $$\rhob : \G_{\Q,S}
\longrightarrow \GL_2(\F_q)$$
be an absolutely irreducible odd Galois representation. Let $R(\rhob)$ denote the universal odd $\G_{\Q,S}$-deformation ring
of $\rhob$ in the sense of Mazur and let $\cX(\rhob)$ be its 
analytic generic fiber : it is the connected component of $\cX_2$ parameterizing the $\rho$ 
with residual representation $\rhob$. In general $\cX(\rhob)$ is a rather complicated space, 
and Mazur first studied in the {\it unobstructed case\footnote{The philosophy of special 
values of $L$-functions suggests that this unobstructed case is the generic situation. It 
is now known for instance that for $\rhob=\rhob_\Delta$ attached to Ramanujan's 
$\Delta=\sum_{n>0} \tau(n)q^n$, $p > 13$  and $p\neq 691$, the deformation problem of 
$\rhob$ is unobstructed (Mazur, Weston).}} $H^2(\G_{\Q,S},{\rm ad}(\rhob))=0$,
for which class field theory shows that $R(\rhob) \simeq \Z_q[[x,y,z]]$, hence $\cX(\rhob)$ 
is the open unit ball of dimension $3$ over $\Q_q$. In this case, Gouv\^ea and Mazur showed that 
{\it $\cXc$ is Zariski-dense in $\cX(\rhob)$}. They actually show that the subset 
$$\cX^{\rm mod}\subset  \cX(\rhob)$$ of $p$-adic Galois representations $\rho_f(m)$ 
attached to an eigenform $f \in S_k(\SL_2(\Z))$ for some weight $k$, and some $m \in \Z$, 
is Zariski-dense in $\cX(\rhob)$. This subset is non empty as $\rhob$ is modular (Khare-Wintenberger).

Their proof relies heavily on the
theory of $p$-adic families of modular eigenforms due to
Coleman, extending pioneering works of Hida, that we briefly recall. 
Let $f = q+a_2q^2+a_3q^3+ \dots \in \Qpb[[q]]$ be a classical modular 
eigenform of level $1$, of some weight $k$, and such that $\rhob_f=\rhob$; 
$\rho_f$ corresponds to some $x_f \in \cX(\rhob)(\Qpb)$. Attached to $f$ we 
have two $p$-Weil numbers of weight $k-1$ which are the roots of the 
polynomial $X^2-a_pX+p^{k-1}$. The main result of Coleman asserts that 
if $\varphi$ is one of them, we can attach to $(f,\varphi)$ an affinoid 
subcurve $$C_{(f,\varphi)} \subset \cX(\rhob)$$ with the property that 
$C_{(f,\varphi)}$ contains a Zariski-dense subset of modular points 
$x_{f'}$, where furthermore $f'$ has a $p$-Weil number with the same $p$-adic 
valuation as $\varphi$. The infinite fern of Gouv\^ea-Mazur is by definition 
the union of all the $C_{(f,\varphi)}$ for all $f$ and choice of $\varphi$. A 
simple but important observation made by Gouv\^ea and Mazur is that 
$C_{(f,\varphi)}\cap C_{(f,\varphi')}$ necessarily is finite if $\varphi$ 
and $\varphi'$ have different valuations: this follows from the previous 
property and the fact that the ``weight'' varies analytically in Coleman's 
families. From this ``fractal'' picture it follows at once that the Zariski-closure of the 
modular points, or which is the same the Zariski-closure of the fern, 
has dimension at least $2$ inside $\cX(\rhob)$, and a simple argument of Tate-twists 
using Sen's theory gives then the result. The story does not quite end here as some 
years later, Coleman and Mazur defined a wonderful object, {\it the eigencurve}, which 
shed new lights on the infinite fern. They define a {\it refined modular point} as a pair 
$$(x,\varphi) \in  \cX(\rhob) \times {\mathbb G}_m$$
where $x=x_f$ is modular and $\varphi$ is a $p$-Weil number of $f$. The interesting fact is that the Zariski-closure of the refined modular points in $\cX(\rhob)\times \mathbb{G}_m$ is an {\it equidimensional curve}, the so-called $\rhob$-{\it eigencurve}. Its image in $\cX(\rhob)$ is the {\it complete infinite fern}, which simultaneously analytically continues each leaf of the infinite fern itself. This picture for $\cX(\rhob)$ provides a rather satisfactory answer towards the first of the main questions above, even though very little is known about the geometry of the eigencurve at present. It is believed that $\cE(\rhob)$ has only finitely 
many irreducible components. An amazing consequence of this conjecture would be that for a well chosen modular point $x \in \cX(\rhob)$ the analytic continuation of a well chosen 
leaf at $x$ would be Zarski-dense in $\cX(\rhob)$ ! However, as far as we know there is no non-trivial case in which this conjecture or its variants in other dimensions $>1$ is known.

Let us mention that when $S$ is general and $\rhob$ is possibly obstructed,
the approach above of Gouv\^ea-Mazur still shows that each irreducible
component of the Zariski-closure of the modular points in $\cX(\rhob)$ has
dimension at least $3$. In a somehow opposite direction, it is conjectured that in all cases 
$\cX(\rhob)$ is equidimensionnal of dimension $3$, and that each of its irreducible components
contains a smooth modular point. This conjecture, combined with the result of
Gouv\^ea-Mazur, implies the Zariski-density of the modular points in
$\cX(\rhob)$ for each odd $\rhob$ (say absolutely irreducible). Relying
on $R=T$ theorems of Taylor-Wiles, Diamond et al., Boeckle was able to show
that conjecture under some rather mild assumption on $\rhob$, hence the
Zariski-density in most cases: we refer to~\cite{boeckle} for the precise
statements.

Our main aim in this paper is to study a generalization of this picture to the higher dimensional case. Our most complete results will concern some pieces of $\cX_3$ satisfying some sort of self-duality condition. Let $E$ be a CM field\footnote{Throughout this paper, a CM field is assumed to be imaginary.},
and let $$\rhob : \G_{E,S}
\longrightarrow \GL_3(\F_q)$$
be an absolutely irreducible Galois representation such that $\rhob^* \simeq
\rhob^c$, where $c$ is a generator of $\Gal(E/F)$ and $F$ the maximal totally real subfield of $E$. Let $\cX(\rhob) \subset
\cX_3$ denote the closed subspace of $x \in \cX_3$ such that 
$\rho_x^* \simeq \rho_x^c$ and $\overline{\rho_x} \simeq \rhob$. This $\cX(\rhob)$ has conjectural
equidimension $6[F:\Q]$, and under an unobstructedness assumption similar to Mazur's one it is 
actually an open unit ball over $\Q_q$ in that number of variables. There is a natural notion of 
modular points in $\cX(\rhob)$: they are the $x$ such that $\rho_x$ is a $p$-adic Galois representation
attached to a cuspidal automorphic representation $\Pi$ of $\GL_3(\AAA_E)$ such that $\Pi_\infty$ is 
cohomological, $\Pi^\vee \simeq \Pi^c$, and such that for $v$ finite dividing $p$ or outside $S$, 
$\Pi_v$ is unramified. Those Galois representations have been constructed by
Rogawski, they are cut out from the \'etale cohomology of (some abelian
varieties over) the Picard modular surfaces
and they are related to automorphic forms on unitary groups in $3$ variables associated to $E/F$. 

\bigskip

{\bf Theorem A:} {\it Assume that $p$ is totally split in $E$. Then each irreducible component of the Zariski-closure of the
modular points in $\cX(\rhob)$ has dimension at least $6[F:\Q]$. 

In particular, in the unobstructed case the set of modular points of $\cX(\rhob)$ is Zariski-dense if it is non-empty.}

\bigskip

In the appendix, we give several examples of elliptic curves $A$ over $\Q$
such that the deformation problem {\it of type $U(3,E/\Q)$} of $\rhob := {\rm Sym^2} A[p]_{|\G_E}(-1)$
is unobstructed for $p=5$ and $E=\Q(i)$. As is the work of \cite{boeckle} in the $\GL(2,\Q)$ case, we expect that combining Theorem A with suitable $R=T$ theorems (as in~\cite{cht} and their followers), 
one should be able to remove the unobstructnedness assumption under suitable assumptions on $\rhob$. 
However, we postpone this to a subsequent study.

As in the work of Gouv\^ea-Mazur, a very important ``constructive'' ingredient of our proof is the theory 
of families of modular forms (for $U(3,E/F)$ here), or better, the related eigenvarieties. They can be 
quickly defined as follows, from the notion of {\it refined modular points}. Assume $F=\Q$ for simplicity 
and fix some prime $v$ of $E$ above $p$, so $E_v=\Q_p$ by assumption. For each modular point 
$x \in \cX(\rhob)$, it is known that $\rho_{x,p}:={\rho_x}_{|\G_{E_v}}$ is crystalline 
with distinct Hodge-Tate weights, say $k_1(x)<k_2(x)<k_3(x)$. 
Define a refined modular point as a pair $(x,(\widetilde{\varphi}_i(x))) \in \cX(\rhob) \times \mathbb{G}_m^3$ such 
that $x$ is a modular point and such that
$$(p^{k_{1}(x)}\widetilde{\varphi}_{1}(x),\,p^{k_{2}(x)}\widetilde{\varphi}_{2}(x),p^{k_{3}(x)}\widetilde{\varphi}_{3}(x))$$ is an ordering 
of the eigenvalues of the crystalline Frobenius of $\Dc(\rho_{x,p})$; there are up to $6$ ways to refine 
a given modular point. We define the $\rhob$-eigenvariety 
$$\cE(\rhob) \subset \cX(\rhob) \times \mathbb{G}_m^{3}$$
as the Zariski-closure of the refined modular points\footnote{Let us warn the reader that it is actually 
not quite the right definition (for instance Theorem D below may not hold with this one), although we 
shall content ourselves with it in this introduction and refer to Def.\ref{defeigen} for the right one.  }; 
the {\it complete infinite fern of type $U(3,E/\Q)$} is the subset $\mathcal{F}(\rhob) \subset \cX(\rhob)$ 
image of $\cE(\rhob)$ under the first projection. By a former result of the author, this eigenvariety 
turns out to be equidimensional of dimension $3$, and has some additional properties. The analogues of Coleman's arcs through a modular point 
$x \in \cX(\rhob)$ are now $3$-dimensional locally closed subspaces $C_{x,\varphi}$ (the
``leaves of the fern 
at $x$'') indexed by each refinements $(x,\varphi)$ of $x$. The germ of $C_{x,\varphi}$ at $x$ is canonical, 
and the modular points are actually Zariski-dense in $C_{x,\varphi}$.  The {\it $\rhob$-infinite fern of type 
$U(3,E/\Q)$} is the union of all the leaves $C_{x,\varphi}$ constructed this way when $(x,\varphi)$ runs over 
all the refined modular points: see below for the picture so far. 

\begin{figure}
\centering
\includegraphics{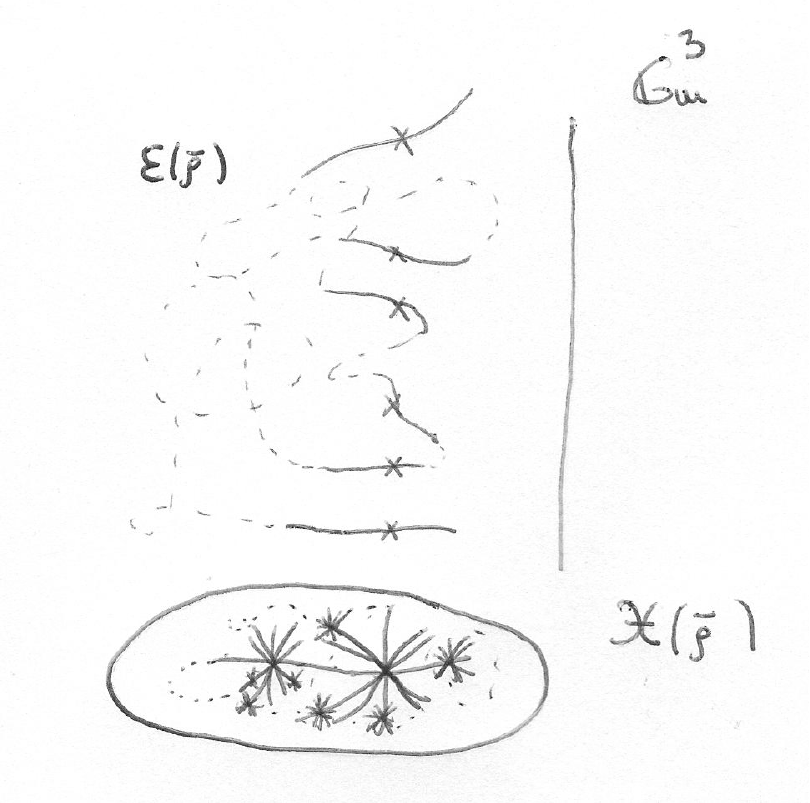}
\end{figure}

For dimension reasons, this situation is subtler than the one of
Gouv\^ea-Mazur. Indeed, assume that $\rhob$ is 
unobstructed to simplify, so $\cX(\rhob)$ is the open unit ball of dimension $6$. 
Remark that there is a one dimensional family of twists by Galois characters of type 
$\U(1)$, and $\cE(\rhob)$ is stable by this family of twists as well, so we should rather think
of $\cX(\rhob)$ as a $5$-dimensional space in which evolves the
$2$-dimensional fern $\Ref(\rhob)$ (with a Zariski-dense accumulating subset of
six times multiple points). The problem is that even if we knew that any two of the six leaves are 
transversal at each modular point, this would not be enough to exclude the possibility of a $4$-dimensional 
Zariski-closure. The situation is even worse when $[F:\Q]>1$ as in this case we have to pass from dimension 
$3[F:\Q]$ to $6[F:\Q]$. 

The idea of our proof is to study the relative positions of the tangent spaces of the local leaves $C_{x,\varphi}$, when 
$\varphi$ varies, at a given modular point $x \in \cX(\rhob)$. A key
intermediate result is the following:

{\bf Theorem B:} {\it There exist modular points $x \in \cX(\rhob)$ such that $\rho_{x,p}$ is absolutely irreducible and has distinct crystalline Frobenius eigenvalues in $k(x)$. If $x$ is such a point,  then the image of the natural map on tangent spaces $$\bigoplus_{y \mapsto x, y \in \mathcal{E}(\rhob)} T_y(\mathcal{E}(\rhob)) \longrightarrow T_x(\cX(\rhob))$$

has dimension $\geq 6$.}
\par\medskip\noindent

It is not difficult to show that Theorem $B$ actually implies Theorem $A$
(see~\S\ref{stratpar}). The first part of Theorem B is a simple application of eigenvarieties, 
but its second part is rather deep. It relies on two general results of independent interest whose 
proofs form the technical heart of this paper. The idea is to study the
image of the map in Theorem B in the tangent space of the deformation space of the local representation 
$\rho_{x,p}$, on which we will have a sufficiently efficient control as
we shall see. As an indication of this, recall that since the works of Kisin, Colmez and Bella\"iche-Chenevier, it is known that the 
restriction at $\G_{\Q_p}$ of the family of Galois representations over $\mathcal{E}(\rhob)$  has some 
very specific properties: they are {\it trianguline} in the neighborhood of the
``non-critical and regular''
modular points.

The first important result we prove is a purely local theorem about the deformation space of a given 
crystalline representation of $G_{\Q_p}={\Gal(\Qpb/\Qp)}$. Let $L$ be a finite extension of $\Q_p$ 
and let $V$ be a crystalline representation of $G_{\Q_p}$ of any $L$-dimension $n$. Assume 
that $\End_{\G_{\Q_p}}(V)=L$, that $V$ has distinct Hodge-Tate numbers, and that the eigenvalues $\varphi_i$ 
of the crystalline Frobenius on $D_{\rm crys}(V)$ belong to $L$ and satisfy 
$\varphi_i\varphi_j^{-1} \neq 1, p$ for all $i\neq j$. Let $\cX_V$ be the deformation functor of $V$ to 
the category of local artinian $L$-algebras with residue field $L$. It is pro-representable and formally 
smooth of dimension $n^2+1$. For each of the $n!$ orderings $\mathcal{F}$ of the $\varphi_i$ (such an ordering is 
called again called a {\it refinement} of $V$, for some obvious reasons), we defined in \cite[\S 2]{bch} the 
$\mathcal{F}$-trianguline deformation subfunctor $\cX_{V,\mathcal{F}} \subset \cX_V$, whose dimension 
is $\frac{n(n+1)}{2}+1$. Roughly speaking, the choice of $\mathcal{F}$ corresponds to a choice of a 
triangulation 
of the $(\varphi,\Gamma)$-module of $V$ over the Robba ring, and $\cX_{V,\mathcal{F}}$ parameterizes 
the deformations such that this triangulation lifts. When the $\varphi$-stable complete flag of 
$D_{\rm cris}(V)$ defined by $\mathcal{F}$ is in general position compared to the Hodge filtration, we 
say that $\mathcal{F}$ is {\it non-critical}. 

\par\medskip

{\bf Theorem C:} {\it Assume that  $n$ ``well-chosen'' refinements of $V$ are non-critical (e.g. all of them), 
or that $n\leq 3$. Then on tangent spaces we have an equality $$\cX_V(L[\varepsilon]) = \sum_{\mathcal F} 
\cX_{V,\mathcal{F}}(L[\varepsilon]).$$}

In other words, {\it any first order deformation of a generic crystalline representation is a 
linear combination of trianguline deformations}. Maybe surprisingly, our proof of this result is
by induction on the dimension of $V$, which requires first to extend the
statement to the world of non-\'etale $\fg$-modules over the 
Robba ring and work in this general setting (as in~\cite{bch}). We also have to extend and complete 
several results of~\cite{bch} to general paraboline deformation functors, and even in the trianguline case 
the proofs given here are actually slightly different from the ones there. Following the ideas of this paper, 
Theorem C has another purely local application to the Zariski-density of the crystalline points in some 
components of the representation space of $\G_{\Q_p}$ (for instance in the residually irreducible ones), that will however be given elsewhere. Last but not least, let us 
remark that the codimension of the crystalline locus in $\cX_V$, namely $\frac{n(n+1)}{2}$, actually 
coincides with the (conjectural) dimension of $\cX(\rhob)$ if $\rhob$ is of type $\U(n)$
(see below). This numerical 
coincidence reminds of course the ones discussed in \cite{cht} and is fundamental to the application to the 
Zariski-density of modular points (e.g. in the unobstructed case of Thm. A).
\smallskip

The second key ingredient for the proof of Theorem $B$ is an important theorem about eigenvarieties that we simply state 
as a slogan in this introduction:

{\bf Theorem D:} {\it Eigenvarieties are \'etale over the weight space at non-critically refined, regular, 
classical points.}

See~Theorems~\ref{thmetalenoncrit} and~\ref{thmetalenoncrit2} for precise statements concerning unitary 
eigenvarieties. In terms of families of $p$-adic automorphic forms, this theorem
means that the natural 
$p$-adic family passing through a non-critically refined (and regular) automorphic representation may be 
parameterized by a polydisc with weight maps for parameters (either automorphic, or Hodge-Tate-Sen weights), 
which is maybe the most natural statement we may expect in the theory of $p$-adic families. For ordinary 
refinements several instances of this slogan have been proved by Hida. Our proof relies
in particular on a number of properties of automorphic representations for unitary groups (including some multiplicity one results of 
Rogawski for $U(3)$ and of Labesse for $U(n)$), on the generalized theta-maps studied by Jones, 
and on a result of Kisin.

Let us go back to Theorem A. We have dealt with cases of absolutely irreducible $\rhob$ to fix ideas, but 
it holds in any residually reducible case, with the same proof, if we work
with pseudocharacters instead of representations (using \cite{chdet} for $p=3$ !). More interestingly, we 
expect that our approach will lead to a proof of the Zariski-density of the infinite fern of type $U(n)$ 
for any $n$. In this case $\cX(\rhob)$ has conjectural dimension $[F:\Q]\frac{n(n+1)}{2}$, $\cE(\rhob)$ 
has dimension $n[F:\Q]$, and there are generically
$n!^{[F:\Q]}$ local leaves in $\cF(\rhob)$ through each modular point. (When $n=3$, the coincidence 
$3!=3(3+1)/2=6$ actually plays no specific role.) Although Theorem $C$ and $D$ apply in this generality, 
we were faced to a new difficulty concerning the existence of sufficiently many global automorphic 
representations which are {\it generic enough} at $p$: we refer to~\S\ref{genoset} for a discussion about it. 
On the other hand, the degree of the base field is rather insensitive to our method, which allowed us for 
instance to extend the result of Gouv\^ea-Mazur to the Hilbert modular case (see~Theorem~\ref{cashbeigen}):

{\bf Theorem E:} {\it Let $F$ be a totally real field of even degree, $p$ an odd prime totally split in $E$, 
$\got{d}$ the Leopold defect of $F$ at $p$, and let $\rhob: G_{F,S} \longrightarrow \GL_2(\F_q)$ an absolutely 
irreducible modular Galois representation.

 Then the irreducible components of the Zariski-closure of the essentially modular points in $\cX(\rhob)$ all 
have dimension at least $1+{\got{d}}+2[F:\Q]$. If $H^2(G_{F,S},{\rm ad}(\rhob))=0$, then the modular points are 
Zariski-dense in $\cX(\rhob)$.}

We end this introduction by stating an interesting general result about {\rm Adjoint'} Selmer 
groups\footnote{`Pronounce ``adjoint primed Selmer groups''.} that follows from the theory of 
eigenvarieties and Theorems $C$ and $D$ above. Let $n\geq 1$ be any integer, let $E/F$ be a CM field, 
and let $\Pi$ be a cohomological cuspidal automorphic representation of $\GL_n(\AAA_E)$ such that 
$\Pi^\vee \simeq \Pi^c$. Assume that $p$ is totally split in $E$, that $\Pi$ is unramified at each 
place of $E$ above $p$, as well as a mild (and conjecturally empty) condition when $n$ is even, 
see~\S\ref{adselmer}. Let $$\rho_\Pi : \G_E \longrightarrow \GL_n(\Qpb)$$ be a $p$-adic Galois 
representation attached to $\Pi$. It turns out that the adjoint representation ${\rm ad}(\rho_\Pi)$ 
has a natural extension to $G_F$ that we denote by $\adr(\rho_\Pi)$, and whose $H^1$ has a natural 
interpretation as the tangent space of the deformation functor of $\rho_\Pi$ of type $U(n)$. This 
interpretation allows to define subspaces 
$H^1_{\Ref}(F,\adr(\rho_\Pi))$ for any choice of refinements of the ${\rho_\Pi}_{|\G_v}$ 
for each place $v$ of $F$ dividing $p$. In the following statement, we use standard notations from the 
theory of Selmer groups: $f$ is the Bloch-Kato condition (automatic outside $p$ for $\adr(\rho_\Pi)$) 
and $H^1_s=H^1/H^1_f$. See also~\S\ref{examplelocalmain} for the precise definition of {\it weakly generic 
and regular}, let us simply say that a crystalline representation of $\G_{\Q_p}$ with sufficiently generic 
crystalline Frobenius eigenvalues and whose refinements are all non-critical has this property (in particular 
it is absolutely indecomposable).

{\bf Theorem F:} {\it Assume that ${\rho_\Pi}_{|\G_v}$ is weakly generic and regular for each $v|p$. Then 
$$H^1(F,\adr(\rho_\Pi)) \longrightarrow \prod_{v|p} H^1_s(F_v,\adr(\rho_\Pi))$$
is surjective. Moreover, $H^1(F,\adr(\rho_\Pi))=\sum_{\Ref} H^1_{\Ref}(F,{\rm \adr(\rho_\Pi)})$ and for 
$S$ big enough we have $\dim H^2(G_{F,S},\adr(\rho_\Pi))=\dim H^1_f(F,\adr(\rho_\Pi))$.
}

The second part of this theorem is a variant for Selmer groups of Theorem C. It is the infinitesimal 
counterpart of the infinite fern of type $U(n,E/F)$ in the Adjoint' Selmer group.  It says in particular that $H^1(F,\adr(\rho_\Pi))$ 
is spanned by elements which are not too far from being geometric. As far as we known, this is the first 
example of such statement. The first part of the theorem is actually deduced from the second one, and the 
assertion on $H^2$ follows from the computation of the sign of $\rho_\Pi$  in \cite{bchsign}.  We believe 
that a variant of Theorem F should hold for a huge class of Selmer groups. Let us mention here that the 
proof of Theorem F also uses some base-change trick to reduce to a situation where the multiplicity one 
results of Labesse (and Theorem D) apply.

The theorem above has a corollary concerning the Galois representations associated to essentially-selfdual 
cuspidal, cohomological, automorphic representations of $\GL_n$ over a totally real field: 
see~Thm.\ref{thmselmertotreel}. Here is a special case. Let $\Pi$ be a selfdual, cohomological, and 
cuspidal, automorphic representation of $\GL_n(\AAA_\Q)$, and let $\rho_\Pi : \G_\Q \longrightarrow 
\GL_n(\Qpb)$ be an associated $p$-adic Galois representation, normalized so that $\rho_\Pi^* \simeq 
\rho_\Pi(n-1)$. Recall that $\rho_\Pi$ is symplectic if $n$ is even and orthogonal otherwise
by~\cite{bchsign}.  Assume again that $\Pi_p$ is unramified, as well as the same mild assumption as 
above for $n$ even. 

{\bf Theorem G:} {\it Assume that ${\rho_\Pi}_{|\G_{\Q_p}}$ is weakly generic and regular, and set 
$V=\Symm^2(\rho_\Pi)(n-1)$ if $n$ is even, $V=\Lambda^2(\rho_\Pi)(n-1)$ otherwise. Then the restriction
map $H^1(\Q,V) \longrightarrow H^1_s(\Q_p,V)$ is surjective. }

\ps \medskip
 Let $f=q+a_2q^2+\dots$ be a classical modular eigenform of weight $k>1$, level prime to $p$, and let $\rho_f : \G_{\Q,S} \rightarrow \GL_2(\Qpb)$ be an associated $p$-adic Galois representation.

{\bf Theorem H:} {\it Assume that ${\rho_f}_{|\G_{\Q_p}}$ is absolutely indecomposable and that the ratio of 
the two roots of the Hecke polynomial $X^2-a_pX+\varepsilon(p)p^{k-1}$ of
$f$ at $p$ is not a root of unity. 
Set $V_n=\Symm^n(\rho_f)\det(\rho_f)^{-n/2}$. 

For $n=2$ and $n=6$, the restriction $H^1(\Q,V_n) \longrightarrow  H^1_s(\Q_p,V_n)$ is surjective 
and $\dim H^1_f(\Q,V_n)=\dim H^2(G_{\Q,S},V_n)$. }

As the reader may have noticed, all the global statements of this paper concern number 
fields which are totally split at $p$. It is clear from our method that this assumption could be removed 
but we simply did not have the energy to do so. One reason for this is our reliance on the local results 
of \cite{kisin} and \cite{bch}, which are only written at the moment for the base field $\Q_p$. Another 
slightly annoying assumption we have is the 
regularity condition in Theorems D and F. Again, we expect that it could be removed if the technical 
problems discussed in~\cite[\S 4]{bch} and occuring in studying 
``refined families'' were solved. For instance, these two improvements
together, combined with the recent developpements on the Sato-Tate
conjecture, would lead to an extension of Theorem H to any $n \equiv 2 \bmod 4$.

\medskip

{\sc Aknowledgements:} We are grateful to Gebhard B\"ockle for some stimulating discussions at 
the origin of this paper. Part of our proof relies on my book with Jo\"el Bella\"iche, it is a 
pleasure to thank him here for the many fruitful exchanges we had the past
years. We thank also Laurent Berger 
and Kevin Buzzard for some useful discussions, as well as Owen Jones for having made his preprint 
available. Some results of this paper have been exposed by the author in a Cours Peccot at the 
Coll\`ege de France in March-April 2008 \cite{peccot}, it is a pleasure to thank this institution as well as the Fondation 
Peccot for their support.

\section{The universal Galois representation of type $\U(n)$}
\newcommand{\trb}{\widetilde{\rho}}
\subsection{Galois representations of type $\U(n)$ with $n$ odd and their deformations} \label{prelimdefo}
Let $E$ be a CM field, $F$ its maximal totally real subfield and
fix $c \in \G_F$ a complex conjugation. Let $p$ be an odd prime number. Let $S$ be a finite set of primes of 
$E$ containing the primes ramified
over $F$ and the primes above $p$. We assume that $c(S)=S$ and we shall sometimes view $S$ 
as a set of places of $F$ containing the subset $S_\infty$ of all the archimedean primes. 
We fix an odd integer $n$, eventually we shall take $n=3$. \ps

Let $A$ be a commutative ring and let $\rho : \G_{E,S} \longrightarrow
\GL_n(A)$ be an absolutely irreducible representation, in the sense that 
$A[\rho(\G_{E,S})]=M_n(A)$. We say that $\rho$
is {\it of type $\U(n)$} if $$\rho^\ast \simeq \rho^c$$ as $A[\G_{E,S}]$-modules,
or which is the same if ${\rm trace}\rho(g^{-1})={\rm trace}\rho(cgc)$ for all $g \in \G_E$. 
Let $\Gal(E/F)$ act on $\GL_n$ by $c(g):={}^tg^{-1}$ and view $\mathcal{G}_n:=\GL_n
\rtimes \Gal(E/F)$ as a group scheme over $\Z$. The terminology {\it $\rho$ is of type
$\U(n)$} comes from the fact that $\mathcal{G}_n$ may be viewed as
the reduced Langlands dual group of the unitary groups in $n$ variables associated to
$E/F$ and the following lemma, which is a variant of \cite[\S 1]{cht}. \ps

\begin{lemma}\label{deftrb} $\rho$ is of type $\U(n)$ if, and only if, there
exists a representation $$\trb : \G_{F,S} \longrightarrow \GL_n(A)
\rtimes \Gal(E/F)$$
where the induced map $\G_{F,S} \rightarrow \Gal(E/F)$ is the natural
map and $\trb_{|\G_{E,S}}=\rho$. \end{lemma}

Indeed, $\trb(c):=Pc$, with $P \in \GL_n(A)$, defines an extension of
$\rho$ as in the statement if, and only if, ${}^t\!P=P$ and 
${}^t\rho(g)^{-1}=P^{-1}\rho(cgc)P$ for all $g \in \G_{E,S}$. This last
condition actually implies ${}^t\!P=\mu P$ for some $\mu \in A^*$ as $\rho$ is absolutely irreducible, and $\mu=1$ as 
$n$ is odd. Note that the extension $\trb$ is not unique, but $\trb'$ is another one if, and only if, 
$\trb'(c)=P'c$ with $P'=\lambda P$ and $\lambda \in A^*$.
\ps\newcommand{\cD}{\mathcal D}
Let $q$ be a power of $p$ and let us fix from now on $$\rhob : \G_{E,S} \longrightarrow
\GL_n(\F_q)$$
an absolutely irreducible, continuous, representation of type $\U(n)$. 
Let $\Z_q$ be the Witt vectors over $\F_q$ and let  $\mathcal{C}$ be the category
of artinian local $\Z_q$-algebras with residue field $A/m_A=\F_q$. Consider the
deformation functor $$\cD : \mathcal{C} \rightarrow {\rm Sets}$$
defined as follows : for any object $A$ of $\mathcal{C}$, $\cD(A)$ is the set of $A$-isomorphism classes of 
continuous representations $\rho_A : \G_{E,S} \longrightarrow \GL_n(A)$ of type $\U(n)$ such that $\rho_A \otimes_A \F_q$ 
is isomorphic to $\rhob$. As $\rhob$ is absolutely irreducible, this functor is pro-representable by a
complete local noetherian $\Z_q$-algebra $R(\rhob)$ with residue field
$\F_q$. Indeed, it is the quotient of Mazur's universal $\G_{E,S}$-deformation ring $R'(\rhob)$ by the ideal generated by  
the elements ${\rm trace}\,\,\rho^u(g^{-1})-{\rm trace}\,\,\rho^u(cgc)$, for all $g \in \G_{E,S}$, where $\rho^u$ is Mazur's universal $R'(\rhob)$-valued deformation. (Note that the argument so far did not use the assumption $p$ odd.)\ps

Using class field theory and arguments of Mazur we will see below that the Krull dimension of 
$R(\rhob)/(p)$ is always $\geq [F:\Q]\frac{n(n+1)}{2}$. The precise structure of $R(\rhob)$ is presumably 
extremely complicated in general, however 
we shall not be interested in these kind of complications in this paper.\ps

 Let $\mathfrak{g}_n \simeq M_n(\F_q)$ be the Lie algebra of $\mathcal{G}_n$ 
over $\F_q$ viewed as a representation of $\mathcal{G}_n(\F_q)$. Fix once and for all some $\trb$ lifting $\rhob$ as in Lemma~\ref{deftrb}. Thanks to 
this choice we may view $\mathfrak{g}_n$ as a $\G_{F,S}$-module, that we shall denote by $\arb$. This latter module does not depend on the choice of $\trb$. 
By definition, $\arb_{|\G_E}$ is the usual adjoint representation of $\rhob$. Moreover, if 
$\trb(c)=Pc$, then $c$ acts on $\arb=M_n(\F_q)$ 
as\footnote{As an exercise, the reader can check that up to conjugating $\rhob$ if necessary, we may choose a $\trb$ such that $\trb(c)=c$, 
i.e. $P={\rm id}$ (use that $n$ and $p$ are odd).} $X \mapsto -P{}^t\!XP^{-1}$. In 
particular, $H^0(\G_{F,S},\arb)=0$ as $\rhob$ is absolutely irreducible and $p$ is odd. 

\begin{definition}\label{defunobs} We say that the $\U(n)$-deformation theory of $\rhob$ is {\it unobstructed}, or for short that $\rhob$ is unobstructed, if  
$H^2(\G_{F,S},\arb)=0$.
\end{definition}
By Tate's global and local duality theorems, $\rhob$ is unobstructed if, and only if, $$\Sha^1_S(\G_{F,S},\arb(1))=0 \mbox{\, \,and\,\, } 
H^0(F_v,\arb(1))=0 \mbox{\, for all finite prime $v \in S$}$$ 
(note that $\mathfrak{g}_n$ is a selfdual $\F_q[\mathcal{G}_n(\F_q)]$-module and 
see~\cite[Rem. 5.2.(c)]{milne}). 

\begin{prop}\label{propunobs} If $\rhob$ is unobstructed, then $R(\rhob)$ is formally smooth over $\Z_q$ of relative dimension 
$$[F:\Q]\frac{n(n+1)}{2}.$$
\end{prop}

Set $V=\arb$. Tate's global Euler characteristic formula shows  that 
$$\dim  H^1(\G_{F,S},V)=[F:\Q]\dim V^{\trb(c)=-1}+\dim H^0(\G_{F,S},V)+\dim H^2(\G_{F,S},V)$$
(dimensions are over $\F_q$ and recall that $p$ is odd). But $X \mapsto XP$ identifies $M_n(\F_q)^{\trb(c)=-1}$ with the subspace of symmetric matrices in $M_n(\F_q)$, which has dimension $\frac{n(n+1)}{2}$ , 
so $\dim H^1(\G_{F,S},\arb)=[F:\Q]\frac{n(n+1)}{2}$ if $\rhob$ is unobstructed. \par

The proposition follows from this computation and from the following general lemma. If $A$ is an object of $\mathcal{C}$, two group homomorphisms with 
target $\mathcal{G}_n(A)$ are said {\it equivalent} if they are conjugate by an element in $\Ker(\mathcal{G}_n(A) \rightarrow \mathcal{G}_n(\F_q))$. 
Define another functor $\cD' : \mathcal{C} \rightarrow {\rm Sets}$ as follows : for an object $A$ of $\mathcal{C}$, $\cD'(A)$ is 
the set of equivalence classes of continuous homomorphisms $r_A : \G_{F,S} \rightarrow \mathcal{G}_n(A)$ such that 
$r_A \otimes_A \F_q = \trb$ (which has been fixed above). 

\begin{lemma}\label{calctanun} The map $r_A \mapsto {r_A}_{|\G_E}$ induces an isomorphism $\cD' \isomo \cD$. Moreover, 
$\cD(\F_q[\varepsilon]) =H^1(\G_{F,S},\arb)$ and if $H^2(\G_{F,S},\arb)=0$, then $\cD$ is formally smooth over $\Z_q$.
\end{lemma}

\begin{pf} For each $A$, $\cD'(A) \rightarrow \cD(A)$ is surjective by the proof of Lemma~\ref{deftrb}. Let $r_1, r_2 \in \cD'(A)$
be such that ${r_1}_{|\G_E} \simeq {r_2}_{|\G_E}$, up to equivalence we may assume that ${r_1}_{|\G_E}={r_2}_{|\G_E}$. 
Set $r_i(c)=P_ic$, where $P_i \in \GL_n(A)$. Then both $P_1$ and $P_2$ intertwine the conjugate and dual of 
${r_1}_{|\G_E}$, so $P_1=\lambda P_2$ for some $\lambda \in A^*$. As each $r_i$ lifts $\trb$, $\lambda \in 1 + m_A$. 
As $p$ is odd, $\lambda=\mu^2$ for some $\mu$ in $1+m_A$, so $\mu^{-1} r_1 \mu$ is equivalent to $r_1$ and coincides with $r_2$. This
concludes the proof of the first part of the statement, the other assertions follow from this and standard facts from group cohomology.
\end{pf}

%
%\begin{remark} {\rm A similar study can be done with $n$ even. In this case, an absolutely irreducible 
%Galois representation $\rho : G_E \rightarrow \GL_n(A)$ is said of type $\U(n)$ if $\rho^* \simeq \rho^c (1)$ {\it and} 
%if the matrices $P \in \GL_n(A)$ such that ${}^t\!\rho(g)^{-1}=P\rho(cgc^{-1})\omega_{\rm cycl}(g)P^{-1}$ for all $g \in \G_E$ 
%are symmetric. In this case the group $\mathcal{G}_n$ has to be modified but Prop.~\ref{propunobs} still holds.}
%\end{remark}

In the sequel, we will only be interested in the rigid analytic space $$\mathfrak{X}(\rhob)={\rm Sp}(R(\rhob)[1/p])$$ 
over $\Q_q:=\Z_q[1/p]$ associated to the generic fiber of $R(\rhob)$, in the sense of Berthelot (see e.g. \cite[\S 7]{dejong} for the basics on this construction). For example, in the unobstructed case $\mathfrak{X}(\rhob)$ is simply the open unit ball over 
$\Q_q$ of dimension $[F:\Q]\frac{n(n+1)}{2}$. Before stating the universal property of $\mathfrak{X}(\rhob)$ we need to set some notations and review 
some facts concerning families of representations.\ps
Let $Y$ be an affinoid over $\Q_q$ and $\rho : G_{E,S} \rightarrow \GL_n(\OO(Y))$ a continuous representation. For 
$y \in Y$ a closed point, with residue field $k(y)$, we denote by $\rho_y :  G_{E,S} \rightarrow \GL_n(k(y))$ the evaluation of $\rho$ at $y$ and by 
$\rhob_y :  G_{E,S} \rightarrow \GL_n(k_y)$ the semi-simplification of its residual representation. Here $k_y$ denotes the residue field of $k(y)$, it 
comes with a natural morphism $\F_q \rightarrow k_y$. We say that {\it $\rho$ is a lift of $\rhob$} if 
$\rhob_y \simeq \rhob \otimes_{\F_q} k_y$ for all $y \in Y$. \ps
Let $\cD^{\rm an}$ be the (contravariant) functor from the category of 
$\Q_q$-affinoids to {\rm Sets} where $\cD^{\rm an}(Y)$ is the set of $\OO(Y)$-isomorphism classes of continuous homomorphisms 
$\rho_Y : G_{E,S} \rightarrow \GL_n(\OO(Y))$ which are of type $\U(n)$ and lift $\rhob$. As $\rhob$ is absolutely irreducible, $\rho_Y \mapsto {\rm trace} \, \rho_Y$
induces a bijection between $\cD^{\rm an}(Y)$ and the set of continuous pseudocharacters $T_Y : G_{E,S} \rightarrow \OO(Y)$ lifting ${\rm trace}\,\rhob$ (for the obvious definition) 
and of type $\U(n)$, that is such that $T_Y(g^{-1})=T_Y(cgc)$ for all $g \in G_{E,S}$. The following fact is probably well-known, 
it follows for instance from~\cite[Thm. 3.16]{chdet}, but this is much simpler here as $\rhob$ is absolutely irreductible.

\begin{lemma} $\cX(\rhob)$ represents $\cD^{\rm an}$. \end{lemma}

This space $\cX(\rhob)$ is actually a connected component of the $\U(n)$-type locus of
the space $\cX_n$ mentionned in the introduction.

\subsection{Modular Galois representations of type $\U(n)$ and examples}\label{sectmod} We keep the assumptions of \S\ref{prelimdefo}. The
main source of Galois representations of type $\U(n)$ is the degree $n-1$ \'etale
cohomology of the quotients of the complex open unit ball $\sum_{i=1}^{n-1} |z_i|^2<1$ in $\C^{n-1}$ by the arithmetic congruence
subgroups of $\U(n-1,1)(\R)$ attached to $E/F$. When $n=3$ these are also called Picard modular surfaces. Thanks to the advances in the theory of automorphic 
forms we may describe them using $\GL(n)$.

\noindent Let $\Pi$ be a cuspidal automorphic representation of $\GL_n(\AAA_E)$
such that :  

(P1) $\Pi^\vee \simeq \Pi^c$,  

(P2) $\Pi_v$ is cohomological for each archimedean place $v$ of $E$, 

(P3) $\Pi_v$ is unramified whenever $v \notin S$ or $v$ divides $p$. \pn

\noindent Fix once and for all a pair of fields embeddings $\iota :
\overline{\Q} \rightarrow \overline{\Q}_p$ and $\, \overline{\Q} \rightarrow
\C$. Class-field theory for $n=1$, work of Rogawski for $n=3$, and recent works of
Shin~\cite{shin} and of the participants of the book project of the Paris GRFA 
seminar for any $n$~\cite{GRFAbook}, attach to such a $\Pi$ 
and $\iota$ a continuous semi-simple
representation $$\rho_\Pi : \G_{E,S} \rightarrow \GL_n(\Qpb)$$
of type $\U(n)$ which is compatible to the Frobenius semi-simplified local Langlands
correspondance at all primes $v$ not dividing $p$ and which is cristalline at each prime $v$ of
$E$ above $p$. In particular, this representation is geometric in the sense of Fontaine-Mazur. Actually, 
infinitely many suitable real quadratic base changes of it are even geometric in the sense of the introduction by construction 
(all when $[F:\Q]>1$). It is known that $\rho_\Pi$ is irreducible when $n=3$, and conjectured in general. Fix a ring 
homomorphism $\overline{\Z}_p \rightarrow \Fpb$. Denote by $\rhob_\Pi : \G_{E,S} \rightarrow
\GL_n(\Fpb)$ the
isomorphism class of the semi-simplification of $\rho_\Pi$. 

\begin{definition}\label{defmodular} A representation $\rho : G_{E,S} \rightarrow \GL_n(\Qpb)$ is {\rm modular of type $\U(n)$} if there exists a 
$\Pi$ satisfying $(P1)$ to $(P3)$ as above such that $\rho \simeq \rho_\Pi$. A representation 
$\rhob : G_{E,S} \rightarrow \GL_n(\Fpb)$ is {\rm modular of type $\U(n)$} if there exists a
$\Pi$ as above such that $\rhob \simeq \rhob_\Pi$.
\end{definition}

For short, we shall often say {\it modular} for {\it modular of type $\U(n)$}. It is tempting to believe that a variant of Serre's conjecture holds in this
context : 
\begin{center} {\it Is any continuous absolutely irreducible $\rhob$ of type $\U(n)$ modular} ? \end{center} 
When $n=3$, by Rogawski's construction, any modular $\rhob$ occurs in $H^2_{et}(S_{\overline{E}},\Fpb(1))$
where $S$ is a suitable Picard modular surface over $E$. A delicate aspect of this variant of Serre's conjecture is that such an $H^2$ might contain some 
classes\footnote{I do not know any example of such a class that generates an irreducible Galois representation of dimension $3$.} that do not lift to cohomology classes in characteristic $0$: this happens for instance for Mumford's fake projective planes.
Nevertheless, the recent advances on $R=T$ theorems and potential
modularity in this context (e.g. by Clozel, Harris and Taylor and their
co-workers) suggest a positive answer.

\begin{example}\label{exmodf}{\rm  If $f$ is a classical modular eigenform then ${\rm
Symm^2}(\rho_f)\otimes \det(\rho_f)^{-1}_{|\G_{E,S}}$ is of type $\U(3)$. When $F/\Q$ is Galois and solvable, it is modular by
results of Arthur-Clozel (if $F\neq \Q$), Gelbart-Jacquet and Rogawski. By the modularity theorem, this
applies in particular to ${\rm Symm^2}A[p]^*(1)$ for any elliptic curve $A$ over $\Q$.} \end{example}

\subsection{Statement of the theorem}

Let $\rhob$ be a Galois representation of type $\U(n)$ and $\cX(\rhob)$ the generic fiber of its universal deformation ring of type $\U(n)$. A closed point $x \in \cX(\rhob)$ is called {\it modular} if the
Galois representation $\rho_x : G_{E,S} \rightarrow \GL_n(k(x))$ attached to
$x$ is modular. Let $$\cXm \subset \cX(\rhob)$$ 
be the (countable) subset of modular points. By definition, $\rhob$ is
modular if, and only if, $\cXm \neq \emptyset$. We can now state the main result of this paper.

\begin{thm}\label{mainthm} Assume that $\rhob$ is modular of type $\U(3)$, unobstructed, and that $p$
splits in $E$. Then $\cXm$ is Zariski-dense in $\cX(\rhob)$.
\end{thm}

A subset of closed points $Z$ of a rigid space $Y$ is called {\it Zariski-dense} if the only reduced closed subspace of $Y$ containing $Z$ is $Y_{\rm red}$.
When $Y$ quasi-Stein (which is the case of the $\cX(\rhob)$), it is equivalent to ask that any global function $f \in \OO(Y)$ vanishing at all points of $Z$ is locally nilpotent. In the
context of the theorem, $\cX(\rhob)$ is the open unit ball of dimension
$d=[F:\Q]\frac{n(n+1)}{2}$ over $\Q_q$, so the statement simply means that any power series $f \in
\C_p[[t_1,\dots,t_d]]$ converging on the whole open ball $|t_i|<1$ and vanishing on $\cXm$ is identically $0$. \ps

\begin{example}\label{ex}{\rm Set $E=\Q(i)$, $p=5$, $S=\{\infty,2,5,17\}$ and let $A$ be the elliptic curve 
$y^2+xy+y=x^3-x^2-x$ over $\Q$,
which is isogenous to the Jacobian of $X_0(17)$. Then $\rhob=\Symm^2(A[5])(-1)_{|\G_E}$ is
modular of type $\U(3)$ and unobstructed (see the appendix). }
\end{example}  

We do not know if $\cXm$ is dense in $\cX$ for the $p$-adic topology. 
However, we have the following positive corollary : 

\begin{cor} (same assumptions) For any object $A$ of $\mathcal{C}$ and any lift $\rho_A \in \cD(A)$ 
of $\rhob$, there is a finite extension $L$ of $\Q_q$ and a set $\rho_1$, $\rho_2$, \dots, $\rho_r$ of
modular Galois representations $\G_{E,S} \rightarrow \GL_3(\OO_L)$ of type $\U(3)$
such that the $\Z_q[\G_{E,S}]$-module $\rho_A$ is a subquotient of $\oplus_{i=1}^r \rho_i$.
\end{cor}

Indeed, this follows from the theorem and the following general fact due to Chevalley: if $A$ is a complete
local noetherian ring, and if $\{a_i, \,i \in I\}$ is a family 
of proper, cofinite lenght, ideals of $A$ such that the intersection of all 
the $a_i$ is zero, then the topology defined by the $a_i$ is the adic topology of
the maximal ideal of $A$.\ps

We will actually prove the following result, of which Theorem~\ref{mainthm} follows at once by Prop.\ref{propunobs}.

\begin{thm}\label{mainthmobs} Assume that $p$ splits in $E$. Then the irreducible components of the Zariski-closure of $\cXm$ inside $\cX(\rhob)$  all have dimension at least  $6[F:\Q]$.  
\end{thm}

We are led to the following optimistic conjecture (for any integer $n$).

\begin{conj}\label{conjgen} Let $\rhob$ be any Galois representation of type $\U(n)$, then
$\cXm$ is Zariski-dense in $\cX(\rhob)$.
\end{conj}

\section{The infinite fern of Galois representations of type $\U(3)$} 

The main ingredient in the proof of Thm.~\ref{mainthmobs} is the so-called infinite fern inside $\cX(\rhob)$ and its study. 
We use the notations of~\S\ref{prelimdefo} and {\it we assume furthermore that $\rhob$ is a modular Galois representation of type $\U(3)$ and that the odd prime 
$p$ splits in $E$}. We denote by $S_p$ the set of places of $F$ dividing
$v$, $S_0=S\backslash (S_\infty \cup S_p)$, and we fix once and for all, for each $v \in S_p$, a place $\tilde{v}$ 
of $E$ above $p$, so $\Q_p=F_v=E_{\tilde{v}}$. Set $I=\{1, 2,3\}$.

\subsection{The eigenvariety and the complete infinite fern} Let $\WW$ and $\WW_0$ be the rigid
analytic spaces over $\Q_p$ parameterizing respectively the $p$-adic continuous characters
of $\Q_p^*$ and $\Z_p^*$, $\WW_0$ is a finite disjoint union of $1$-dimensional open unit
balls over $\Q_p$ and $\WW \simeq \mathbb{G}_m \times \WW_0$. Consider the subset 
$$\cZ \subset \cX(\rhob)(\Qpb) \times \WW(\Qpb)^{I \times S_p}$$
of pairs $(\rho_\Pi,\delta)$ satisfying (i) and (ii) below : 

(i) $\rho_\Pi$ is the modular Galois representation of type $\U(3)$ associated to $\Pi$.

\noindent In particular $\rho_{\Pi,v}:=(\rho_\Pi)_{|\G_{E_{\tilde{v}}}}$ is
crystalline with distinct Hodge-Tate weights for each $v \in S_p$. We denote by
$k_{1,v} <
k_{2,v} < k_{3,v}$ these weights by increasing order, with the convention
that the cyclotomic character has Hodge-Tate weight $-1$. 

(ii) $\delta=(\delta_{i,v})$ where for
each $(i,v) \in I \times S_p$, the character $\delta_{i,v} : \Q_p^* \rightarrow \Qpb^*$ is
the product of $x \mapsto x^{-k_{i,v}}$ and of the character sending
$\Z_p^*$ to $1$ and $p$ to $\varphi_{i,v}$, where $(\varphi_{i,v})_{i \in I}$ is an ordering of the eigenvalues of the
crystalline Frobenius of $\Dc(\rho_{\Pi,v})$ (recall that $E_{\tilde{v}}=\Q_p$).

\noindent Each such pair is called a {\it refined modular point}. There are
up to $6^{[F:\Q]}$ ways to refine a given modular point, corresponding to the
number of ways to choose an ordering of the eigenvalues of the crystalline
Frobenius of $\Dc(\rho_{\Pi,v})$ for each $v \in S_p$.

\begin{definition}\label{defeigen} The {\rm eigenvariety} of type $\U(3)$ of $\rhob$ is the Zariski-closure $$\cE(\rhob)  \subset
\cX(\rhob) \times \WW^{I \times
S_p}$$ of the subset $\cZ$ inside $\cX(\rhob) \times \WW^{I \times S_p}$. The {\rm complete infinite
fern} of type $\U(3)$ of $\rhob$ $$\cF(\rhob) \subset \cX(\rhob)$$ 
is the set theoretic image of $\cE(\rhob)$ via the first projection
${\rm pr}_1: \cX(\rhob) \times {\WW}^{I \times S_p} \rightarrow \cX(\rhob)$.
\end{definition}

By definition $\cE(\rhob)$ is a reduced analytic space. The {\it weight space} is the space $\WW_0^{I \times S_p}$. The natural map $\WW \rightarrow \WW_0$ and the second projection induce a natural map $$\kappa : \cE(\rhob) \rightarrow \WW_0^{I \times S_p}.$$ This map turns out to refine the Hodge-Tate-Sen map. Indeed, recall that the work of Sen defines, for each $v \in S_p$, a monic polynomial $$P_{{\rm Sen},v}(t) \in \OO(\cX(\rhob))[t]$$ of 
degree $3$ whose evaluation at any $x \in \cX(\rhob)$ is the usual Sen polynomial of $\rho_{x,v}$ 
(the roots of which are the generalized Hodge-Tate weights of $\rho_{x,v}$). The coefficients of this collection of polynomials give rise to a natural map 
$\kappa_{{\rm HT}} : \cX(\rhob) \rightarrow \AAA^{I \times S_p}$. The morphism $\beta: \WW_0^{I \times S_p}  \rightarrow \AAA^{I \times S_p}$ given by the formula $\beta((\delta_{i,v}))_v=\prod_i (t+\left(\frac{\partial \delta_{i,v}}{\partial x}\right)_{|x=1})$ is the composite of a finite covering of $\AAA^{I \times S_p}$ of degree $6^{[F:\Q]}$, \'etale over the locus parametrizing polynomials with nonzero discriminant, with an \'etale morphism of infinite degree (essentially a $p$-adic logarithm). By definition, we have $\kappa_{HT} \, {\rm o}\,{\rm pr}_1 = \beta \, {\rm o}\, \kappa$ on $\cE(\rhob)$,

% By construction again, the composite of the natural map $\cE(\rhob) \rightarrow \Hom(\Z_p^*,\mathbb{G}_m)^{I \times S_p}$ 
%and of the derivative a $1$ map $\Hom(\Z_p^*,\mathbb{G}_m) \rightarrow \AAA^1$....
%
%The work of Sen \cite{sen} defines a {\it Sen polynomial map} $$P_{\rm Sen} : \cX(\rhob) \rightarrow \AAA^{I \times S_p}$$
%where $P_{\rm Sen}(x) $ is the Sen polynomial of the Galois re

\begin{thm}\label{thmeigen} $\cE(\rhob)$ is equidimensional of dimension $3[F:\Q]$. The map $\kappa$ is locally finite: $\cE(\rhob)$ is admissibly covered by the open affinoids $U$ such that $\kappa(U) \subset \WW_0^{I\times S_p}$ is open affinoid and $\kappa_{|U} : U \rightarrow \kappa(U)$ is finite.
\end{thm} 

This follows from the main theorem of the theory of $p$-families of automorphic forms for the 
definite unitary group $\U(3)$ (see \cite{che},\cite{che2}), and from the base
change results of Rogawski~\cite{roglivre}. This point of view on eigenvarieties is a generalization of the one of  Coleman and Mazur in~\cite{colemanmazur}. Concretely, the equations defining $\cE(\rhob)$ inside $\cX(\rhob) \times
\WW^{I \times
S_p}$ are given by characteristic power series of certain compact Hecke
operators acting on the flat family of Banach spaces of $p$-adic automorphic
forms for the definite unitary group $\U(3)$. Let us explain now how to deduce this statement from other ones in the litterature. Our method is actually a bit different from the one in~\cite{colemanmazur}.

\begin{pf} First, consider the unitary group $U/F$ in three variables attached to
the positive definite hermitian norm $(z_1,z_2,z_3) \mapsto \sum_{i=1}^3 {\rm Norm}_{E/F}(z_i)$ on $E^3$.
This group is necessarily quasi-split at all finite places of $F$, as we are
in odd dimension. The work of Rogawski defines by base change a bijection between automorphic 
representations $\pi$ of
$\GL_3(\AAA_E)$ satisfying (P1) and (P2) and stable tempered $L$-packets
$\Pi$ of $U$, which is compatible with a local base change (that Rogawski also defines) at all places
(\cite[\S 13]{roglivre}). The following lemma follows for instance from 
Rogawski's classification~\cite[p. 174, \S 13.1]{roglivre}.

\begin{lemma}\label{choixK} For any finite union of Bernstein component $\cB'$ of $\GL_3(\AAA_{E,S_0})$, there is a finite union $\cB$ of Bernstein components of $U(\AAA_{F,S_0})$ such that for each irreducible $\pi'$ in $\cB'$ and any irreducible $\pi$ of $U(\AAA_{F,S_0})$ whose base-change is $\pi'$, we have $\pi \in \cB$.
\end{lemma}

As $\cX(\rhob)$ has finitely many irreducible components, and by~\cite[Lemma 7.8.17]{bch}, there is a finite union of Bernstein components $\cB'$ of $\GL_3(\AAA_{E,S_0})$ such that for each $\Pi$ with $\rho_\Pi \in \cX(\rhob)$ we have $\Pi_{S_0} \in \cB'$. Choose an associated $\cB$ as in the above lemma. Up to enlarging $\cB$ if necessary, we may find a compact open subgroup 
$K \subset U(\AAA_{F,S_0})$ cutting exactly the union of components $\cB$,
in the sense that an irreducible representation $\tau$ of $U(\AAA_{F,S_0})$
belongs to $\cB$ if, and only if, $\tau^K \neq 0$.

Let $X/\Q_q$ be the $p$-adic eigenvariety of $U$ associated to 
$(\iota_p,\iota_\infty)$, to the tame level $K$ (spherical outside $S$), 
to the set of all places $S_p$ above $p$, and with respect to the spherical
Hecke algebra $\HH$ outside $S$ (see~\cite[Thm. 1.6]{che2}). By Rogawski's results and
\cite[Cor. 7.5.4]{bch} (or \cite[Cor. 7.7.1]{che}),
this eigenvariety carries a continuous pseudocharacter $T : \G_{E,S}
\longrightarrow \OO(X)$ such that for each refined classical point
$x$ associated to some $\Pi$ satisfying (P1)-(P3) the evaluation of $T$ at
$x$ is ${\rm trace}(\rho_\Pi)$. In turn, for each
$x \in X$ we have a natural associated semi-simple Galois representation
$\rho_x : G_{E,S} \rightarrow \GL_3(\overline{k(x)})$, whose residual
semi-simple representation will be denoted by $\rhob_x :  G_{E,S}
\rightarrow \GL_3(\overline{k_x})$, where $k_x$ is the residue field of $k(x)$ 
(a $\mathbb{F}_q$-algebra, as $k(x)$ is a $\mathbb{Q}_q$-algebra). The locus $X(\rhob) \subset X$
of $x$ such that $\rhob_x \simeq \rhob$ is an admissible closed and open subspace
of $X$ (this follows for instance from the Brauer-Nesbitt theorem and 
\cite[Lemma 3.9]{chdet}). By definition of
$\cB$ and the properties of $X(\rhob)$, there is a natural injection 
$j : \cZ \rightarrow X(\rho)(\Qpb)$, and $j(\cZ)$ is a Zariski-dense subset
of $X(\rhob)$. The universal property of $\cX(\rhob)$ defines
a canonical analytic map $\phi : X(\rhob) \rightarrow \cX(\rhob)$
such that $\phi \cdot j = {\rm id}_{\mathcal{Z}}$. The eigenvariety $X$ is also equipped with a
finite analytic map $\nu : X \rightarrow \WW^{I \times S_p}$, and property
(i) of \cite[Thm 1.6]{che2} ensures that $$\phi \times \nu : X(\rhob)
\rightarrow \cX(\rhob) \times
\WW^{I \times S_p}$$ is a closed
immersion. As a consequence, $\phi \times \nu$ induces an isomorphism $X(\rhob) \isomo
\cE(\rhob)$, and the last statement follows from the properties
of $X(\rhob)$.\end{pf}

\begin{remark}\label{bcanyn} {\rm The proof above is pseudo-character theoretic, hence extends verbatim to the case of any semi-simple $\rhob$. Moreover, this theorem would also hold for any odd $n$ by the same argument, if we appeal instead of Rogawski's work to the recent works of Moeglin (definition of the tempered $L$-paquets for quasi-split unitary groups) and Labesse (base change to $\GL(n)$) instead of the work of Rogawski, as long as $[F:\Q] \geq 2$ (assumption occuring in Labesse's base-change at the moment).} 
\end{remark}

An important property of the eigenvariety is the so-called {\it classicity criterion}. Say that $(\delta_{i,v}) \in \WW_0$ is {\it algebraic} if for each $(i,v) \in I\times S_p$ there exists an integer $k_{i,v}$ such that $\delta_{i,v}$ is the character  $u \mapsto u^{-k_{i,v}}$. Say that $(x,\delta) \in \cE(\rhob)$ is 
{\it of algebraic weight} if $\kappa(x)$ is algebraic. In this case, we denote by $k_{i,v}(x)$ the $k_{i,v}$ above and set $\varphi_{i,v}(x):=p^{k_{i,v}(x)}\delta_{i,v}(p)$. The next proposition is~\cite[Thm. 1.6.(vi)]{che2}. 

\begin{prop} \label{classcrit} Let $x \in \cE(\rhob)$ be of algebraic weight. Assume that $\forall v \in S_p$ :

{\rm (i)} $k_{1,v}(x) < k_{2,v}(x) < k_{3,v}(x)$, 
 
{\rm (ii)} ${\bf v}(\varphi_{1,v}(x)) < k_{2,v}(x) < {\bf v}(\varphi_{3,v}(x))$,

{\rm (iii)} $\forall i, j\,\,\,\in I$, $\varphi_{i,v}(x)\varphi_{j,v}(x)^{-1} \neq p$. 

\noindent Then $x \in \cZ$.

In particular, $\mathcal{Z}$ is an accumulation subset of $\cE(\rhob)$. More precisely, if $U \subset \cE(\rhob)$ is an open affinoid as in the statement of Thm.~\ref{thmeigen}, then $\mathcal{Z}$ is Zariski-dense in $U$ if and only if each connected component of $\kappa(U)$ contains an algebraic weight.
\end{prop}

Recall that if $A$ and $B$ are two subsets (of closed points)
of a rigid analytic space, we say that {\it $A$ accumulates at $B$} if for each $b \in B$ and for
each affinoid neighborhood $U$ of $b$, there exists an affinoid neighborhood $V \subset U$ of $b$ such that $A\cap V$ is Zariski-dense in $V$. 
An {\it accumulation subset} is a subset accumulating at itself. For instance the integers 
$\N$ is an accumulation subset of the closed unit disc over $\Q_p$, and the algebraic weights form an accumulation subset of $\WW_0^{I\times S_p}$.

\subsection{The local leaves and the infinite fern}\label{mainpic}

%\label{mainpic} From now on to the end of this section we discuss the main ideas and 
%strategy of the proof of Thm.\ref{mainthmobs}. 
%
%In the non-obstructed case, the space $\cX(\rhob)$ is the open unit ball of dimension $6[F:\Q]$ over $\Q_q$. In general it should still have dimension $6[F:\Q]$. 

Let $x=\rho_\Pi \in \cX(\rhob)$ be a modular point, and choose $\delta$ so that $(x,\delta) \in \mathcal{Z}$. By Theorem~\ref{thmeigen} and Prop.~\ref{classcrit}, we may find a basis of affinoid neighborhoods  $U_{x,\delta} \subset \cE(\rhob)$ of $(x,\delta)$ such 
that $\cZ$ is Zariski-dense in $U_{x,\delta}$, such that $\kappa(U_{x,\delta})$ is an open affinoid, and such that $\kappa: U_{x,\delta} \rightarrow \kappa(U_{x,\delta})$ is finite (necessarily surjective when restricted to any irreducible components of $U_{x,\delta}$). As $\beta$ is \'etale at $\kappa(x)$ we may even choose $U_{x,\delta}$ small enough so that the previous assertion holds with $\beta \,\, {\rm o}\,\, \kappa = \kappa_{HT} \,\,{\rm o}\,\, \pr_1$ instead of $\kappa$. In particular,  if $V=\beta(\kappa (U_{x,\delta}))$ then the induced map $\pr_1: U_{x,\delta} \rightarrow \kappa_{HT}^{-1}(V)$ is a finite map, thus {\it $\pr_1(U_{x,\delta})$ is a locally closed subset of $\cX(\rhob)$}, and even a closed subset of the admissible open $\kappa_{HT}^{-1}(V) \subset \cX(\rhob)$.

We have constructed this way up to
$6^{|S_p|}$ 
locally closed subspaces of $\cX(\rhob)$, inside $\cF(\rhob)$ and containing $x$, namely the  $\pr_1(U_{x,\delta})$, that we will call the {\it leaves of the fern} at 
$x$. Each of these leaves has equidimension $3[F:\Q]$; with this definition, only its germ at $x$ is canonical. By construction, each of these 
leaves also contains a Zariski-dense subset of modular points, 
a Zariski-dense subset of which even has the extra property that the eigenvalues of its cristalline Frobenius at $p$ are distinct (see e.g. the proof below of Lemma~\ref{lemmazd}), so that each of them 
will admit exactly $6^{|S_p|}$ associated refined modular points and so on... The {\it infinite fern} is the union of all the leaves constructed this way, from any modular point, namely 
$$\bigcup_{(x,\delta) \in \cZ} \pr_1(U_{x,\delta}) \subset \cX(\rhob),$$
a picture of which has been given in the introduction. Our main aim from now one will be to bound below the dimension of the Zariski-closure of the infinite fern. As explained in the introduction, the idea of our proof is to study the relative positions of the local leaves $\pr_1(U_{x,\delta})$, when 
$\delta$ varies, at a given modular point $x \in \cX(\rhob)$. 

\subsection{First reduction toward Theorem $A$} \label{stratpar} Let $${\cXm}' \subset \cXm$$ be the subset of modular
points $x$ such that for all $v \in S_p$, $\rho_{\Pi,v}$ is irreducible and the eigenvalues of its crystalline 
Frobenius are distinct and in $k(x)^*$. We claim first that ${\cXm}'$ is dense in $\cXm$:

\begin{lemma}\label{lemmazd} For any $x \in
\cXm$ and any affinoid neighbourghood $U$ of $x$ in $\cX(\rhob)$, we have ${\cXm}'\cap U
\neq \emptyset$. In particular,
${\cXm}'\neq \emptyset$. \end{lemma}

\begin{pf} We use an argument similar to~\cite[\S 7.7]{bch}: starting from any modular point we show that by moving in two at most two leaves of the fern we may 
find such a modular point by an argument of Newton polygon. So let $x \in \cXm$ be any modular point
and $U$ as in the statement. Choose any refinement $\delta$ of $x$ and consider an affinoid neighborhood $U_{x,\delta}$
of $(x,\delta)$ in $\cE(\rhob)$ as in~\S\ref{mainpic}, and such that ${\rm pr}_1(U_{x,\delta}) \subset
U$. Over this affinoid, the maps 
$y \mapsto |{\bf v}(\delta_{i,v}(p)(y))|$ are bounded by some integer $M$ by the maximum modulus principle. We may choose a point 
$y \in U_{x,\delta}$ of algebraic weight such that $k_{2,v}(y)-k_{1,v}(y)$ and $k_{3,v}(y)-k_{2,v}(y)$ are both bigger than $M+1$ for each $v$, and such that
$$\forall (i,v) \in S_p, \, \, \, \, {\bf v}(\delta_{i,v}(p)(y))={\bf v}(\delta_{i,v}(p)(x))={\bf v}(\varphi_{i,v}(x))-k_{i,v}(x).$$
In particular, such an $y$ is modular by Prop. \ref{classcrit}. Moreover, for each $(i,v) \in S_p$ we have 
${\bf v}(\varphi_{i,v}(y)) \neq k_{i+1,v}(y)$ where the index $i$ is taken mod $3$. In 
particular, up to replacing $x$ by $y$ if necessary, we may assume that the point $x$ we started from has this property as well, and also that 
we had chosen a refinement $\delta$ with the property that $$\forall (i,v) \in I \times S_p, \, \, \, {\bf v}(\delta_{i,v}(x)) \neq 0.$$
But then, a modular point $y$ as above has the property that for any $(i,j,v) \in I \times I \times S_p$, ${\bf v}(\varphi_{i,v}(y)) \neq k_{j,v}(y)$, so $\rho_{y,v}$ is absolutely irreducible 
by weak admissibility of its $\Dc$. Its crystalline Frobenius eigenvalues
have distinct valuations, so they are distinct and in $k(y)$. \end{pf}

\begin{remark} {\rm In the examples given in the appendix, the natural modular lift 
$\Symm^2(V_p(A))(-1)_{|\G_{E,S}}$ never belongs to ${\cXm}'$. }\end{remark}

If $Y$ is a rigid analytic space and $y \in Y$ a closed point we denote 
by $T_y(Y)$ the tangent space of $Y$ at $y$ (a $k(y)$-vector space).

\begin{thm}\label{redthm} For any $x \in {\cXm}'$, the image of the natural map $d \pr_1:$
$$\bigoplus_{y=(x,\delta) \in \cE(\rhob)} T_y(\cE(\rhob)) \longrightarrow T_x(\cX(\rhob))$$
has dimension at least $6[F:\Q]$.
\end{thm}

(Note that for any $x \in {\cXm}'$, the residue field of any $(x,\delta) \in \cE(\rhob)$
is $k(x)$ by definition.)

Let us first show that Thm~\ref{redthm} implies Thm~\ref{mainthmobs}. Let $W$ be the Zariski-closure of $\cXm$ in $\cX(\rhob)$, equipped with its reduced structure. Let $C$ be an irreducible component of  $W$. Recall that affinoid algebras are excellent and Jacobson, so the regular locus of any reduced rigid analytic space is a Zariski-open and Zariski-dense subspace (see e.g. \cite[\S 1]{conrad}). Thus may choose some modular point $x_0 \in C \cap \cXm$ such that $x_0$ is a smooth point of $W$ (hence of $C$). Choose $U \subset \cX(\rhob)$ an affinoid neighborhood of $x_0$, which is small enough so that $U \cap W$ is regular, and equal to $U\cap C$. By Lemma~\ref{lemmazd}, we may find a modular point $x$ in $U$ that furthermore belongs to ${\cXm}'$, by construction $x \in U\cap W=U\cap C$. We are going to apply Theorem~\ref{redthm} to this point $x$. Note that for each associated refined
modular point $(x,\delta) \in \cE(\rhob)$, 
\begin{equation}\label{inclusionTan} {\rm Im}\left(\,\,d \pr_1 :  T_{(x,\delta)}(\cE(\rhob)) \longrightarrow T_x(\cX(\rhob))\,\,\right) \subset T_x(C).\end{equation}
Indeed, let $U_{x,\delta} \subset \cE(\rhob)$ be an affinoid neighborhood of $(x,\delta)$ in $\cE(\rhob)$ as in~\S\ref{mainpic}. As $\cZ$ is Zariski-dense in $U_{x,\delta}$ we have $\pr_1(U_{x,\delta}) \subset W$. Thus if $U_{x,\delta}$ is chosen small enough, then $\pr_1(U_{x,\delta}) \subset U\cap W=U \cap C$, and \eqref{inclusionTan} follows. As a consequence, Theorem~\ref{redthm} implies that $\dim T_x(C)\geq 6[F:\Q]$. As $x$ is a smooth point of $C$, it follows that $C$ itself has dimension at least $6[F:\Q]$, and we are done.

\medskip

The end of the paper will be devoted to the proof of Theorem~\ref{redthm}. As explained in the introduction, there 
are two important ingredients, treated in the next two chapters. The first one is a purely local result on the 
deformation space of a crystalline representation, and the other one a geometric 
property of eigenvarieties at non-critical classical points. As a motivation, the reader may already have a look to the end of the proof in~\S\ref{mainpf}.

\section{The linear span of trianguline deformations}\label{tridef}
\newcommand{\Gp}{{G_{\Q_p}}}

\subsection{The setting}\label{setdef} Let $L$ be a finite extension of $\Q_p$, 
$\Gp=\Gal(\Qpb/\Q_p)$, $n\geq 1$ an integer and let $V$ be a continuous
$L$-linear representation 
of $\Gp$ of dimension $n$. If $m \in \Z$, we set $V(m)=V\otimes \chi^m$ where $\chi : \Gp \rightarrow \Z_p^*$ is the cyclotomic character.

Denote by $\Cc$ the category whose objects are the finite dimensional local $\Qp$-algebras 
$A$ equipped with an isomorphism $\pi: A/m_A \isomo L$ (so $A$ has a unique
structure of $L$-algebra such that $\pi$ is $L$-linear), and whose morphisms are the local 
$L$-algebra homomorphisms. Let $$\cX_V : \Cc \rightarrow {\rm Sets}$$ be the deformation functor 
of $V$ to $\Cc$. Recall that for an object $A$ of $\Cc$, $\cX_V(A)$ is the set of isomorphism 
classes of pairs $(V_A,\pi)$ where $V_A$ is a free $A$-module of rank $n$ equipped with a 
continuous $A$-linear representation of $\Gp$ and $\pi : V_A \otimes_A L \isomo V$ is
an 
$L[\Gp]$-isomorphism. Assume that $\End_\Gp(V)=L$ and that $\Hom_\Gp(V,V(-1))=0$, so that $H^i(\Gp,\End_L(V))=0$ 
if $i\geq 2$ and $\dim_L H^1(\Gp,\End_L(V))=n^2+1$ by Tate's results on Galois cohomology of 
local fields. By standard results of Mazur~\cite{mazurdef}, $\cX_V$ is pro-representable, formally 
smooth over $L$, of dimension $n^2+1$:

\begin{prop}
$\cX_V \simeq \Spf(L[[X_0,\cdots,X_{n^2}]])$. 
\end{prop}
Our main aim in this section will be the study and comparison
of a collection of subfunctors of $\cX_V$ when $V$ is crystalline in the sense of Fontaine. Two of them are the subfunctors of crystalline and 
Hodge-Tate deformations, but we shall actually mostly be interested in an additional collection of 
subfunctors introduced in~\cite[Chap. 2]{bch} 
under the name of {\it trianguline deformation functors}. They are some sorts of analogues of the 
ordinary deformation functors of an ordinary representation defined by Mazur
but that apply to any crystalline 
representation; they depend on the datum of a {\it refinement} of the crystalline representation.
In rank $2$ they are close to the $h$-deformation functors previously defined by Kisin in~\cite{kisin}. They belong to an even more general collection of {\it paraboline deformation
functors} that we define and study below using the theory of $\fg$-modules over the Robba ring.

\newcommand{\Dr}{D_{\rm rig}}

\newcommand{\Fil}{{\rm Fil}}

\subsection{Paraboline deformation functors of $\fg$-modules over the Robba
ring}\label{paradef} Let $A$ be a finite dimensional $\Q_p$-algebra, equipped with the
topology given by any norm of $\Q_p$-vector space on $A$. Recall that the Robba ring with coefficients in $A$ is
the $A$-algebra $\Ro_A$ of power series $$f=\sum_{n\in \Z} a_n
(z-1)^n, \, \, \, \, \, a_n
\in A,$$ converging on some annulus of $\C_p$ of the form $r(f) \leq |z-1|<1$.
It is a topological $A$-algebra, namely an inductive limit of projective limits of
affinoid algebras. It is equipped with $A$-linear commuting actions of $\varphi$ and
$\Gamma:=\Z_p^*$ given by the formulae $\varphi(f)(z)=f(z^p)$ and
$\gamma(f)(z)=f(z^\gamma)$, the action of $\Gamma$ being continuous. We set
$\Ro=\Ro_{\Q_p}$, so $\Ro_A = \Ro \otimes_{\Qp} A$.

A $\fg$-module over $\Ro_A$ is a finite free
$\Ro_A$-module $D$ equipped with commuting semi-linear actions of $\varphi$ and
$\Gamma$, such that $\Ro\varphi(D)=D$ (``non degeneracy of $\varphi$'') and such that the action of $\Gamma$
on $D$ is continuous.\footnote{It means the following. Fix an
$\Ro$-basis of $D$ and for $\gamma \in \Gamma$ let $M_\gamma$ be the matrix of
$\gamma$. Then the coefficients of all the $M_\gamma$, $\gamma \in \Gamma$,
have to converge on some fixed annulus $r < |z-1| < 1$, and for any $r< r' < 1$
in $p^{\Q}$ the map $\gamma \mapsto
M_\gamma$, $\Gamma \rightarrow M_n(\OO[r,r'])$ is continuous, where
$\OO[r,r']$ is the affinoid algebra of analytic functions on the annulus
$r \leq |z-1| \leq r'$ for the sup. norm. } We refer to \cite[\S 2.2]{bch} for 
a summary of the basic facts concerning the theory of
$\fg$-modules over the Robba ring. If $V$ is a $\Q_p$-representation of $\Gp$, denote by $\Dr(V)$ its associated $\fg$-module
over $\Ro$. Recall that the 
functor $V \mapsto \Dr(V)$ is a tensor equivalence between the category of $\Q_p$-representations of $\Gp$ and 
\'etale $\fg$-modules over $\Ro$ (Fontaine, Colmez-Cherbonnier, Kedlaya). In
this equivalence, if $V$ is
a $\Q_p$-representation of $\Gp$ equipped with a $\Q_p[\Gp]$-linear action of $A$, then $V$ is
free as $A$-module if and only if $\Dr(V)$ is a $\fg$-module over $\Ro_A$
by~\cite{bch} Lemma 2.2.7.

Let $D$ be a $\fg$-module over $\Ro_L$ and 
let $\mathcal{P}=(\Fil_i(D))_{i \in \Z}$ be an increasing filtration of $D$ by $\fg$-submodules of $D$, 
each $\Fil_i(D)$ being a direct summand as $\Ro_L$-submodule. As a convention, we shall
always assume that such filtrations are normalized so that $\Fil_i(D)=0$ if $i\leq 0$ and $\Fil_i(D)=D$ if $\Fil_i(D)=\Fil_{i-1}(D)$. 
Let us define the {\it paraboline deformation functor} of $D$ associated to $\mathcal{P}$ $$\cX_{D,\mathcal{P}} : \mathcal{C} \rightarrow {\rm Sets}$$
as follows. For an object $A$ of $\mathcal{C}$ define $\cX_{D,\mathcal{P}}(A)$ as the 
set of isomorphism classes of triples $(D_A, \Fil_i(D_A),\pi)$ where $D_A$ is a 
$\fg$-module over $\Ro_A$, $\Fil_i(D_A)$ is  an increasing filtration of $D_A$ by  
$\fg$-submodules over $\Ro_A$ of $D$, each $\Fil_i(D_A)$ being a direct summand as $\Ro_A$-submodule, and 
where $\pi : D_A \otimes_A L \isomo D$ is a $\fg$-module $L$-isomorphism such that 
for each $i\geq 0$ we have $\pi(\Fil_i(D_A))=\Fil_i(D)$. To be explicit, by an 
isomorphism between two such triples we mean an $\Ro_A$-linear isomorphism 
$\psi : D_A  \rightarrow D'_A$ commuting with $\varphi$ and $\Gamma$, 
mapping $\Fil_i$ onto $\Fil_i$, and such that $\pi=\pi' \cdot  \psi$.

When $\Fil_1(D)=D$ then $\cX_{D,\mathcal{P}}$ is simply the
{\it deformation functor} of $D$, and we shall simply denote it by $\cX_D$. 
When $\Fil_i$ is a complete flag in $D$, {\it i.e.}
$\rk_{\Ro_L}(\Fil_i(D))=i$ for $i=1,\dots,\rk_{\Ro_L}(D)$, $\cX_{D,\mathcal{P}}$ has been studied 
in \cite[\S 2.3]{bch} under the name of trianguline or triangular deformation functor of $D$. In 
the case $D=D_{\rm rig}(V)$ is \'etale, then $- \mapsto D_{\rm rig}(-)$ induces an 
equivalence $\cX_V \isomo \cX_{D}$ by \cite[Prop. 2.3.13]{bch}, and we shall set as 
well $$\cX_{V,\mathcal P}:=\cX_V \times_{\cX_D} \cX_{D,\mathcal P}$$ 
( $\isomo \cX_{D,\mathcal{P}}$). When furthermore each $\Fil_i(D)$ is \'etale, hence 
is the $D_{\rm rig}$ of a unique Galois stable filtration $\Fil_i(V)$, then 
$\cX_{V,{\mathcal P}}(A)$ coincides with the isomorphism classes of deformations 
of $V$ to $A$ equipped with some Galois stable filtrations over $A$ lifting 
$\Fil_i(V)$; however, there is no such description for a general $\mathcal{P}$.

By a well-known result of Mazur~\cite[Prop. 1]{mazurdef}, $\cX_D$ is pro-representable if $\End_{\fg}(D)=L$ 
and $D$ is \'etale. Let $\End_{\fg,\mathcal{P}}(D)$ be the $L$-vector space of 
$\fg$-module endomorphisms of $D$ preserving $\Fil_i(D)$ for each $i$. 

\begin{prop}\label{representability} $\cX_{D,\mathcal{P}}$ admits a versal pro-deformation. If
${\rm End}_{\fg,\mathcal{P}}(D)=L$, then it is pro-representable. 
\end{prop}

\begin{pf} Set $F=\cX_{D,\mathcal{P}}$. We obviously have
$|F(L)|=1$. Let $A' \rightarrow A$ and $A'' \rightarrow A$
two morphisms in $\mathcal{C}$. By Schlessinger's criterion~\cite{schless} we need to check that the natural map 
\begin{equation}\label{grotcrit} F(A' \times_A A'') \rightarrow
F(A') \times_{F(A)} F(A')\end{equation} is bijective whenever $A'' \rightarrow A$ is
surjective.

Let $E : \mathcal{C} \rightarrow {\rm Ens}$ be the following
functor: for an object
$A$ of $\mathcal{C}$, let $E(A)$ be the set of $\fg$-module structures
over $\Ro_A$ on the $\Ro_A$-module 
$D\otimes_L A$ that preserve the {\it constant filtration}
$\Fil_i(D)\otimes_L A$ and that induce the 
$\fg$-module structure of $D$ via the given map $A \rightarrow L$. For all $A,A'$ and $A''$,
we claim that (\ref{grotcrit}) is bijective when $F$ is replaced by
$E$. Indeed, $\Ro_{A'\times_A A''} = \Ro_A' \times_{\Ro_A} \Ro_{A''}$, and
if $R(A)$ denotes the ``parabolic'' $\Ro_A$-algebra of $\Ro_A$-module endomorphisms of
$D\otimes_L A$
preserving $\Fil_i(D)\otimes_L A$, then $R$ commutes with finite fiber
products in $\mathcal{C}$ as well. Moreover, if we fix some basis $(e_i)_{i=1}^n$ of $D$ adapted to the
filtration $\Fil_i(D)$, then it is equivalent to give an element of $E(A)$ and
a collection of matrices $M_\ast$, for $\ast=\varphi$ or in $\Gamma$, lifting the matrices in $M_n(\Ro_L)$ defined by the $\fg$-module structure of $D$, and satisfying $y(M_x)M_y=x(M_y)M_x$ for all $x, y \in \{\varphi, \gamma \in \Gamma\}$. Indeed, the non degeneracy of $\varphi$ and the continuity of $\Gamma$ follow automatically from these properties for $D$.

Set $R^1(A)=\Ker(R(A)^*\rightarrow R(L)^*)$, the group $R^1(A)$ acts by conjugacy on $E(A)$ and we have $F(A)=R^1(A)\backslash E(A)$ by definition. If $A'' \rightarrow A$ is any surjection in $\mathcal{C}$, then $R^1(A'') \rightarrow R^1(A)$ is surjective as well, as parabolic subgroups of $\GL_n$ are smooth over $\Q_p$. Moreover, $R^1$ obviously commutes with finite fiber products in $\mathcal{C}$. It follows at once that (\ref{grotcrit}) is surjective whenever $f: A'' \rightarrow A$ is. It is even bijective if for all $x \in E(A'')$, the group homomorphism induced on stabilizers 
\begin{equation}\label{stabmap} R^1(A'')_x \rightarrow R^1(A)_{E(f)(x)}\end{equation} is surjective (these observations really are Mazur's, see~\cite[p. 390]{mazurdef}). This is clearly satisfied if $A=L$ as the latter stabilizer is trivial. As a consequence, conditions (H1) and (H2) of Schlessinger hold. In particular $F(L[\varepsilon])$ is an $L$-vector space and
we shall see in Prop.~\ref{esptang} below that it is always finite dimensional (contrary to $E(L[\varepsilon])$
in this $\fg$-module context), hence (H3) and the
first part of the statement follow. When $\End_{\fg,\mathcal{P}}(D)=L$, then $R^1(A)_x=A^*$ for any $A$ and $x \in E(A)$ by Remark~\ref{forgetpi} below, hence (\ref{stabmap}) is surjective, and (H4) holds.
\end{pf}

\begin{remark}\label{forgetpi}{\rm  If $\End_{\fg,\mathcal{P}}(D)=L$, then a
simple induction on the lenght of $A$ shows that for each element $(D_A,\Fil_i,\pi)$ of
$\cX_{D,\mathcal{P}}(A)$, we have $\End_{\fg,(\Fil_i)}(D_A)=A$. In particular, 
$\cX_{D,\mathcal{P}}(A)$ is in canonical bijection with the isomorphism classes 
of pairs $(D_A,\Fil_i)$ lifting $D$ with its filtration (forgetting the $\pi$).}
\end{remark}

If $D_1$ and $D_2$ are two $\fg$-modules over $\Ro_L$, then $\Hom_{\Ro_L}(D_1,D_2)$ 
has a natural structure of $\fg$-module as follows: for any 
$u \in \Hom_{\Ro_L}(D_1,D_2)$, set $\gamma(u)(x):=\gamma(u(\gamma^{-1}(x)))$ and 
$\varphi(u)(\varphi(x)):=\varphi(u(x))$. This last formula defines a unique element 
$\varphi(u) \in  \Hom_{\Ro_L}(D_1,D_2)$ as $D_1$ admits an $\Ro_L$-basis in 
$\varphi(D_1)$ by definition. In particular, $\End_{\Ro_L}(D)$ is a $\fg$-module of 
rank $n^2$ over $\Ro_L$ and we check at once with the formulas above that the
$\Ro_L$-submodule $$\End_{\mathcal{P}}(D):=\{ u \in
\End_{\Ro_L}(D), u(\Fil_i(D)) \subset \Fil_i(D) \, \, \, \, \, \forall i \in \N\}$$
is a $\fg$-submodule of $\End_{\Ro_L}(D)$ (which is actually a direct
summand as $\Ro_L$-submodule). For $i\in \Z$ we set 
$\gr_i(-):=\Fil_i(-)/\Fil_{i-1}(-)$. If $n_i=\rk_{\Ro_L}(\gr_i(D))$, we obviously have 
\begin{equation}\label{rkep} \rk_{\Ro_L}(\End_{\mathcal{P}}(D)) = \sum_{i \leq j} n_in_j.\end{equation}

We refer to~\cite[\S 2]{colmeztri} and~\cite{liu} for the main properties of the cohomology of $\fg$-modules over $\Ro$, that we shall 
denote by $H^i_{\fg}(-)$. In particular,
$H^0_{\fg}(\End_{\mathcal{P}}(D))=\End_{\fg,\mathcal{P}}(D)$ and $H^2_{\fg}(\End_{\mathcal{P}}(D))$ is dual to the 
$L$-vector space of $\fg$-morphisms $D \rightarrow D(-1)$ preserving $\mathcal{P}$.

\begin{prop}\label{esptang} \begin{itemize}
\item[(i)] If $\Hom_{\fg}(\gr_i(D),D/\Fil_i(D))=0$  for each $i$, then $\cX_{D,\mathcal{P}}$ is a subfunctor of $\cX_D$.
\item[(ii)] There is a natural isomorphism $\cX_{D,\mathcal{P}}(L[\varepsilon]) \isomo H^1_{\fg}(\End_{\mathcal{P}}(D))$ and
$$\dim_L \cX_{D,\mathcal{P}}(L[\varepsilon]) = \dim_L H^0_{\fg}(\End_{\mathcal{P}}(D))+ \dim_L H^2_{\fg}(\End_{\mathcal{P}}(D)) +  \sum_{i \leq j} n_in_j.$$
\item[(iii)] If $H^2_{\fg}(\End_{\mathcal{P}}(D))=0$ then $\cX_{D,\mathcal{P}}$ is formally smooth over $L$. \end{itemize}
\end{prop}

\begin{pf} Part (i) follows as in \cite[Lemma 2.3.7]{bch}. The second part of (ii) follows from the first one, (\ref{rkep}) 
and the Euler characteristic formula of Liu~\cite{liu}. Let us check the first part of (ii). It is a semi-linear analogue of a 
well-known fact in the context of group representations. Using the notations
of the proof of Prop.~\ref{representability}, for all $A$ we have a natural
identification $\cX_{D,\mathcal{P}}(A)= R^1(A)\backslash E(A)$. 

Consider the $\Ro_L[\varepsilon]$-module $D_0:=D\otimes_L L[\varepsilon]=D \oplus \varepsilon D$ equipped with its
constant filtration as in the proof of Prop.~\ref{representability}. Any
element of $E(L[\varepsilon])$ is given by unique elements $c_\ast \in \End_\Ro(D)$, where $\ast=\varphi$ or 
$\gamma \in \Gamma$, satisfying 
the formulas
$$\widetilde{\gamma}(x):=\gamma(x) +c_\gamma(\gamma(x))\varepsilon , \,\,\, \widetilde{\varphi}(x):=\varphi(x)+ c_\varphi(\varphi(x)) \varepsilon, \, \, \, \forall x\in D \subset D\oplus \varepsilon D.$$
Note that the map $c_\varphi$ (resp. $c_\gamma$) is well-defined as $D$ admits an
$\Ro_L$-basis in $\varphi(D)$ (resp. as $\gamma$ is bijective on $D$). A straightforward computation shows 
that the commutation of $\widetilde{\varphi}$ and $\widetilde{\gamma}$ is equivalent to the relation $(\varphi-1)c_\gamma=(\gamma-1)c_\varphi$ in $\End_\Ro(D)$, for the $\fg$-module structure recalled above on $\End_\Ro(D)$. The continuity of the action of $\Gamma$ is automatic, as well as the non degeneracy of $\widetilde{\varphi}$. Moreover, $*$ 
preserves each $\Fil_i D \otimes_L L[\varepsilon]$ if, and only if, $c_\ast \in \End_{\mathcal{P}}(D)$. 
In other words, $c_\ast$ is a $1$-cocycle in $\End_{\mathcal{P}}(D)$ and
$E(L[\varepsilon])$ coincides with this space of cocycles. Another straightforward computation shows that 
two such lifts are isomorphic in
$\cX_{D,\mathcal{P}}(L[\varepsilon])=E(L[\varepsilon])/R^1(L[\varepsilon])$ if and only if the associated 
$1$-cocycles $c_\ast$ and $c'_\ast$ differ by $(\ast-1)u$ for some $u \in \End_{\mathcal{P}}(D)$ independent 
of $\ast$, hence a coboundary, which concludes the proof of (ii) (the $L$-linearity part of the statement is 
immediate). 

For part (iii), let $A$ be an object of $\mathcal{C}$, $m$ its maximal ideal, and $I$ an ideal of $A$ such 
that $mI=0$. Fix an element of $E(A/I)$, hence a $\Ro_{A/I}$-linear $\fg$-module structure 
on $D\otimes_L A/I$ preserving $\Fil_i(D)\otimes_L A/I$, and consider the problem of lifting this structure to 
an element of $E(A)$. First, we may certainly lift the actions of $\varphi$ and of $\Gamma$ to $D\otimes_L A$ in such a 
way that they preserve $\Fil_i(D)\otimes_L A$. If $\widetilde{\varphi}$ and
$\widetilde{\gamma}$ denote such a lift, consider the $\Ro_L$-linear map $\langle \widetilde{\varphi}, 
\widetilde{\gamma}\rangle \in \End_{\mathcal{P}}(D) \otimes_L I$ defined by the formula $$\forall x \in D\otimes_L A, (\widetilde{\varphi}\widetilde{\gamma}-\widetilde{\gamma}\widetilde{\varphi})(x)=\langle \widetilde{\varphi}, \widetilde{\gamma}\rangle(\varphi\gamma(\overline{x})),$$
where $\overline{x}$ is the image of $x$ in $D$. This map $\langle \widetilde{\varphi}, \widetilde{\gamma}\rangle$ is well-defined as $\gamma$ is bijective on $D$ and as $D$ admits an $\Ro$-basis in $\varphi(D)$. The elements $\widetilde{\varphi}$ and $\widetilde{\gamma}$ are uniquely defined up to adding any elements $u(\varphi(\overline{\cdot}))$ and $v(\gamma(\overline{\cdot}))$ where $u,v \in \End_{\mathcal{P}}(D) \otimes_L I$, and $$\langle \widetilde{\varphi}+u(\varphi(\bar\cdot)), \widetilde{\gamma}+v(\gamma(\bar\cdot))\rangle = \langle \widetilde{\varphi}, \widetilde{\gamma}\rangle + (\varphi-1)v-(\gamma-1)u.$$
This concludes the proof as $H^2_{\fg}(\End_\mathcal{P}(D))=\End_\mathcal{P}(D)/(\varphi-1,\gamma-1)\End_\mathcal{P}(D)$. Note that as in the classical case, it is straightforward to check that the set of lifts is either empty or an affine space under $H^1_{\fg}(\End_{\mathcal{P}}(D))\otimes_L I$. \end{pf}
\newcommand{\Aut}{{\rm Aut}}
Let $I \subset \N$ be the (finite) subset of jumps of the filtration
$\Fil_i$, {\it i.e.} the integers $i$ such that $\gr_i(D)\neq 0$. There is a natural functor morphism $\cX_{D,\mathcal{P}} \rightarrow \prod_{i \in I} \cX_{\gr_i(D)}$ sending an object $(D_A,\Fil_i,\pi)$ to $(\gr_i(D_A),\gr_i(\pi))_{i\in I}$. 

\begin{prop}\label{relsmooth} If $H^2_{\fg}(\Hom_\Ro(D/\Fil_i(D),\gr_i(D)))=0$ for each $i$, 
then $\cX_{D,\mathcal{P}} \rightarrow \prod_{i \in I} \cX_{\gr_i(D)}$ is formally smooth. In particular, it is surjective on points.
\end{prop}

\begin{pf} The proposition is obvious when $|I|\leq 1$ so assume that
$|I|>1$. Considering the filtration $\Fil'$ on $D$ such that ${\Fil'}_1(D)=\Fil_j(D)$ for $j$ the biggest 
integer such that $\Fil_j(D) \neq D$, $\Fil'_2(D)=D$, and arguing by induction on $|I|$, we may assume that 
$I=\{1,2\}$. So $D$ is an extension of $D_1:=D/\Fil_1(D)$ by $D_2:=\Fil_1(D)$, which defines a unique 
class $c \in H^1_{\fg}(\Hom_\Ro(D_1,D_2))$. Let $A$ be an objet of $\mathcal{C}$, $m$ its maximal ideal, 
and $I$ an ideal of $A$ such that $mI=0$. Fix $(U,\Fil_i(U),\pi) \in \cX_{{\mathcal{P}},D}(A/I)$, it defines 
as above an element $c_{A/I} \in H^1_{\fg}(\Hom_\Ro(U/\Fil_1(U),\Fil_1(U)))$
that maps to $c$ modulo $m$. For $i=1,2$, choose $(D'_i,\pi_i) \in \cX_{D_i}(A)$ lifting respectively 
$U/\Fil_1(U)$ and $\Fil_1(U)$. We have to show that the natural map 
$$H^1_{\fg}(\Hom_\Ro(D'_1,D'_2)) \longrightarrow
H^1_{\fg}(\Hom_\Ro(U/\Fil_1(U),\Fil_1(U)))$$
is surjective. But its cokernel injects into $H^2_{\fg}(\Hom_{\Ro}(D_1,D_2))$.
\end{pf}

\begin{example} (The rank $1$ case) {\rm Let $A$ be an object of $\mathcal{C}$
and let $\delta : \Q_p^* \rightarrow A^*$ be a continuous character.  We denote by
$\Ro_A(\delta)$ the $\fg$-module of rank $1$ over $\Ro_A$ having a basis $e$ such
that $\varphi(e)=\delta(p)e$ and $\gamma(e)=\delta(\gamma)e$ for all $\gamma
\in \Gamma$. By~\cite[Prop. 2.3.1]{bch}, each $\fg$-module of rank $1$ over $\Ro_A$ has 
the form $\Ro_A(\delta)$ for a unique character $\delta$ as above. In
particular, if $D$ is a $\fg$-module of rank $1$ over $\Ro_L$, then $\cX_{D}
\simeq \Spf(L[[X,Y]])$.}
\end{example}

\subsection{Crystalline $\fg$-modules}\label{triref}

If $D$ is any $\fg$-module over $\Ro$, we set $\cDc(D) = (D[1/t])^\Gamma$, 
where $t={\rm log}(z) \in \Ro$. As $\Ro$ is a domain and ${\rm Frac}(\Ro)^\Gamma=\Q_p$, a standard 
argument shows that 
\begin{equation}\label{injectdc}
\cDc(D) \otimes_L \Ro_L[1/t] \rightarrow D[1/t]\end{equation}
is injective for any $D$, so we have $\dim_{\Q_p} \cDc(D) \leq {\rm rk}_\Ro
D$. Mimicking Fontaine, we say that $D$ is {\it crystalline} if the equality holds. By left exactness of 
$- \mapsto \cDc(-)$, any $\fg$-module subquotient of a crystalline $D$ is also crystalline. Note 
that $\cDc(D)$ has a $\Q_p$-linear invertible action of $\varphi$.

\begin{lemma}\label{trivprop} Let $D$ be a $\fg$-module over $\Ro_L$. The map $$D' \mapsto \cDc(D')$$
induces a bijection 
between the set of crystalline $\fg$-submodules of $D$ which are direct summand as $\Ro_L$-module, and the set of $\varphi$-stable subspaces of
$\cDc(D)$. The inverse bijection is $W \mapsto (\Ro[1/t]\cdot W)\cap D$. 
\end{lemma}

Indeed, a $\fg$-submodule $D' \subset D$ is a direct summand as 
$\Ro_L$-module if, and only if $D'=D'[1/t]\cap D$ by \cite[Prop. 2.2.2]{bch}, so the
lemma follows from (\ref{injectdc}). In particular, there is always 
a biggest crystalline $\fg$-submodule of $D$, namely $(\Ro[1/t] \cdot D[1/t]^{\Gamma})\cap D$. 
Furthermore, when $D$ is crystalline then Lemma~\ref{trivprop} induces a rank preserving bijection 
between the set of filtrations of $D$ as in \S\ref{paradef} and the set of increasing filtrations 
$(\Fil_i)_{i \in \Z}$ of $\cDc(D)$ by $\varphi$-stable $\Q_p$-vector space
(normalized with the same conventions as before), we shall often identify 
such filtrations below. In particular, the paraboline deformation functors of $D$ can
(and will) be viewed as attached to such filtrations. 

\noindent Consider the following properties of a crystalline $\fg$-module of rank $n$ over $\Ro_L$:\ps
\begin{itemize}
\item[(i)] $\varphi$ has $n$ distinct eigenvalues in $L^*$ on $\cDc(D)$.\ps
\item[(ii)] For any two such eigenvalues $\phi, \phi'$, we have 
$\phi'\phi^{-1} \neq p$. \ps
\item[(iii)] $\End_{\fg}(D)=L$.\ps
\end{itemize}

\begin{cor}\label{corparacrys} Let $D$ be a crystalline $\fg$-module over $\Ro_L$ satisfying
(i), (ii) and (iii). Then each paraboline deformation functor $\cX_{D,\mathcal P}$ of $D$
is a pro-representable subfunctor of $\cX_D$, formally smooth over $L$ of dimension 
$\rk_{\Ro_L}(\End_{\mathcal P}(D))+1$.
\end{cor}

\begin{pf} The pro-representability follows from Prop.~\ref{representability} and
(iii). As $D$ is crystalline, remark that 
$\Hom_{\fg, L}(D,D') \subset \Hom_{L[\varphi]}(\cDc(D),\cDc(D'))$ for any
$\fg$-module $D'$ over $\Ro_L$. By
Prop.~\ref{esptang}, the subfunctor property follows from (i) and the formal
smoothness from (ii). \end{pf}

For any $D$, Berger defines as well a decreasing, exhaustive, filtration on $\cDc(D)$ called the {\it Hodge-filtration}, 
that we shall denote with upper indices $\Fil^i$ for $i\in \Z$
(see~\cite{berger1},~\cite{berger2}, as well as \cite[\S 2.2.7]{bch}). 
A first important result of Berger is that when $D=\Dr(V)$, the filter $\varphi$-module $\cD_\ast(D)$ 
is canonically isomorphic to the classical $D_\ast(V)$ defined by Fontaine. In particular, $V$ is 
crystalline if, and only if, $D$ is. Another beautiful result of Berger~\cite{berger2} is that 
$D \mapsto \cDc(D)$ induces an exact $\otimes$-equivalence of category between the category of 
crystalline $\fg$-modules over $\Ro$ and the category of filtered $\varphi$-vector
space, \'etale corresponding to weakly admissible.

Let $D$ be a crystalline $\fg$-module over $\Ro_L$. Define 
$\cX_{D,{\rm crys}} \subset \cX_D$ as the subfunctor of deformations $(D_A,\pi)$ such that $D_A$ is 
crystalline, viewed as $\fg$-module over $\Ro$. As the crystalline $\fg$-modules are stable 
by direct sums, subquotients, and tensor products, $\cX_{D,{\rm crys}} \rightarrow \cX_{D}$ is 
relatively representable (Ramakrishnan's criterion). 

\begin{prop}\label{desccrys} Let $D$ be a crystalline $\fg$-module over $\Ro_L$.
\begin{itemize}\item[(a)] If $D$ is equipped with an action of an object $A$ of $\mathcal{C}$, then $D$ is free as $\Ro_A$-module if and only if $\Fil^i(\cDc(D))$ is free as $A$-module for each $i \in \Z$.\ps

\item[(b)] The functor $-\mapsto \cDc(-)$ induces a canonical bijection between $\cX_{D,{\rm crys}}(A)$ and the set of isomorphism classes of pairs 
$(E_A,\pi)$ where $E_A$ is a filtered $A[\varphi]$-module such that $\Fil^i
E_A$ 
is free and direct summand as $A$-module and $\pi : E_A\otimes_A L \isomo \cDc(D)$ is an $A$-linear 
isomorphism in the category of filtered $\varphi$-modules.

\item[(c)] If $D$ satisfies properties (i) and (iii) above, then $\cX_{D,{\rm crys}}$ is formally smooth over $L$ of dimension $$1+\sum_{i < j} n_in_j,$$ where $n_i = \dim_L \Fil^i(\cDc(D))/\Fil^{i+1}(\cDc(D))$ for all $i \in \Z$.
\end{itemize}
\end{prop}

\begin{pf} Part (a) holds by the same proof as~\cite[Lemma 2.2.7]{bch}. 
Part (b) follows from (a) and Berger's equivalence. Using (b), part (c) is a simple exercise on deformations of filtered $\varphi$-modules
using (i) and (iii) that is left to the reader.\end{pf}

Let $D$ be a crystalline $\fg$-module over $\Ro_L$ and let $W \subset \cDc(D)$ be an $L$-subvectorspace. 
We say that $W$ is {\it non-critical} (in $\cDc(D)$) if there exists an integer $j \in \Z$ such that $\cDc(D)=W\oplus 
\Fil^j(\cDc(D))$. If $\mathcal{P}=(\Fil_i)$ is a $\varphi$-stable filtration of $\cDc(D)$ as above, we 
say that $\mathcal{P}$ is non-critical if $\Fil_i$ is non-critical for all $i\geq 0$.

\begin{prop}\label{crysnoncrit} Let $D$ be a crystalline $\fg$-module over $\Ro_L$ satisfying (i) and (iii). For any non-critical filtration $\mathcal{P}$, we have $\cX_{D,{\rm crys}} 
\subset \cX_{D,\mathcal{P}}$. 
\end{prop}

\begin{pf} By (i), $\cX_{D,\mathcal{P}}$ is a subfunctor of $\cX_{D}$. Fix 
$(D_A,\pi) \in \cX_{D,{\rm cris}}(A)$ and let $W \subset \cDc(D)$ be an element of $\mathcal{P}$. 
By (i) again, there is a unique $\varphi$-stable $A$-submodule $W_A \subset \cDc(D_A)$ which is free over 
$A$ and such that $\pi(W_A)=W$. Note that if $W \subset W'$ then $W_A \subset W'_A$.
Set $T_A=(\Ro[1/t]W_A)\cap D$, we will show that it is free over $\Ro_A$ and lifts $(\Ro[1/t]W)\cap D$. 

By Prop.~\ref{desccrys} (a) and Berger's equivalence, it is equivalent to check that if we equip $W_A$ with the induced 
Hodge-filtration of $\cDc(D_A)$, then $\Fil^i(W_A)$ is free over $A$ for each $i$ and lifts 
$\Fil^i(W)$. But if $j$ is such that $\Fil^j(\cDc(D))\oplus W = \cDc(D)$, then 
$\Fil^j(\cDc(D_A))\oplus W_A = \cDc(D_A)$ by Nakayama's lemma, and the result follows.
\end{pf}

Let $D$ be a crystalline $\fg$-module over $\Ro_L$. For the applications of this paper, the most 
important paraboline deformation functors of $D$ will be the ones associated to the $L[\varphi]$-stable 
filtrations $\mathcal{P}=(\Fil_i)$ of $\cDc(D)$ which are {\it complete
flags}, {\it i.e.} $\dim_L\Fil_i=i$ for $0\leq i \leq \rk_{\Ro_L}(D)$. Such a filtration is called a 
{\it refinement} of $D$, and we shall usually denote them by the letter $\Ref$. Via Lemma~\ref{trivprop}, 
refinements of $D$ corresponds to triangulations of $D$ in the sense of
Colmez~\cite{colmeztri}. Note that $D$ admits a refinement (or a
triangulation) if and only if the characteristic polynomial of $\varphi$ on
$\cDc(D)$ splits in $L$. Before stating the 
main result about the trianguline deformation functors of $D$, we need some preliminary remarks about the 
extension of Sen's theory to $\fg$-modules.

For any $\fg$-module $D$ over $\Ro$, it makes sense to consider its evaluation at a $p^n$-th root 
of unity for any sufficiently big $n$, which defines a semi-linear $\Gamma$-module $\cD_{\rm Sen}$ of rank $\rk_{\Ro}D$ over 
$\Q_p(\mu_{p^{\infty}})$ by extending the scalars (see e.g.~\cite[\S 2.2.7]{bch}). The characteristic polynomial of the Sen operator of $\cDS(D)$
will be denoted by $\PSen(D)$, its eigenvalues are called the (generalized) Hodge-Tate weights of $D$. As in the classical case, 
let us say that $D$ is Hodge-Tate if its Sen operator is semi-simple with integral eigenvalues. If $D$ 
is crystalline, it is Hodge-Tate and its Hodge-Tate weights are the jumps (with multiplicity) of its 
Hodge-filtration. 
When $D$ is a $\fg$-module over $\Ro_A$, then $\cDS(D)$ is a free
$A\otimes_{\Q_p} \Q_p(\mu_p^{\infty})$-module, hence we may also consider the
relative characteristic
polynomial $\PSen_{/A}(D)$ of the Sen operator viewed as an $A$-linear
endomorphism.

Let $D$ be a Hodge-Tate $\fg$-module over $\Ro_L$. There is a relatively representable subfunctor 
$\cX_{D,{\rm Sen}} \subset \cX_{D}$ parameterizing deformations $D_A$ which are Hodge-Tate, viewed 
as $\fg$-module over $\Ro$. For any $D$, the relative Sen polynomial induces a natural functor 
$\PSen : \cX_D \rightarrow (\fga)^d$, this latter space being identified with the completion at 
$\PSen(D)$ of the (affine) space of monic polynomial of degree $d$.  The following lemma is easy (\cite[Lemma 2.2.11, Prop. 2.3.3]{bch}).

\begin{lemma}\label{lemmesen}\begin{itemize} \item[(a)] If $0 \rightarrow D \rightarrow D' \rightarrow D'' \rightarrow 0$ is an exact
sequence of $\fg$-modules over $\Ro_A$, then $\PSen_{/A}(D')=\PSen_{/A}(D)\PSen_{/A}(D'')$.
\item[(b)] If $\delta : \Qp^* \rightarrow A^*$ is a continuous character, $\PSen_{/A}(\Ro_A(\delta))=T+\left(\frac{\partial\delta}{\delta \gamma}\right)_{|\gamma=1}$. In particular, if $D$ has rank $1$ over $\Ro_L$ then $\PSen : \cX_{D} \rightarrow \fga$ is formally smooth.
\end{itemize}
\end{lemma}

\noindent Consider the following property:
\begin{itemize}
\item[(iv)] $\PSen_{/L}(D)$ has $\rk_{\Ro_L}$ distinct roots in $\Z$.
\end{itemize}
Under (iv), $\cX_{D,{\rm Sen}}$ is simply the locus of $\cX_D$ defined by $\PSen=\PSen(D)$. When $D$ is crystalline, $\cX_{D,{\rm crys}} \subset \cX_{D,{\rm Sen}}$.

\begin{prop}\label{maindefref} Let $D$ be a crystalline $\fg$ module of rank $n$ over $\Ro_L$ 
satisfying (i), (ii), (iii) and (iv) above. Let $\Ref$ be a non-critical refinement 
of $D$. 

Then $\cX_{D,{\rm crys}}$ and $\cX_{D,\mathcal{F}}$ are pro-representable subfunctors of $\cX_{D}$, they are formally smooth over $L$ of respective dimension $\frac{n(n-1)}{2}+1$ and $\frac{n(n+1)}{2}+1$.
The morphism $\pht$ is formally smooth and $\cX_{D,{\rm crys}}=\cX_{D,{\rm
Sen}}\times_{\cX_D} \cX_{D,\mathcal{F}}$. In particular, we have an exact sequence on tangent spaces:
$$0 \rightarrow \cX_{D,{\rm crys}}(L[\varepsilon]) \rightarrow
\cX_{D,\mathcal{F}}(L[\varepsilon]) \overset{\pht}{\rightarrow}
L^n \rightarrow 0.$$
\end{prop} 

\begin{pf} This is~\cite[Thm. 2.5.1, Thm. 2.5.10]{bch}, here is a slightly different proof. The first assertion follows from Lemma~\ref{corparacrys} and 
Prop.~\ref{crysnoncrit}, as well as the inclusion $$\cX_{D,{\rm crys}}
\subset \cX':=\cX_{D,{\rm Sen}}\times_{\cX_D} \cX_{D,\mathcal{F}}.$$ The formal smoothness of $\pht$ follows from Lemma~\ref{lemmesen} and Prop.~\ref{relsmooth}. It implies the surjectivity of $\pht$ in the
exact sequence of the statement, as well as the exactness of the whole
sequence for dimensions reasons. It also implies that $\cX'$ is formally 
smooth over $L$, hence that $\cX=\cX'$ as they have the same dimension. 
\end{pf}

The other kind of paraboline deformation functor that will play a role in
the sequel is the ``miraboline'' case. Assume that $D$ is a crystalline
$\fg$-module over $\Ro_L$ satisfying (i) and let $\phi \in L^*$ 
be an eigenvalue of $\varphi$ on $\cDc(D)$. By Lemma~\ref{trivprop}, there is a unique $\fg$-submodule 
$D_\phi$ of rank $1$, as well as a unique $\fg$-module $D^\phi$ of rank $n-1$, both direct summand as 
$\Ro_L$-modules, such that $$\cDc(D_\phi)^{\varphi=\phi} \neq 0\,\, \, \, {\rm and} \,\,\,\, \cDc(D/D^\phi)^{\varphi=\phi} \neq 0.$$
Consider the filtration $\mathcal{P}_\phi$ (resp. $\mathcal{P}^\phi$) of $D$ whose unique proper subspace is $\Fil_1=D_\phi$ (resp. $\Fil_1=D^\phi$). It will be convenient to modify a little bit those functors by fixing some Hodge-Tate weight. Assume that (iv) holds and let $k_1 < k_2 < \dots < k_n$ be the Hodge-Tate weights of
$D$. We define $\cX_{D,\phi \downarrow}$ (resp. $\cX_{D,\phi \uparrow}$) as the
subfunctor of $\cX_{D,\mathcal{P}_\phi}$ (resp. $\cX_{D,\mathcal{P}^\phi}$) defined by $\pht(k_1)=0$ (resp. $\pht(k_n)=0$).

\begin{prop}\label{defmira} Assume that $D$ is a crystalline $\fg$-module of rank $n$ over $\Ro_L$ satisfying (i), (ii), (iii) and
(iv) and let $\phi$ be an eigenvalue of $\varphi$ on $\cDc(D)$. 
Then $\cX_{D,\phi \downarrow}$ is a pro-representable subfunctor of $\cX_D$, formally smooth over $L$ 
of dimension $n^2-n+1$.

Assume that $\cDc(D)^{\varphi=\phi}$ is non-critical in $\cDc(D)$. Then $\cX_{D,{\rm crys}} \subset \cX_{D,\phi \downarrow}$.
Furthermore, for any object $A$ in $\mathcal{C}$, $\cX_{D,\phi \downarrow}(A)$ is the subset of
$(D_A,\pi) \in  
\cX_{D}(A)$ such that there exists $\widetilde{\phi} \in A^*$ lifting $\phi$ and such that 
$\cDc(D_A)^{\varphi=\widetilde{\phi}}$ is free of rank $1$ over $A$.
\end{prop}

\begin{pf} By Lemma~\ref{corparacrys} (i), $\cX_{D,\mathcal{P}_\phi}$ is a subfunctor
of $\cX_D$, formally smooth over $L$ of dimension $n^2-n+2$. The assertion
on $\cX_{D,\phi \downarrow}$ follows from this, property (ii), Prop.~\ref{relsmooth} and Lemma~\ref{lemmesen}.
The inclusion $\cX_{D,{\rm crys}} \subset \cX_{D,\phi \downarrow}$ follows from Prop.\ref{crysnoncrit}. 

Before checking the last assertion, remark that for an element $(D_A,\Fil_i,\pi) \in \cX_{D,\mathcal{P}_\phi}(A)$, and $\widetilde{\phi}\in A^*$ lifting $\phi$, the natural inclusion induces a bijection 
\begin{equation}\label{devisscrys} \cDc(\Fil_1(D_A))^{\varphi=\widetilde{\phi}} = \cDc(D_A)^{\varphi=\widetilde{\phi}}.
\end{equation}
Indeed, $-\mapsto \cDc(-)^{\varphi=\widetilde{\phi}}$ is left exact over the category of $\fg$-modules over $\Ro_A$ and $\cDc(D_A/\Fil_1)^{\varphi=\widetilde{\phi}}$ vanishes by a devissage as $\cDc(D/\Fil_1(D))^{\varphi=\phi}=0$. 

Let us check the last assertion, define $\cX'(A)$ as the subset of $(D_A,\pi)
\in \cX_{D}(A)$ having the property given in the statement. Assume first that
$(D_A,\Fil_i,\pi) \in \cX_{D,\phi \downarrow}(A)$. Write $\Fil_1 \simeq \Ro_A(\delta)$ for some character $\delta : \Q_p^* \rightarrow A^*$. As $\cDc(D_\phi)$ is non-critical $\PSen_{/L}(D_{\phi})=T-k_1$, so the reduction $\overline{\delta} : \Q_p^* \rightarrow L^*$ of $\delta$ modulo $m_A$ is the algebraic character $x \mapsto x^{-k_1}$ over $\Z_p^*$ (Lemma~\ref{lemmesen} (b)). As $k_1$ is a simple root of
$\PSen_L(D)$ by (iv), the assumption $\PSen_{/A}(D_A)(k_1)=0$ and 
Lemma~\ref{lemmesen} ensure that $$\PSen_{/A}(\Fil_1(D_A))=\PSen_{/A}(\Ro_A(\delta))=T-k_1.$$
This implies that $\delta(x)=x^{-k_1}$ for any $x \in \Z_p^*$, so $\Fil_1(D_A)$ is crystalline, hence a crystalline deformation of
$D_\phi$. By (\ref{devisscrys}), this implies that $(D_A,\pi) \in \cX'(A)$ (even with
$\widetilde{\phi}=\delta(p)p^{-k_1}$). 

We now check the other inclusion $\cX'(A) \subset
\cX_{D,\phi \downarrow}(A)$. Fix $(D_A,\pi) \in \cX'(A)$ and let
$\widetilde{\phi} \in A^*$ be as in the statement. We claim first that 
$$\Fil^i(\cDc(D_A)^{\varphi=\widetilde{\phi}})=0 \,
\,\,\,\,\, {\rm for \, \, all}\, \, \, \,  i>k_1.$$ Indeed, let us argue by induction on $\dim_L A$. It holds when
$\dim_L A=1$ as $D_\phi$ is non-critical. When $\dim_L A>1$ we may find an
ideal $I \subset A$ of $L$-dimension $1$, and as $D_A$ is free over $\Ro_A$
we may consider the exact sequence $0
\rightarrow ID_A=D \rightarrow D_A \rightarrow D_{A/I} \rightarrow 0$ of
$\fg$-modules over $\Ro_L$ (and of $\Ro_A$-modules). Applying
$\cDc(-)^{\varphi=\widetilde{\phi}}$ we obtain for any $i$ an exact sequence 
$0 \rightarrow \Fil^i\cDc(D)^{\varphi=\phi} \rightarrow
\Fil^i \cDc(D_A)^{\varphi=\widetilde{\phi}} \rightarrow
\Fil^i \cDc(D_{A/I})^{\varphi=\widetilde{\phi}}$
from which we get the claim by induction. As a consequence, the crystalline $\fg$-submodule $D_A' \subset D_A$ corresponding to the filtered submodule $\cDc(D_A)^{\varphi=\widetilde{\phi}}$ is free over $\Ro_A$ by Lemma~\ref{desccrys}. The natural map
$\pi : D_A' \rightarrow D$ surjects onto $D_\phi$ as it is so after applying
$\cDc(-)$, and we are done. (For an alternative argument, see~\cite[Lemma 2.5.2]{bch}).
\end{pf}

\begin{remark}\label{dualityrem} {\rm The duality $E \mapsto E^\vee=\Hom_{\Ro_A}(E,\Ro_A)$ on $\fg$-modules $E$
over $\Ro_A$ induces an isomorphism $\iota : \cX_D \isomo \cX_{D^\vee}$. Note that $D^\vee$ satisfies (i),
(ii), (iii), or (iv) if $D$ does. If $\mathcal{P}$ is a filtration on $D$,
then there is a natural dual filtration $\mathcal{P}^\vee$ on $D^\vee$. When
$D$ is crystalline, $\mathcal{P}$ is non-critical if and only if
$\mathcal{P}^\vee$ is. The isomorphism $\iota$ induces an 
isomorphism $\cX_{D,\mathcal{P}} \isomo \cX_{D^\vee,\mathcal{P}^\vee}$, as
well as an isomorphism $\cX_{D,\phi
\uparrow} \isomo \cX_{D^\vee,\phi^{-1} \downarrow}$ (when it makes sense). 
In particular, the proposition above has
an obvious analogue for $\cX_{D,\phi \uparrow}$ whose statement is left as
an exercise to the reader. }
\end{remark}

\subsection{The main theorem} ${}^{}$\ps
\label{sectmainthmloc}

Let $D$ be a crystalline $\fg$-module of rank $n\geq 1$ over $\Ro_L$. We say
that $D$ is {\it generic} if it satisfies (i), (ii) of \S~\ref{triref}, and if {\it all of 
its $n!$ refinements are non-critical}. As we shall see below, such a $D$ necessarily satisfies (iii) and
(iv) as well. In the algebraic variety of all filtered $\varphi$-modules
of dimension $n$ over $L$ whose filtration admits $n$ jumps, the generic condition is Zariski-open 
and Zariski-dense.

Let us denote by $\tg_\ast$ the tangent space $\cX_{D,\ast}(L[\varepsilon])$.
By Lemma~\ref{corparacrys}, $\tg$ contains the $L$-subvector-space $\tg_\Ref$ for any of the $n!$ 
refinements $\Ref$ of $D$. 

\begin{thm}\label{mainlocal} If $D$ is generic then $\tg=\sum_{\Ref}
\tg_{\Ref}$. 
\end{thm}

In other words, {\it any first order deformation of a generic crystalline representation is a linear 
combination of trianguline deformations}. 

We will actually prove a more general statement, whose formulation requires
some more definitions. Let $I \subset
\{1,\dots,n\}$ be a non-empty subset and denote by $i_1 < \dots < i_r$
its elements. An element $\sigma$ in $\got{S}_n$ will be called an {\it ordered cycle} on $I$
if either $\sigma = (i_1, i_2, \dots, i_r)$ or $\sigma = (i_r, i_{r-1}, \dots,
i_1)$. By a {\it nested} sequence of intervals in $\{1,\dots,n\}$ we mean a
decreasing sequence of subsets $$I_0=\{1,\dots,n\} \supset I_1 \supset \dots \supset
I_{n-1}$$
such that for each $1 \leq i \leq n-1$, $I_i=I_{i-1} \backslash \{y_i\}$ where
$y_i \in \{ {\rm Min}(I_{i-1}), \, {\rm Max}(I_{i-1})\}$. Such a sequence is
uniquely determined by the sequence $(y_i)$ entering in its definition. For
each $1 \leq i \leq n-1$ we define as well an element $y_i^* \in \{1,\dots,n\}$
by $\{y_i,y_i^*\}=\{{\rm Min}(I_{i-1}),{\rm Max}(I_{i-1})\}$.

We will say that a sequence of $n$ permutations $\sigma_0,\sigma_1,\dots,\sigma_{n-1} \in
\got{S}_n$ is {\it nested} if $\sigma_0={\rm id}$ and if there exists a nested sequence of intervals
$(I_i)$ of $\{1,\dots,n\}$ such that $$ \forall i \in \{1,\dots,n\!-\!1\},\, \, \, \, \, \, \sigma_i = c_i
\sigma_{i-1},$$ 
where $c_i$ is the ordered cycle on $I_{i-1}$ such that $c_i(y_i^*)=y_i$
(the $y_i$ and $y_i^*$ being associated to $I_i$ as above). For $n>1$ there are exactly $2^{n-2}$ 
nested sequences in $\got{S}_n$. For instance, $id$, $(1,\, 2,\, 3)$, $(2,\,3)(1,\,2,\,3)$ and ${ id}$, $(3,2,1)$,
$(1,\,2)(3,\,2,\,1)$ are the two nested sequences in $\got{S}_3$.

We will say that a sequence of $2n-1$ permutations
$$\sigma_0,\sigma_1,\sigma_1^*,\sigma_2,\sigma_2^*,\dots,\sigma_{n-1},\sigma_{n-1}^* \in
\got{S}_n$$ is {\it weakly nested} if $\sigma_0={\rm id}$ and if there exists a nested
sequence of intervals
$(I_i)$ of $\{1,\dots,n\}$ such that $$ \forall i \in \{1,\dots,n-1\},\, \, \, \, \,
\, \sigma_i = \tau_i 
\sigma_{i-1}\, \, \, \,  {\rm and}\, \, \, \,  \sigma_i^*=c_i \sigma_{i-1},$$ 
where $c_i$ is the ordered cycle on $I_{i-1}$ such that $c_i(y_i^*)=y_i$, and
where $\tau_i \in \got{S}_n$ is any element such that $\tau_i(y_i^*)=y_i$
and such that ${\tau_i}(j)=j$ for all $j \in \{1,\dots,n\}\backslash I_{i-1}$.
For example, the nested sequence of intervals defined by $y_i=i$ for each
$i$, and the elements $\tau_i=(i, n)$ for $1\leq i \leq n-1$, define a unique weakly nested sequence 
in $\got{S}_n$.

Of course, a nested sequence defines a weakly nested sequence if we set
$\sigma_i^*=\sigma_i$, but there are much more weakly nested sequences in general.
When $\tau_i(y_i) \neq y_i^*$
for $1\leq i<n-1$ (so for most weakly nested sequences, but not for the one
in the example above), note that the elements $\sigma_i$ determine uniquely
the nested sequence of intervals entering in their definition, hence the $\sigma_i^*$
as well. 

%Moreover, if $\sigma_i,\sigma_i^*$ is a weakly nested sequence of
%permutations of $\got{S}_n$, remark that the elements $${\rm id},
%\sigma_2\sigma_1^{-1},
%\sigma_2^*\sigma_1^{-1},\dots,\sigma_{n-1}\sigma_1^{-1},\sigma_{n-1}^*\sigma_1^{-1}$$
%define a weakly nested sequence of
%permutations of $\got{S}(I_1)=\got{S}_{n-1}$, where $\{1,\dots,n-1\}$ is
%identified with $I_1$ by the unique increasing bijection. 

\medskip

Let $D$ be a crystalline $\fg$-module of rank $n$ over $\Ro_L$ satisfying (i). Any
refinement $\Ref$ of $D$ determines an ordering
$(\varphi_1,\dots,\varphi_n)$ of the eigenvalues of $\varphi$ on $\cDc(D)$
defined by $\det(\varphi_{|\Fil_i(\cDc(D))})=\prod_{j=1}^i\varphi_j$, which of
course determines $\Ref$, so we shall also write
$\Ref=(\varphi_1,\dots,\varphi_n)$. In particular $\got{S}_n$ acts on the set of refinements of $D$, via $\sigma((\varphi_i))=
(\varphi_{\sigma^{-1}(i)})$, this action being uniquely transitive. We say that a
sequence of refinements $\Ref_0,\dots,\Ref_{n-1}$ is nested if there exists
a nested sequence of permutations $\sigma_i$ of $\got{S}_n$ such that
$\Ref_i=\sigma_i(\Ref_0)$ for each $i$. We say that a
sequence of refinements
$\Ref_0,\Ref_1,\Ref_1^*,\dots,\Ref_{n-1},\Ref_{n-1}^*$ is weakly nested if there exists
a weakly nested sequence of permutations $\sigma_i,\sigma_i^*$ of $\got{S}_n$ such that
$\Ref_i=\sigma_i(\Ref_0)$ and $\Ref_i^*=\sigma_i^*(\Ref_0)$ for each
$i=1,\dots,n-1$.

\begin{thm}\label{mainlocalbis} Let $D$ be a crystalline $\fg$-module
of rank $n$ over $\Ro_L$ satisfying (i) and (ii).  If
$\Ref_0,\Ref_1,\Ref_1^*,\dots,\Ref_{n-1},\Ref_{n-1}^*$ is a
weakly nested sequence of non-critical refinements
of $D$, then $\tg=\sum_{i=0}^{n-1}\tg_{\Ref_i}$.

In particular, if $\Ref_0$,\dots,$\Ref_{n-1}$ is a nested sequence of non-critical refinements
of $D$, then $\tg=\sum_{i=0}^{n-1}\tg_{\Ref_i}$. 
\end{thm}

Let us settle first the underlined representability questions (property (iii)).

\begin{lemma}\label{indnest} Let $D$ be a crystalline $\fg$-module satisfying (i). If $D$
admits a weakly nested sequence of non-critical refinements, then (iii) and (iv)
hold. Any subquotient of a generic $\fg$-module is generic. 
\end{lemma} 

\begin{pf} If $D$ admits a non critical refinement, then (iv) holds by
definition. By (i), $\End_{\fg}(D) \subset \End_{L[\varphi]}(\cDc(D))\simeq L^n$ (a diagonal 
$L$-algebra), so (iii) for $D$ is equivalent to the fact that $D$ is not the direct sum of two $\fg$-modules over $\Ro_L$. \par
	Assume that $D=D_1\oplus D_2$ and let $k_1<k_2<\dots<k_n$ be the Hodge-Tate weights of 
$D$. This defines a partition $\{1,2,\dots,n\}=A \coprod B$ by $i \in A$ if, and only if, $k_i$ 
is a Hodge-Tate weight of $D_1$. Let $(\varphi_1,\dots,\varphi_n)$ be a non-critical refinement 
of $D$. Observe that for all $i$, $\varphi_i$ is an eigenvalue of $\varphi$ on $\cDc(D_1)$ if and 
only if $i \in A$: indeed, if $(\Fil_i(D))$ is the triangulation of $D$ associated to $\Ref$, then for 
all $j$ the jump of the Hodge filtration on the line $\cDc(D/\Fil_{j-1}(D))^{\varphi=\varphi_j}$ 
is exactly $k_j$ as $\Ref$ is non-critical, and we are done by induction on $j$. As a consequence, 
if $\Ref$ and $\sigma(\Ref)$ are both non-critical, then $\sigma(A) \subset A$. If $\sigma$ is a cycle, $A=\emptyset$ 
or $A=\{1,\dots,n\}$, which concludes the proof of the first part. The
second part follows at once. 
\end{pf}

The remaining of this section is devoted to 
the proof of Theorem~\ref{mainlocalbis}. Recall that for each non critical refinement $\Ref$ we have an inclusion 
$$\tg_{\rm crys} \subset \tg_{\Ref} \subset \tg,$$
with respective dimensions $\frac{n(n-1)}{2}+1$, $\frac{n(n+1)}{2}+1$ and $n^2+1$. So 
$\tg_{\Ref}/\tg_{\rm crys}$ is a subspace of dimension $n$ inside $\tg/\tg_{\rm crys}$, which has 
dimension $\frac{n(n+1)}{2}$. The idea of the proof is to show that such an $\Ref$ being given, the tangent space 
of a suitable paraboline deformation functor of type $(n-1,1)$ or $(1,n-1)$ is in direct sum with $\tg_\Ref$ 
{\it modulo } $\tg_{\rm crys}$, and then argue by induction on $\rk_{\Ro_L}(D)$ in the $(n-1) \times (n-1)$ 
square. Note that even if we start from an \'etale $D$, we will loose this property in the induction process, which justifies the generality adopted in~\S\ref{triref}.

In the following lemma, $D$ is a crystalline $\fg$-module over $\Ro_L$ satsfying (i),
(ii), (iii) and (iv). Fix $\Ref$ a refinement of $D$ and set 
$c=(1,\,2,\,\dots,\,n) \in \got{S}_n$. 

\begin{lemma}\label{mainsplit} If $\Ref$ and $c(\Ref)$ are non critical, and
if $\phi$ is the last element of $\Ref$, $$\tg = \tg_{\Ref} +
\tg_{\phi \downarrow}\, \, \, \, {\rm and}\, \, \, \tg_{\Ref}\cap \tg_{\phi
\downarrow}=\tg_{\rm crys}.$$
Similarly, if $\Ref$ and $c^{-1}(\Ref)$ are non critical, and if $\phi$ is
the first element of $\Ref$, then $\tg=\tg_{\Ref}+\tg_{\phi \uparrow}$ and $\tg_{\Ref}\cap \tg_{\phi
\uparrow}=\tg_{\rm crys}$.
\end{lemma}

\begin{pf} The second part follows from the first one by duality by Remark~\ref{dualityrem}, so we focus on the first part. First of all, 
by Prop.~\ref{maindefref} and~\ref{defmira}, $\tg_{\rm crys} \subset \tg_{\Ref}\cap
\tg_{\phi\downarrow}$ and 
\begin{equation}\dim(\tg_{\Ref}/\tg_{\rm crys})+\dim(\tg_{\phi \downarrow}/\tg_{\rm
crys})=\frac{n(n+1)}{2}=\dim(\tg/\tg_{\rm crys}).
\end{equation}
(Prop.~\ref{maindefref} applies as $\cDc(D)^{\varphi=\phi}$ is non-critical,
since $\phi$ is the first eigenvalue of the non critical refinement
$c(\Ref)$.) Thus it only remains to show that
$\tg_{\Ref}\cap \tg_{\phi\downarrow} \subset \tg_{\rm crys}$. As $\Ref$
is non-critical, it is enough to show that 
\begin{equation} \cX_{D,\Ref}\cap \cX_{D,\phi\downarrow} \subset \cX_{D,{\rm
Sen}}
\end{equation}
by Prop.~\ref{maindefref}. Fix an object $A$ of $\mathcal{C}$, let 
$(D_A,\Fil_i,\pi) \in \cX_{D,\Ref}(A)$, and assume that $(D_A,\pi) \in \cX_{D,\phi}(A)$
as well. By property (iv), we may write $$\PSen_{/A}(D)=\prod_{i=1}^n(T-\kappa_i) \in A[T],$$ where $\kappa_i \in A$ lifts $k_i$. Choose $\widetilde{\phi} \in A^*$ lifting $\phi$ such that 
$\cDc(D_A)^{\varphi=\widetilde{\phi}}$ is free of rank $1$ over $A$ and let
$0 \leq j \leq n-1$. As $\phi$ is the last element of $\Ref$, $\cDc(\Fil_j(D))^{\varphi=\phi}=0$
and an immediate d\'evissage shows that
$\cDc(\Fil_j(D_A))^{\varphi=\widetilde{\phi}}=0$. Applying the left exact functor $\cDc(-)^{\varphi=\widetilde{\phi}}$ to the exact sequence $0 \rightarrow \Fil_j(D_A) \rightarrow D_A
\rightarrow D_A/\Fil_j(D_A) \rightarrow 0$ of $\fg$-modules over $\Ro_A$, we obtain an
$A[\varphi]$-linear injection 
\begin{equation}\label{injquot} \cDc(D_A)^{\varphi=\widetilde{\phi}} \hookrightarrow 
\cDc(D_A/\Fil_j(D_A))^{\varphi=\widetilde{\phi}}.
\end{equation}
On the other hand, as $\phi$ is an eigenvalue of multiplicity $1$ in
$\cDc(D)$, another immediate d\'evissage shows that the length of the
$A$-module $\cDc(D_A/\Fil_j(D_A))^{\varphi=\widetilde{\phi}}$ is $\leq {\rm length}(A)$, thus (\ref{injquot}) is an isomorphism. As $c(\Ref)$ is non critical, the Hodge-Tate weights of $D/\Fil_j(D)$ are $k_{j+1},\, \dots,\, k_n$ and $\cDc(D/\Fil_j(D))^{\varphi=\phi}$ is non-critical in $\cDc(D/\Fil_j(D))$. Thus Prop.~\ref{defmira} applies, and shows that for the obvious $\pi'$ we have $(D_A/\Fil_j(D_A),\pi') \in \cX_{D/\Fil_j(D),\phi\downarrow}(A)$. In particular, $$\PSen_{/A}(D_A/\Fil_j(D_A))(k_{j+1})=0.$$ As $\PSen_{/A}(D_A/\Fil_j(D_A))=\prod_{i=j+1}^n(T-\kappa_i)$ by Lemma~\ref{lemmesen}, we conclude that $\kappa_{j+1}=k_{j+1}$ for all $0 \leq j \leq n-1$: $D_A$ is Hodge-Tate, and we are done.\end{pf}

We now prove Theorem~\ref{mainlocalbis} by induction on $n=\rk_{\Ro_L}D$. When $n=1$ the
theorem is obvious as $\cX_{D,\Ref}=\cX_{D}$, so assume $n>1$. 
Let $\Ref_0, \Ref_1, \Ref_1^*,\dots, \Ref_{n-1},\Ref_{n-1}^*$ be a
weakly nested sequence of non-critical
refinements of $D$. Applying the above lemma to $\Ref=\Ref_0$,
we obtain a fortiori
\begin{equation}\label{fin1} \tg=\tg_{\Ref_0}+\tg_{\mathcal P}\end{equation}
where either $\mathcal{P}={\mathcal{P}_\phi}$ and $\phi$ is the last
element of $\Ref_0$ (case $\sigma_1^*=c$), or $\mathcal{P}=\mathcal{P}^\phi$ and $\phi$ is the first element of $\Ref_0$ (case
$\sigma_1^*=c^{-1}$). Set $D'=D/D_\phi$ in the first case and $D'=D^\phi$ in
the second case. By definition, for $i=1,\dots,n-1$ the refinements $\Ref_i$ and $\Ref_i^*$ 
have the form $(\phi,\Ref'_i)$ and $(\phi,{\Ref'_i}^*)$ (resp. $(\Ref'_i,\phi)$
and $({\Ref'_i}^*,\phi)$) in the first case
(resp. second case). Moreover, 
$$\Ref'_1,\Ref'_2,{\Ref'_2}^*,\dots,\Ref'_{n-1},{\Ref'_{n-1}}^*$$
is a weakly nested sequence of refinements of $D'$. Each $\Ref'_i$ (resp. 
${\Ref'_i}^*$) is non-critical for $D'$ as $\Ref_i$ (resp. $\Ref_i^*$) is non critical for
$D$. As $D'$ obviously still satisfies (i) and (ii), we obtain by induction:
\begin{equation}\label{fin2} \tg'=\sum_{i=1}^{n-1}
\tg'_{\Ref'_i},\end{equation}
where $\tg'_\ast=\cX_{D',\ast}(L[\varepsilon])$. Note that for
$i=1,\dots,n-1$ we have $\cX_{D,\Ref_i} \subset \cX_{D,\mathcal{P}}$.
By Prop.\ref{relsmooth} and property (ii), the natural map
$\cX_{D,\mathcal{P}} \rightarrow \cX_{D'}$ induces isomorphisms
\begin{equation}\label{fin3}\cX_{D,\Ref_i} = \cX_{D,\mathcal{P}} \times_{\cX_{D'}}
\cX_{D',\Ref'_i}.
\end{equation}
By (\ref{fin2}) and
(\ref{fin3}) we have $\tg_{\mathcal{P}} =  \sum_{i=1}^{n-1} \tg_{\Ref_i}$, and we are done by (\ref{fin1}). $\square$

Another consequence of the proof of Lemma~\ref{mainsplit} is the following proposition.

\begin{prop}\label{transversality} ({\it Transversality of trianguline deformation functors}) Assume that $D$ satisfies 
(i), (ii) and let $\Ref$ and $\Ref'$ be two refinements of $D$. If $\Ref'$ starts with the last 
element of $\Ref$, and if
$\Ref$, $c(\Ref)$ and $\Ref'$ are non-critical, $\cX_{D,\Ref}\cap \cX_{D,\Ref'}$ is exactly the subfunctor of
deformations of $D$ which are crystalline up to a twist. 
\end{prop}

%\begin{remark} We leave as an exercice to the reader to check that for $D$ generic of rank $>2$, if $\Ref_0$ is any refinement of
%$D$ and if $\Ref_i:=(i,\, i+1)(\Ref_{i-1})$ for $i\geq 1$, then $\sum_{i=0}^{n-1} \tg_{\Ref_i} \subsetneq
%\tg$. 
%\end{remark}

\subsection{Some examples}\label{examplelocalmain}

Note that any $\fg$-module of rank $2$ satisfying (i), (ii), (iii) and (iv) is generic. The following other example plays some role in the main global result
of this paper.

\begin{example}\label{dim3well} Assume that $D=\D_{\rm rig}(V)$ for some representation $V$ of dimension $3$ over $L$ satisfying
(i), (iii) and (iv). Then $D$ admits a nested sequence of non-critical refinements.
\end{example}

\begin{pf} Let $X$ be the set of Frobenius eigenvalues of $\varphi$ on $\Dc(V)$, and let 
$k_1 < k_2 < k_3$ be the Hodge-Tate weights of $V$. A line $L \subset \Dc(V)$ (resp. a plane $P$) is critical if $L \subset \Fil^{k_2}(\Dc(V))$ (resp. $P \supset \Fil^{k_3}(\Dc(V))$).

Assume first that there is some $\phi \in X$ such that the line $L_\phi:=\Dc(V)^{\varphi=\phi}$ is critical. Then we claim that there is a unique such $\phi$ in $X$ and that no $\varphi$-stable plane 
$P \subset \Dc(V)$ is critical. In particular, the critical refinements of $D$ are exactly the two ones starting with $\phi$. This implies the lemma in this case as each refinement
of the nested sequence
$$(\phi',\phi,\phi''), \, \, \, (\phi'',\phi',\phi), \, \, \,
(\phi'',\phi,\phi')$$
are non critical. If there is no critical line $L_\phi$ as above, then either $D$ is generic or there is a $\varphi$-stable critical plane $L_{\phi'}\oplus L_{\phi''}$. It follows from the preceding claim applied to $V^*$ that in the latter case such a plane is 
unique, hence the critical refinements of $D$ are the ones ending by $\phi$ and 
$$(\phi',\phi,\phi''), \, \, \, (\phi,\phi'',\phi'), \, \, \,
(\phi'',\phi,\phi')$$
is a nested sequence of non-critical refinements, and we are done again.
%The claim is a simple exercise using the weak admissibility of $D$ and the irreducibility of $V$, that we leave as an exercise to the reader.

Let us check the claim. As $V$ is indecomposable, $\Dc(V)$ is not the direct sum of an admissible line and of an admissible plane. Write $X=\{\phi,\phi',\phi''\}$ and do not 
assume anything on $\phi$ for the moment. If $L_\phi\oplus L_{\phi'}={\rm Fil}^{k_2}(\Dc(V))$, the weak admissibility property of 
$\Dc(V)$ implies that $L_{\phi''}$ and $L_\phi\oplus L_{\phi'}$ are admissible, which is absurd. Similarily, 
$L_\phi \neq {\rm Fil}^{k_3}(\Dc(V))$ otherwise $L_\phi$ and $L_{\phi'}\oplus L_{\phi''}$ would be admissible. Assume now that 
$L_\phi$ is critical. If $P$ is a $\varphi$-stable critical plane, then $P \neq {\rm Fil}^{k_2}(\Dc(V))$ and $P$ does not contain 
$L_\phi$ by what we just proved, so $P=L_{\phi'}\oplus L_{\phi''}$ and the jump indices of the 
Hodge filtration on $P$ are $k_1$ and $k_3$. But weak admissibility implies that $L_\phi$ and $P$
are admissible, and we are done.\end{pf}

\begin{example}\label{casexsym} Let $D$ be a crystalline $\fg$-module of rank $2$ over $\Ro_L$ satisfying (i), (iii) and such that the ratio $r$ between its two Frobenius eigenvalues satisfies $r^i \neq 1$ for $i=1,\dots,n$. 

Then ${\rm Symm^n_{\Ro_L}(D)}$ satisfies (i) and all of its refinements are non critical.
\end{example}

\begin{pf} We leave as an exercise to the reader to check that the statement follows from the following claim: if $I \subset \{0,\dots,n\}$ has $i$ elements and if we have a polynomial identity $Q(x)(1+x)^i = \sum_{i \in I} a_i x^i$ where $a_i \in L$ and $Q \in L[x]$ has degree $\leq n+1-i$, then $Q=0$. To check the claim, consider the $i$ successive derivatives at $-1$ of the right hand side. This is a linear system in the unknown $(a_i)_{i \in I}$ whose determinant (in absolute value) is the Vandermonde determinant over the $i$ elements of $I$. 
\end{pf}

Let us end this paragraph by introducing some {\it regularity conditions} that will arise for technical reasons in the subsequent global applications. Let $D$ be a crystalline $\fg$-module over $\Ro_L$ satisfying (i). We say that a refinement $\Ref=(\varphi_1,\dots,\varphi_n)$ is {\it regular} if for
all $1\leq j \leq n$, the element $\varphi_1\varphi_2\cdots\varphi_j$
is a simple root of the polynomial $\prod_{I}(T-\prod_{i \in I}
\varphi_i)$, the product being over the $I \subset \{1,\dots,n\}$ with $|I|=j$.\footnote{When $n\leq 3$, this condition is a consequence of property (i) of $D$.}

\begin{definition}\label{defwgr}{\it We say that $D$ is {\it weakly generic and regular} if it satisfies (i), (ii), and if it possesses a weakly nested sequence of non-critical refinements $\{\Ref_i, \Ref_i^*\}$ such that $\Ref_i$ is regular for each $i=0,\dots,n-1$. If $D=\Dr(V)$, we say that $V$ is regular if $D$ is.} \end{definition}

The condition of being weakly generic and regular is Zariski-open and Zariski-dense in the algebraic variety of all filtered $\varphi$-modules of dimension $n$ over $L$ (whose filtration admits $n$ jumps). This is actually still true if we restrict to filtered $\varphi$-modules satisfying a self-duality condition (orthogonal or symplectic, possibly with a similitude factor, in which cases the filtration is fixed by a Lagrangian in $L^n$). Indeed, set $k=[n/2]$, let $(X_i)_{i=1,\dots,k}$ be some indeterminates and consider the sequence $s=(X_1,\dots,X_k,X_k^{-1},\dots,X_1^{-1})$ if $n=2k$, or $s=(X_1,\dots,X_k,1,X_k^{-1},\dots,X_1^{-1})$ if $n=2k+1$. Observe that $s$ is regular (for the obvious definition), moreover each permutation of $s$ whose first $k$ terms do not contain both some $X_i$ and its inverse $X_i^{-1}$ is regular as well (for instance, there are $2^k k!$ such permutations when $n$ is even). We conclude as we may find a weakly nested sequence of permutations of $\got{S}_n$, associated to the nested sequence of intervals $y_i=i$, whose $\sigma_i$ have the property above: for instance take $\tau_i = (i, n-i, n)$ for all $i\leq k$ if $n$ is even. 

\begin{example}\label{symtrois} Let $D$ be generic of rank $4$ over $\Ro_L$ and assume that it admits a refinement of the form $\Ref=(\mu,\mu x,\mu x^2,\mu x^3)$ for some $\mu,x\in L^*$. If $x^j\neq 1$ for $1\leq j \leq 4$, then $D$ is weakly generic regular.
\end{example}

\begin{pf} $\Ref_0$, $\Ref_1=(\mu x^3, \mu x, \mu, \mu x^2)$, $\Ref_2=(\mu x^3, \mu x^2, \mu x, \mu)$ and $\Ref_3=(\mu x^3, \mu x^2, \mu, \mu x)$ is a weakly nested sequence of regular refinements of $D$ (associated to $y_i=i$). 
\end{pf}

\section{The eigenvariety at non-critical tempered classical points}

Our main goal in this section is to show that the eigenvarieties of $U(3)$ or even $U(n)$ are \'etale over 
the weight space at the (stable, tempered) non-critical classical points. 

%enonce un corollaire sur les familles

\subsection{An infinitesimal classicity criterion}\label{infclasscrit}

Let $n\geq 1$ be any integer and let $U/F$ be a unitary group in $n$ variables attached to $E/F$. 
We need to recall some definitions of the theory of $p$-adic automorphic forms for $U$. The reader 
may consult~\cite[\S 2]{che2} for a detailled discussion and complete proofs
of the statements of this section, as well as~\cite{che} and~\cite[\S 7.3]{bch}. 
We assume that $U(F_v)$ is compact for all archimedean places $v$ and that $U(F_v)\simeq \GL_n(\Q_p)$ if $v$ 
divides $p$. In particular $p$ splits in $E$ and for each
$v$ in the set $S_p$ of finite places of $F$ dividing $p$, we fix such an
isomorphism, as well as a place $\widetilde{v}$ of $E$ above $p$ (so
$\Q_p=F_v=E_{\widetilde{v}}$).  

We need to introduce some group theoretic notations concerning
$\GL_n(\Q_p)^{S_p}$, that we shall view as an algebraic group over $\Q_p$. Let $B$ be its 
upper triangular Borel subgroup, $T$ its diagonal torus, $X^*(T)$ the algebraic characters
of $T$, $I$ the subgroup of 
$\GL_n(\Z_p)^{S_p}$ whose elements are upper triangular modulo $p$, $\Nob$ the subgroup of lower 
triangular elements of $I$, $T^0=T\cap I$, $T^- \subset T$ the submonoid whose elements $t$ satisfy 
$t^{-1}\Nob t \subset \Nob$, and $M=IT^{-}I$ (it is a submonoid of $\GL_n(\Q_p)^{S_p}$).
% and that there is a unique multiplicative homomorphism $\tau : M \rightarrow T^-/T^0$ such that $m \in I\tau(m)I$ for all $m \in M$. 
%In particular, we may and will see characters of $T/T^0$ as characters of $M$ by composition with $\tau$. 

For each $B$-dominant weight $\chi \in X^*(T)$, let $W_\chi$ be the irreducible algebraic 
$\Q_p$-representation of $\GL_n(\Q_p)^{S_p}$ with highest weight $\chi$. For any $\chi \in X^*(T)$, 
let $\Cc_{\chi}$ be the standard analytic principal series of the Iwahori subgroup $I$. It may be 
defined as follows. If $J$ a finite set, a function $f : \Z_p^J \rightarrow \Q_p$ is said analytic 
if $f((x_j))$ belongs to the Tate algebra $\Q_p\langle (x_j)_j \rangle$. Let $J$ be the set of 
triples $(i,j,v)$ with $1\leq i < j \leq n$ and $v \in S_p$. Identifying $\Z_p^J$ with $\Nob$ 
via the bijection $(x_{(i,j,v)}) \mapsto (px_{i,j})_v$, we obtain a notion of analytic function 
on $\Nob$ and set:

$$\Cc_\chi = \left\{ \begin{array}{lll}
f: IB \longrightarrow \Q_p,\, \, f(xb)=\chi(b)f(x)\, \, \, \, \forall\, \,  x \in IB, \, \, b \in B, \\ \\
f_{|\Nob}\, {\rm\, is\,}\, {\rm analytic}.
\end{array} \right\} .$$

The product in $\GL_n(\Q_p)^{S_p}$ induces an isomorphism $\Nob \times B \isomo I\,B$. 
Moreover, as $t^{-1}\Nob t \subset \Nob$ for $t \in T^-$, we have $M^{-1}IB \subset IB$, 
thus left translations $(m.f)(x):=f(m^{-1}x)$ defines a representation of $M$ on $\Cc_\chi$. 
When $\chi$ is dominant, and if $v$ is a highest weight vector in $W_\chi$,
the map $\varphi \mapsto (g\mapsto \varphi(g(v)))$ defines a $\Q_p[M]$-equivariant
embedding $W_\chi^\vee \rightarrow \Cc_\chi$, whose image is the 
subspace of polynomial elements in $\Cc_\chi$.\ps

The spaces of $p$-adic automorphic forms for $U$ are built from these representations $\Cc_\chi$. 
Let $K$ be any compact open subgroup of $U(\AAA_{F,f})$ of the form $K=I \times K^{S_p}$. For any 
$\Q_p[M]$-module $W$, consider the $\Q_p$-vector space
$$F(W) := \left\{ \begin{array}{lll} f: \U(F)\backslash \U(\AAA_{F,f}) \longrightarrow W,\\ \\
f(gk)=(\prod_{v|p} k_v)^{-1}f(g),\, \, \forall g\in \U(\AAA_{F,f}),\, \, \forall\, k
\in K. \end{array}\right\}$$
This defines an exact functor from the $\Q_p[M]$-modules to the $\Q_p$-vector spaces. Moreover, 
$F(W)$ is in a natural way a module over the Atkin-Lehner algebra $\ATL^-$ which is the $\Q$-subalgebra 
of the Iwahori Hecke-algebra of $\GL_n(\Q_p)^{S_p}$ relative to $I$, consisting of functions with support 
in $M$. Recall that the map $T^-/T^0 \rightarrow \ATL^-$, $t \mapsto \mathbbm{1}_{ItI}$, is a 
multiplicative homomorphism inducing an isomorphism $$\Q[T^{-}/T^0] \isomo
\ATL^-$$ (in particular 
$\ATL^-$ is commutative). Using this latter isomorphism we will view elements of $T^-$ as elements of 
$\ATL^-$. Last but not least, $F(W)$ also admits a natural structure of module over the Hecke algebra 
of $(U(\AAA_{F,f}^{S_p}), K^{S_p})$ that commutes to $\ATL^-$. We fix a commutative subring $\HH$ of 
this latter Hecke-algebra that contains the spherical Hecke-algebra for almost all finite primes $v$ 
of $F$.

Let $\chi \in X^*(T)$. The space of $p$-adic analytic automorphic forms of $U$ of weight $\chi$ and 
level $K$ is the space $F(\Cc_\chi)$. It is a $\Q_p$-Banach space in a natural way, on which the 
Hecke-operators $\ATL^- \otimes \HH$ act continuously. Moreover, if $T^{--} \subset T^-$ denotes 
the submonoid whose elements $t$ are such that $t^{-1}\Nob t$ vanishes mod $p$, then any $t \in T^{--}$ 
acts compactly on $F(\Cc_\chi)$. 

When $\chi$ is dominant, $F(\Cc_\chi)$ contains as $\ATL^- \otimes \HH$-submodule the subspace 
$F(W_\chi^\vee)$ of automorphic forms of weight $\chi$ and level $K$, often called the subspace of 
{\it classical} $p$-adic automorphic forms. If we fix a $\iota$ as in
\S\ref{sectmod}, this subspace is a canonical 
$\Q_p$-structure of the complex $\Q[U(\AAA_{F,f})//K]$-module 
\begin{equation}\label{isomclass}\bigoplus_\Pi m(\Pi) \Pi_f^K,\end{equation}
where $\Pi$ varies in the set of irreducible automorphic representations of $U(\AAA_F)$ such that 
$\Pi_\infty \simeq W_\chi$, and with multiplicity $m(\Pi)$ in $L^2(U(F)\backslash \U(\AAA_F))$.

We shall also need to assume that the global base-change from $U$ to $\GL_n(\AAA_E)$ is known, 
or rather its corollary that for any automorphic representation $\Pi$ of $U$ there is a semi-simple 
continuous Galois representation $\rho_\Pi : G_E \longrightarrow \GL_n(\Qpb)$ which is
$\iota$-compatible with 
the unramified, Frobenius semi-simplified, local Langlands correspondence at the places where
$\Pi$, $F$ and $U$ are unramified (split places would even be enough), and normalized so that 
$\rho_\Pi^{\vee,c}\simeq \rho_\Pi(n-1)$ (note that $n$ is not necessarily odd in this paragraph). 
This is currently known in the following cases :  
\begin{itemize}
\item[-] $n\leq 3$ (Rogawski) or, 
\item[-] there is a finite place $v$ such that $U(F_v)$ is the group of invertible elements 
of a central division algebra over $F_v$ (Clozel-Labesse, Harris-Labesse) or, 
\item[-] $[F:\Q]\geq 2$ (Labesse).
\end{itemize}

Before stating the {\it infinitesimal classicality criterion} we need to recall
the relation between refinements and eigenvectors of $\ATL^-$. Let $f \neq 0 \in
F(W_\chi^\vee)\otimes_{\Q_p} \Qpb$ be a 
common eigenvector for all the elements of $\ATL^- \otimes \HH$. Let $\Pi$ be any automorphic 
representation of $U$ containing $f$, the Galois representation
$\rho_f:=\rho_\Pi$ does not depend on the choice of $\Pi$ by Cebotarev's theorem. 
According to the receipe described in \cite[\S 6.4]{bch}, the action of $\ATL^-$ on $f$ determines 
for each $v \in S_p$ a canonical ordering $(\varphi_{i,v})_{i=1}^n$ of the Langlands conjugacy class associated to 
the representation $\Pi_v$ (which has Iwahori invariants by construction). 
When $\varphi_{i,v}\varphi_{j,v}^{-1} \neq p$ for all $i \neq j$ then $\Pi_v$ is unramified, 
so $\rho_{f,\widetilde{v}}$ is crystalline. Moreover, the characteristic polynomial of its
crystalline Frobenius coincides via $\iota$ with the characteristic
polynomial of Langlands' conjugacy class of $\Pi_v$ (\cite{chharris}). If furthermore $\varphi_{i,v}\varphi_{j,v}^{-1} \neq 1$ for all $i \neq j$, then $(\varphi_{i,v})_{i=1}^n$ defines a refinement of
$\rho_{f,\widetilde{v}}$ 
in the sense of~\S\ref{triref} (see the discussion preceding Thm~\ref{mainlocalbis}).

%Let us say that a sequence $(\varphi_1,\dots,\varphi_n)$, with $\varphi_i$ in
%$\Qpb^*$, is {\it regular}
%if for all $1\leq j \leq n$, the element $\varphi_1\varphi_2\cdots\varphi_j$
%is a simple root of the polynomial $\prod_{I}(x-\prod_{i \in I} \varphi_i)$,
%the product being over the $I \subset \{1,\dots,n\}$ with $|I|=j$. In
%particular, the $\varphi_i$ are distinct (and this is sufficient if $n\leq
%3$). 

\begin{prop}\label{controlethm} Let $f \in F(W_\chi^\vee) \otimes_{\Q_p} \Qpb$ be an eigenvector for all the 
elements of $\ATL^- \otimes \HH$. Assume that for each $v \in S_p$, the sequence $(\varphi_{i,v})_{i=1}^n$ defined above satisfies $\varphi_{i,v}\varphi_{j,v}^{-1}\neq 1,p$ for any $i\neq j$, and defines a non-critical, regular, refinement of $\rho_{f,\widetilde{v}}$. 

Then the generalized $\ATL^-\otimes \HH$-eigenspace of $f$ inside
$F(\Cc_\chi)\otimes_{\Q_p}\Qpb$ is included in
$F(W_\chi^\vee)\otimes_{\Qp} \Qpb$.
\end{prop}

When $(\varphi_{i,v})_{i=1}^n$ is numerically non-critical, this result follows from the standard classicity 
criterion (compare~\cite[Rem. 2.4.6]{bch} with~\cite[Thm. 1.6 (vi)]{che2}). However, this simple case is not 
enough for our purpose in this paper, and the general case is much deeper. It is the analogue of the 
classical fact that an ordinary modular eigenform which is in the image of the
theta map is split at $p$ 
(Mazur-Wiles).

The first ingredient of the proof is the following result of Jones, which is a (locally) analytic version 
of the Bernstein-Gelfand-Gelfand resolution, and plays the role of the theta map in the context of $p$-adic 
modular forms. 

\begin{lemma}{\rm (\cite[Thm. 22]{jones})} Assume that $\chi \in X^*(T)$ is dominant. There exists an exact sequence 
of $\Q_p[M]$-modules : 
$$ 0 \longrightarrow W_{\chi}^\vee \longrightarrow \Cc_{\chi} \longrightarrow
\prod_{l(\sigma)=1} \Cc_{\sigma(\chi+\rho)-\rho}.$$
\end{lemma}

In this statement, $\sigma$ is an element in the Weyl group
$\got{S}=\got{S}_n^{S_p}$ of $T$, 
$\rho$ is the half-sum of the positive roots (this is not quite an element of $X^*(T)$, 
but $\sigma(\chi+\rho)-\rho$ is anyway). By exactness of the functor $F$ we obtain an exact 
sequence of $\ATL^-\otimes\HH$-modules $0 \rightarrow F(W_{\chi}^\vee) \rightarrow F(\Cc_{\chi}) 
\rightarrow \prod_{l(\sigma)=1} F(\Cc_{\sigma(\chi+\rho)-\rho})$. Let $\psi : \ATL^-\otimes \HH 
\rightarrow \Qpb$ be the $\Q$-algebra homomorphism defined by $f$: $a(f)=\psi(a)f$ for all $a \in 
\ATL^-\otimes \HH$. If $E$ is any $\Q_p\otimes\ATL^- \otimes \HH$-module we denote by
$E[\psi]$ the generalized 
eigenspace associated to $\psi$: $$E[\psi]=\{ e\in E\otimes_{\Q_p}\Qpb, \, \, \exists n \in \N, \, 
\, \forall a \in \ATL^- \otimes \HH, \, \,\,\, (a-\psi(a){\rm id})^ne=0\}.$$ As $f$ is classical, 
$\psi(T^-) \subset \Qpb^*$, and by compactness of any $t \in T^{--}$ we obtain 
an exact sequence of finite dimensional $\Qpb$-vector spaces :
$$0 \longrightarrow F(W_{\chi}^\vee)[\psi] \longrightarrow F(\Cc_{\chi})[\psi] \longrightarrow 
\prod_{l(\sigma)=1} F(\Cc_{\sigma(\chi+\rho)-\rho})[\psi].$$
We claim that for each simple reflexion $\sigma=(i, i+1)_v$, we have
$F(\Cc_{\sigma(\chi+\rho)-\rho})[\psi]=0$. 
Indeed, it is enough to show that there is no $g\neq 0 \in
F(\Cc_{\sigma(\chi+\rho)-\rho})[\psi]$ that is an 
eigenvector for $\ATL^-\otimes \HH$. Assume by absurdum that there are such $i$, $v \in S_p$ and $g$.

As $g$ and $f$ share the same system of Hecke-eigenvalues under $\ATL^-\otimes \HH$, the form $g$ is of finite slope and its associated $p$-adic Galois representation $\rho_g : \G_E \rightarrow \GL_n(\Qpb)$ coincides with $\rho_f$.
%,as well as the same associated refinements 
%$(\varphi_{i,v})_{i=1}^n$ at the places $v \in S_p$. 
However, the form $f$ has weight 
$\chi$ whereas its ``companion form'' $g$ has weight $\sigma(\chi+\rho)-\rho$. Let us denote by 
$$k_{1,v} < k_{2,v} < \cdots < k_{n,v}$$
the Hodge-Tate weights of $\rho_{f,\widetilde{v}}$. They are related to the weight $\chi$ of $f$ by the usual receipe, namely an identification $X^*(T) \simeq (\Z^n)^{S_p}$ making the dominant elements of $X^*(T)$ correspond to increasing sequences as above. According to the same receipe, the ordering of the Hodge-Tate numbers of 
$\rho_{g,\widetilde{v}}$ corresponding to $\sigma(\chi+\rho)-\rho$ is the same one
except that $k_{i,v}$ and $k_{i+1,v}$ are interchanged. As the actions of $\ATL^-$ on $g$ and $f$ coincide, Kisin's result \cite[\S 5]{kisin} on
the continuity of crystalline periods in refined families ensures that $\rho_g$ has the following property\footnote{To deduce this result from the litterature, we may refer as follows to~\cite{bch}. First, choose by Lemma~7.8.11 an affinoid neighbourghood $\Omega$ of the point $x$ corresponding to $g$ in the eigenvariety of $U$ of level $K$. There is a surjective alteration $\Omega' \rightarrow \Omega$ and a locally free $\OO_{\Omega}$-module $\mathcal{M}$ with a continuous $\OO_{\Omega'}$-linear action of $\G_E$ which is generically semi-simple and whose trace is the pull-back of the natural family of pseudo-characters on $\Omega$. Apply Theorem 3.3.3 to $\Lambda^i_{\OO_{\Omega'}}(\mathcal{M})$ and any point $x' \in \Omega'$ above $x$, its assumptions are satisfied by~Prop. 7.5.13 and~\S 4.2.4. We obtain (\ref{jump}) with $\overline{\mathcal{M}}_{x'}$ instead of $\rho_{g,\widetilde{v}}$. But $\rho_{g,\widetilde{v}}$ is the semi-simplification of $\overline{\mathcal{M}}_{x'} \otimes \Qpb$, so (\ref{jump}) follows from the left-exactness of the functor $\Fil^k(\Dc(-)^{\varphi=\phi})$.}:
\begin{equation}\label{jump} {\rm Fil}^{k_{1,v}+k_{2,v}+\cdots+k_{i-1,v}+k_{i+1,v}}(D_{\rm crys}(\Lambda^i
\rho_{g,\widetilde{v}})^{\varphi=\varphi_{1,v}\varphi_{2,v}\cdots\varphi_{i,v}}) \neq 0.
\end{equation}

\noindent As $(\varphi_{i,v})_{i=1}^n$ is regular, if 
$F_i=\oplus_{j=1}^i\Dc(\rho_{g,\widetilde{v}})^{\varphi=\varphi_{j,v}}
\subset \Dc(\rho_{g,\widetilde{v}})$ we
have
$$D_{\rm crys}(\Lambda^i
\rho_{g,\widetilde{v}})^{\varphi=\varphi_{1,v}\varphi_{2,v}\cdots\varphi_{i,v} }= \Lambda^i
(F_i) \subset \Lambda^i D_{\rm crys}(\rho_{g,\widetilde{v}}).$$ 
On the other hand, $(\varphi_{i,v})_{i=1}^n$ is non critical so the jumps
of the filtration on $F_i$ induced by the Hodge filtration on 
$\Dc(\rho_{g,\widetilde{v}})$ are $k_{1,v},k_{2,v},\dots,k_{i,v}$, which
contradicts (\ref{jump}). $\square$

\subsection{Non critical classical points are \'etale over the weight space}\label{noncritetalepar} The preceding result, combined with some properties of automorphic forms, 
allows to show that eigenvarieties are \'etale over the weight space at non critical classical points. \ps

Let $\ell \neq p$ be a prime and $E_\ell$ a finite extension of $\Q_\ell$. We denote by $W_{E_\ell}$ the Weil-group of $E_\ell$ and 
$I_{E_\ell} \subset W_{E_\ell}$ its inertia subgroup. If $\rho : W_{E_\ell} \rightarrow \GL_n(\Qpb)$ is continuous, 
it admits an associated Frobenius semi-simple Weil-Deligne representations $(r,N)$. Recall that we have fixed 
some embeddings $\iota_p, \iota_\infty$. We say that $\rho$ is defined over $\Qb$ if $r$ is. In this case, 
the local Langlands correspondence associates to $(r,N)$ (and $\iota_p$, $\iota_\infty$) a canonical irreducible 
smooth complex representation of $\GL_n(E_\ell)$. For two continuous representations 
$\rho, \rho' : W_{E_\ell} \longrightarrow \GL_n(\Qpb)$, with associated Weil-Deligne representations $(r,N)$ and 
$(r',N')$, we write $\rho \prec_I \rho'$ if $r_{|I_{E_\ell}} \simeq r'_{|I_{E_\ell}}$ and if {\it $\rho$ has less 
monodromy than $\rho'$} in the following sense : for each irreducible $\Qpb$-representation $\tau$ of $I_{E_\ell}$, 
the nilpotent conjugacy class of the monodromy operator is greater on the $\tau$-isotypic component of $r'$ than 
on the $\tau$-isotypic component of $r$ (for the dominance ordering, the conjugacy class of $0$ being  by convention 
the smallest one). Recall that if $r=\sum_{i=1}^s r_i$ and $r'=\sum_{j=1}^{s'}{r'_j}$ are irreducible decompositions 
of $r$ and $r'$, then $r_{|I_{E_\ell}} \simeq r'_{|I_{E_\ell}}$ if and only if, $s=s'$ and, up to renumbering the $r_i$, 
$r'_i$ is an unramified twist of $r_i$ for each $i \in \{1,\dots,s\}$ (see e.g.\cite[Lemme 3.14 (i)]{che2}).

\begin{lemma}\label{lemmehorsp} Let $X$ be an irreducible affinoid over $\Q_p$, $Z \subset X(\Qpb)$ a Zariski-dense 
subset, and let $\rho : W_{E_\ell} \longrightarrow \GL_n(\OO(X))$ be a continuous representation.
\begin{itemize}
\item[(i)] There is a Zariski-dense subset $Z' \subset Z$ such that $\rho_x \prec_I \rho_z$ and 
${\rho_z}_{|I_{E_\ell}} \simeq {\rho_{z'}}_{|I_{E_\ell}}$ for all $z, z' \in Z'$ and $x \in X$.\ps
\item[(ii)] For each $z \in Z$ assume that $\rho_z$ is defined over $\Qb$ and that its Langlands 
correspondent $\pi_z$ is tempered. Then for any $z, z' \in Z$, we have 
$${\rho_z}_{|I_{E_\ell}} \simeq {\rho_{z'}}_{|I_{E_\ell}} \hspace{5mm} {\it and} \hspace{5mm}
{\pi_z}_{|\GL_n(\OO_{E_\ell})} 
\simeq {\pi_{z'}}_{|\GL_n(\OO_{E_\ell})}.$$
\end{itemize}
\end{lemma}

\begin{pf} Part (i) is \cite[Prop. 7.8.19]{bch} (in the notations there, apply it to the tautological 
pseudo-character $T : \GL_n(\OO(X)) \rightarrow \OO(X)$ and to the morphism $W_{E_\ell} \rightarrow 
\GL_n(\OO(X))$ of the statement). Let us check Part (ii). 

Let $\pi$ and $\pi'$ be two irreducible {\it tempered} representations of $\GL_n(E_\ell)$, with 
associated Weil-Deligne representations $(r,N)$ and $(r',N')$. If $(r_{|I},N)$ and $(r'_{|I},N')$ 
are isomorphic, then $\pi$ and $\pi'$ are fully induced from discrete series $\Delta$ and $\Delta'$ 
of a same Levi subgroup of $\GL_n(E_\ell)$, and $\Delta$ and $\Delta'$ only differ by a unitary 
unramified twist (Langlands, Zelevinski). In particular, $\pi$ and $\pi'$ are isomorphic restricted 
to $\GL_n(\OO_{E_\ell})$. As a consequence, the second assertion in part (ii) follows from the first 
one, that we check now. We may assume $z' \in Z'$. Moreover, by the inequality in part (i) it is enough 
to check that  ${\rho_z}_{|I_{E_\ell}} \simeq {\rho_{z'}}_{|I_{E_\ell}}$
holds after a finite base change. As base change preserves temperedness, we may assume that for 
each $z$ in $Z$, if $(r_z,N_z)$ is the Weil-Deligne representation associated to $\rho_z$, then 
$r_z$ is trivial over $I_{E_\ell}$.

Let $\varphi \in W_{{E_\ell}}$ be a lift of a geometric Frobenius and let $P \in
\OO(X)[T]$ be the 
characteristic polynomial of $\rho(F)$. Up to replacing $X$ by a finite covering and $Z$ by its 
inverse image, we may assume that $P$ splits in $\OO(X)[T]$; let $F_1$,\dots,$F_n \in \OO(X)$ be 
its roots.  By part (i) the conjugacy class of $N_z$ does not depend on $z \in Z'$, let us call 
it $N$; it is determined by a partition $n_1+n_2+\dots+n_s=n$. For each $z \in Z'$, there is a 
renumbering $F'_i$ of the $F_i$ such that whenever $1 \leq j <n_1$ or $n_1+n_2+\dots+n_r+1\leq j < n_1+n_2+\dots+n_{r+1}$ we have $F'_j(z)=qF'_{j+1}(z)$, where $q$ is the cardinal of the residue field of $E_\ell$. For $z$ in a Zariski-dense subset $Z''\subset Z'$, this numbering will be the same, hence up to renumbering the $F_i$ once and for all we may assume that if $1 \leq j <n_1$ or if $n_1+n_2+\dots+n_r+1\leq j < n_1+n_2+\dots+n_{r+1}$, 
\begin{equation}\label{eqtempn} F_j(z)=qF_{j+1}(z)\hspace{1cm}\,\,\forall z \in Z''.
\end{equation} 

Thus $F_{j+1}=qF_j$ in $\OO(X)$ (we may assume $X$ is reduced) for any $j$ as above, hence \eqref{eqtempn} 
holds as well for any $z \in Z$. But the temperedness of $\pi_z$ and the inequality $N_z \prec N$ 
imply in turn that $N_z$ is conjugate to $N$. (Note the funny interplay between the complex and 
$p$-adic sides in this proof.)
\end{pf}

Assume now that $E_\ell$ is a quadratic extension of $F_\ell$ ($\ell$ is still prime to $p$) and let $\U(F_\ell)$ be the quasi-split unitary group in $n$ variables attached to $E_\ell/F_\ell$. For $n \leq 3$, the standard base-change between $\U(F_\ell)$ and $\GL_n(E_\ell)$, conjectured by Langlands and proved 
by Rogaswki in \cite[Prop. 11.4.1, Thm. 13.2.1]{roglivre}, is an injection from the set of $L$-packets of 
$\U(F_\ell)$ to the irreducible smooth representations $\pi$ of $\GL_n(E_\ell)$ such that 
$\pi^{\vee,w} \simeq \pi$ ($w \in W_{F_\ell}\backslash W_{E_\ell}$) and whose central character is trivial over $F_\ell^*$ 
if $n$ is odd. 

\begin{lemma}\label{lemmehorspbis} We keep the assumptions of Lemma~\ref{lemmehorsp} and assume that $n \leq 3$.

For each $z \in Z$ assume that $\rho_z$ is defined over $\Qb$ and that there is a tempered $L$-packet $\Pi_z$ of $U(F_\ell)$ whose base-change to $\GL_3(E_\ell)$ corresponds to $\rho_z$. Then for any compact open subgroup $K_v$ of $U(F_\ell)$ and any $z,z' \in Z$, we have 
$$(\bigoplus_{\pi \in \Pi_z} \pi)_{|K_v} \simeq (\bigoplus_{\pi \in \Pi_{z'}} \pi)_{|K_v}.$$
\end{lemma}

If $\chi$ is a character of $E_\ell^*$, we shall denote by $\chi^\bot$ the character $x \mapsto \chi(c(x))^{-1}$ 
where $c$ generates $\Gal(E_\ell/F_\ell)$.

\begin{pf} Assume $n=3$. The lemma is obvious when the isomorphism class of $\rho_z$ is constant when $z$ varies 
in $Z$, as then $\Pi_z$ does not depend on $z$. We claim that this occurs in particular when there is some $z \in Z$ 
such that $\rho_z \otimes \Qpb$ is either irreducible, of the form $r \oplus \chi$ with $r$ irreducible of 
dimension $2$, or of the form $\chi_1\oplus\chi_2\oplus \chi_3$ where the $\chi_i$ are distinct characters such 
that $\chi_i^{\bot}=\chi_i$ for each $i$ (the associated $L$-packet of $\U(F_\ell)$ is discrete and has respectively 
$1$, $2$ or $4$ elements by \cite[\S 13]{roglivre}). 

Indeed, in the first case part (i) shows that for any $x \in X(\Qpb)$ there is an unramified character $\chi_x$ such 
that $\rho_x \otimes \Qpb = (\rho_z \otimes \Qpb)\otimes \chi_x$, and the self-duality like condition forces $\chi_x$ 
to vary in a finite set when $x$ varies. This implies that for all $g \in W_{E_\ell}$, ${\rm trace}(\rho(g))$ is 
{\it constant}, {\it i.e.} belongs to the biggest finite extension of $\Qp$ in the domain $\OO(X)$, hence the 
statement. The second case is similar. In the last case, part (i) implies that ${\rm trace}(\rho)$ factors through 
a pseudo-character of $W_{E_\ell}^{\rm ab}=E_\ell^*$. We may find a normal affinoid $Y$, a finite map $Y \rightarrow X$, 
as well as continuous characters $\tilde{\chi_i} : W_{E_\ell} \rightarrow \OO(Y)^*$ such that 
\begin{equation}\label{argext}{\rm trace}(\rho)=\sum_{i=1}^3 \tilde{\chi_i},\end{equation}
thus up to replacing $X$ by $Y$ and $Z$ by its inverse image in $Y$ we may assume that $X=Y$. 
The self-duality like condition and the assumption that the evaluation of the $\tilde{\chi_i}$ at $z$ are all 
self-$\bot$ and distinct imply that $\tilde{\chi_i}^\bot=\tilde{\chi_i}$ for each $i$. As they are constant on 
$I_{E_\ell}$ by part (i), we see again that there is only a finite number of possibilities for the $\rho_x\otimes \Qpb$ 
with $x \in X$, and we conclude as before.

In the non-constant case, then for any $z$ in $Z$ we have $$\rho_z^{\rm ss}\otimes \Qpb =\chi_z \oplus \nu_z \oplus 
\nu_z^\bot$$ for some characters $\chi_z, \nu_z : E_\ell^* \rightarrow \Qpb^*$, with $\chi_z^\bot=\chi_z$. This implies 
that $\chi_z$ is trivial over $F_\ell^*$ (as $\chi_z\nu_z\nu_z^\bot$ is), hence is the base change to $E_\ell^*$ of a unique 
character $\tilde{\chi}_z$ of $U(1)$. This case is more subtle as there are in general four possible kinds of 
$L$-packets for $\Pi_z$, labelled as (1), (2), (6) and (8) in \cite[p. 174]{roglivre} and to which we shall refer 
several times. We shall argue according to the generic monodromy of $\rho$ :

(a) If the monodromy of $\rho_z$ has nilpotent index $3$, then from that list we see that $\Pi_z$ is some twist of 
the Steinberg representation (case (8)). If this occurs for some $z \in Z'$, Lemma~\ref{lemmehorsp} (i) implies that 
this also occurs for all $z \in Z'$, so that $\Pi_z$ is independent of $z \in Z'$ in this case (the twist alluded 
above does not depend on $z$ by an argument similar to the ones above). Assume it is so and fix $x \in Z$. This 
implies 
\begin{equation}\label{casa} \{ \nu_x, \nu_x^\bot\}=\{\chi_x|.|, \chi_x|.|^{-1}\}.\end{equation} 
We claim that Rogawski's list shows that $\Pi_x$ is the twist of the Steinberg representation by $\tilde{\chi}_x$, 
so $\Pi_x=\Pi_z$ for any $z \in Z'$. This would conclude the proof in case (a). To check the claim note that 
property \eqref{casa} of $\nu_x$ excludes the packets of type (1) ({\it irreducible principal series}), as well 
as those of type (6) ({\it l.d.s. $L$-packets}), so it only remains to exclude the type (2) ({\it endoscopic 
transfert of a twist of the Steinberg representation of $U(2)\times U(1)$}, an $L$-packet with $2$ elements). 
But the Weil-Deligne representation $(r,N)$, $r : W_{E_\ell} \rightarrow \GL_3(\C)$, of the base-change of an 
$L$-packet of type (2) has a non trivial $N$ and $r$ is the sum of three characters
\begin{equation}\label{casendosc} \chi\oplus \nu|.|^{-1/2} \oplus \nu|.|^{1/2} \end{equation}
where $\chi=\chi^{\bot}$, $\nu=\nu^\bot$, $\chi_{|F_\ell^*}=1$ and $\nu_{|F_\ell^*} \neq 1$. These conditions 
contradict \eqref{casa}.\ps

(b) If the monodromy of $\rho_z$ has nilpotent index $2$, the unique possibility from the list is that 
the $L$-packet $\Pi_z$ is of type (2). If this occurs for some $z \in Z'$, an argument similar to the one 
in case (a) shows again that $\Pi_z$ is independent of $z \in Z$. \ps

(c) In the remaining case, $\rho_z$ has no monodromy for all $z \in Z$, so we see from Rogawski's 
description that $\Pi_z$ is the set of irreducible constituents of the parabolic induced representation 
${\rm Ind}(\eta_z)$ of the character $\eta_z=\nu_z\chi_z^{-1}\times {\tilde{\chi}}_z$ of the diagonal 
torus $E_\ell^* \times U(1)$ of $U(F_\ell)$ (this covers the two possibilities (1) and (6)). Arguing as for 
\eqref{argext}, we may assume that $\chi_z$ and $\nu_z$ have been chosen so that they vary analytically 
when $z \in Z \subset X$. But then ${\eta_z}_{|\OO_{E_\ell}^*\times U(1)}$ does not depend on $z \in Z$ by 
part (i), so we check at once that ${\rm Ind}(\eta_x)_{|K}={\rm Ind}(\eta_z)_{|K}$ for all $x, z \in Z$.

This concludes the case $n=3$. The case $n=2$ is similiar and only easier (use~the description of $L$-packets 
given in \cite[Chap. 11]{roglivre}) and $n=1$ is trivial.
\end{pf}

\begin{lemma}\label{lemmehorspter} Lemmas~\ref{lemmehorsp} and~\ref{lemmehorspbis} still hold if $\rho$ is replaced by the restriction to $W_{E_\ell}$ of a continuous pseudo-character $T : {\rm Gal}(\overline{E}/E)  \rightarrow \OO(X)$ of dimension $n$.
\end{lemma}

In this setting, if $x \in X(\Qpb)$ we mean by $\rho_x$ the restriction to $W_{E_\ell}$ of the (unique) semi-simple representation ${\rm Gal}(\overline{E}/E) \rightarrow \GL_n(\Qpb)$ whose trace is the evaluation of $T$ at $x$.

\begin{pf} Indeed, Lemmas~\ref{lemmehorsp} (ii) and~\ref{lemmehorspbis} only relied on the existence of a $\rho$ such that ${\rm trace}(\rho)=T$ via its corollary Lemma~\ref{lemmehorsp} (i). But part (i) in this more general setting follows again from~\cite[Prop. 7.8.19]{bch}.
\end{pf}

We now prove the main results of this section. Recall that $p$ splits in $E$ and that $U$ is a unitary group in $n$ variables attached to $E/F$ as in~\S\ref{infclasscrit}. Recall also that we have fixed a compact open subgroup $K \subset \U(\AAA_{F,f})$ that we assume now of the form $\prod_v K_v$. Let $S$ be a finite set of finite primes of $F$ containing the primes ramified in $E$ and such that $K_v$ is maximal hyperspecial for each $v \notin S$. Let $\mathcal{E}$ be the $p$-adic eigenvariety of $U$ of level $K$ associated to the set $S_p$ of all places above $p$ and to the global Hecke-algebra $\HH$ unramified outside $S$. 

\begin{thm}\label{thmetalenoncrit} ($n\leq 3$) Let $\pi$ be an automorphic representation of $U$ which is unramified above $p$ and such that $\pi^K \neq 0$. We assume that the base-change of $\pi$ to $\GL_n(\AAA_E)$ is cuspidal and that for each $v \in S_p$ the crystalline Frobenius has $n$ distinct eigenvalues on $\Dc(\rho_{\pi,\widetilde{v}})$.

Let $\{\Ref_v\}$ be a collection of refinements of the $\rho_{\pi,\widetilde{v}}$ for $v \in S_p$ and let $x \in \mathcal{E}$ be the point associated to $(\pi,\{\Ref_v\})$. If $\Ref_v$ is non-critical and regular for each $v \in S_p$, then $\mathcal{E}$ is \'etale at $x$ over the weight space.
\end{thm}

This solves conjecture (CRIT) of \cite{bch} \S 7.6 for $n\leq 3$, and to which we refer for a more complete discussion and motivations. The {\it non-critical} assumption is important, contrary to the {\it regular} one which should be removable (in Thm.~\ref{controlethm} already). 

%Recall that the accessible refinements of a non-tempered $\pi$ are automatically critical whereas any refinement of a tempered $\pi$ is accessible (see~\cite[\S 6.4.4]{bch} for this subtelty).

\begin{pf} For each closed point $y \in \mathcal{E}$ (resp. $y \in \mathcal{E}(\Qpb)$) we denote by $k(y)$ its residue field and by $\psi_y : \ATL^-\otimes \HH \rightarrow k(y)$ (resp. $\psi_y : \ATL^-\otimes \HH \rightarrow \Qpb$) the evaluation at $y$ of the structural ring homomorphism $\ATL^-\otimes \HH \rightarrow \OO(\mathcal{E})$ at $y$. 

By construction of the eigenvariety, we may find a connected affinoid neighborhood $B$ of the weight $\chi_0$ of $x$ (as defined in the statement) in the weight space and a finite locally free $\OO(B)$-module $M$ equipped with an $\OO(B)$-linear action of $\ATL^-\otimes \HH$ such that : 
\begin{itemize}
\item[(a)] The affinoid spectrum $V$ of ${\rm Im}\left( \OO(B)\otimes \ATL^- \otimes \HH \rightarrow \End_{\OO(B)}(M)\right)$ is an affinoid neighborhood of $x$ in $\mathcal{E}$. Denote by $\kappa : V \rightarrow B$ the natural map (so $\kappa(x)=\chi_0$). \ps
\item[(b)] For each algebraic weight $\chi \in B$, there is an isomorphism of $\Q_p\otimes\ATL^-\otimes \HH$-modules $$M \otimes_{\OO(B)} k(y) = \bigoplus_{y \in \kappa^{-1}(\chi)} F(\Cc_\chi)[\psi_y],$$
up to some normalizing twist depending on $\chi$ for the action of $\ATL^-$.
\item[(c)] $\kappa^{-1}(\kappa(x))^{\rm red}=\{x\}$ (in particular, $V$ is connected) and the natural surjection $\OO(V) \rightarrow k(x)$ has a section.
\end{itemize}

The classicity criterion implies that for $\chi \in B$ algebraic, dominant and sufficiently far from the walls of its Weyl chamber, $M \otimes_{\OO(B)} k(\chi) \subset F(W^\vee_\chi)$, so for $z \in \kappa^{-1}(\chi)$ we have  \begin{equation}\label{firstfactz}F(W_{\chi}^\vee)[\psi_z]=F(\Cc_{\chi})[\psi_z].\end{equation}
Let $Z_0 \subset B$ be a Zariski-dense subset of such weights, and let $Z \subset V$ be the union of $\{x\}$ and $\kappa^{-1}(Z_0)$, this a Zariski-dense subset of classical points of $V$. Up to reducing $Z_0$ if necessary, we may assume that the refinements associated to any $z \in Z$ and any $v \in S_p$ are regular, and that $\kappa$ is \'etale at each point is $Z\backslash\{x\}$.

The first important fact is that \begin{equation}\label{firstfact}M \otimes_{\OO(B)} k(\chi_0) = F(W_{\chi_0}^\vee)[\psi_x]=F(\Cc_{\chi_0})[\psi_x].\end{equation}
Indeed, by (b) and (c) above it is enough to check that $F(\Cc_{\chi_0})[\psi_x] \subset F(W_{\chi_0}^\vee)$ but this is Proposition~\ref{controlethm} as the $\Ref_v$ are regular and non-critical (note that $\varphi_{i,v}\varphi_{j,v}^{-1} \neq p$ as $\pi_v$ is tempered and unramified).

For each $z \in Z(\Qpb)$, let $\Pi_z$ be the global discrete $A$-packet of representations of $\U(F_v)$ which are unramified outside $S\backslash S_p$ and whose system of $\ATL^-\otimes \HH$-eigenvalues is $\iota_\infty\iota_p^{-1}(\psi_z)$. As the weights in $Z_0$ have been chosen outside the walls, and as $\pi$ is tempered by assumption, each $\Pi_z$ is actually a tempered $L$-packet. 
We deduce from \eqref{isomclass} that for each $z \in Z(\Qpb)$ 
\begin{equation}\label{eqdim} {\dim_{\Qpb}}
 \delta(z):=F(\Cc_{\kappa(z)})[\psi_z] = \sum_{\pi \in \Pi_z} m(\pi) \dim_{\C}(\pi_{S\backslash S_p}^{K_{S\backslash S_p}}),
\end{equation}
where $m(\pi)$ denotes the multiplicity of $\pi$ in the discrete spectrum of $U$. To obtain the formula above, 
we have used that $\pi_v^{K_v}$ has dimension $1$ for $v \notin S$ and that for $v \in S_p$ the generalized 
eigenspace of ${\psi_z}_{|\ATL^-}$ on $\pi_v^{K_v}$ is also one-dimensional as the $\varphi_{i,v}$ are distinct 
for each $v \in S_p$ by assumption.

We are going to show that for each $\chi \in Z_0$, $\kappa^{-1}(\chi)$ has a single element, whose residue 
field is $k(x)$. By (a), (b) and (c) this implies $\OO(V)=\OO(B) \otimes_{\Qp} k(x)$ and would conclude 
the proof. Putting together \eqref{firstfactz}, \eqref{firstfact}, the local
freeness of $M$ over $\OO(B)$, (b), (c), and \eqref{eqdim} above, it is 
enough to show that $\delta(z)$ does not depend on $z \in Z(\Qpb)$, or at least that $\delta(z)\geq \delta(x)$ 
for any $z$ in $Z(\Qpb)$. We shall use for this the natural $n$-dimensional pseudo-character $T_Y : 
{\rm Gal}(\overline{E}/E) \rightarrow \OO(Y)$ for each irreducible component $Y$ of $V$; these components all 
meet at $x$ by (c). 

We have not used so far that $\pi$ has a cuspidal base change to $\GL_n(\AAA_E)$. This implies 
that $\rho_\pi=\rho_x$ is absolutely irreducible (by Ribet, Blasius and Rogawski). As the absolutely 
irreducible locus is Zariski open in $V$, up to reducing $Z_0$ if necessary we may assume that for 
all $z \in Z(\Qpb)$ the packet $\Pi_z$ has a cuspidal base change as well. In this case, all the 
elements of $\Pi_z$ are automorphic and have multiplicity one by \cite[Thm. 13.3.3 (c), Thm. 14.6.1]{roglivre} 
(when $n=2$ this is due to Labesse-Langlands). By Lemmas~\ref{lemmehorsp}, \ref{lemmehorspbis} 
and~\ref{lemmehorspter}, we obtain that $\delta(z)$ is independent of $z \in Z(\Qpb)$, which 
concludes the proof. These lemmas apply as for each $\pi'$ with cuspidal base-change $\Pi'$, 
$\pi'$ is tempered and $\Pi' \mapsto \rho_{\pi'}$ is compatible with the Frobenius semi-simplified 
local Langlands correspondence at all primes prime to $p$.

Note that strictly, we have not considered the case of a split place $v \in S$ such that $U(F_v)$ is 
the group of units of a central simple algebra over $F_v$. But in this case the statement of 
Lemma~\ref{lemmehorspbis} still makes sense and holds: it follows from Lemma~\ref{lemmehorsp} (ii) and 
the fact that if $\pi_1$ and $\pi_2$ are two represesentations of $U(F_v)$ whose Weil-Deligne 
representations are inertially equivalent, then $\pi_1$ and $\pi_2$ differ by an unramified twist, 
hence coincide over the maximal compact subgroups of $U(F_v)$.\end{pf}

\begin{remark} {\rm A similar statement is probably true in many other cases (all?) when $\pi$ is only 
assumed to be tempered, but this would have forced us to look at Rogawski's multiplicity fomula in 
unpleasant detail. This is however simple enough if $S \backslash S_p$ does not contain any non split 
prime. Indeed, in this case the multiplicity formula~\cite[Thm. 14.6.5]{roglivre} implies that if $\Pi$ 
is a discrete tempered (possibily endoscopic) packet of $U$, then all the elements $\pi \in \Pi$ which are 
unramified at $p$ and such that $\pi^K \neq 0$ have multiplicity one, so the argument above still applies.} 
\end{remark}

We end this paragraph with a result valid for any $n$. Assume that $E/F$ is unramified everywhere and that 
$U$ is quasi-split at each finite place of $F$ (so $[F:\Q]$ is even if $n$ is even, by Hasse's principle 
for unitary groups). Assume that $S$ only contains primes which split in $E$. Let again $\mathcal{E}$ be 
the $p$-adic eigenvariety of $U$ of level $K$ associated to the set $S_p$ of all places above $p$ and to 
the global Hecke-algebra $\HH$ unramified outside $S$.

Let $\pi$ be an automorphic representation of $U$, and let $\rho_\pi$ the Galois representation associated 
to the base change $\Pi$ of $\pi$ to $\GL_n(\AAA_E)$. Consider the property:

(P)      $\Pi \mapsto \rho_{\pi}$ is compatible with the Frobenius semi-simplified local Langlands 
correspondence at all finite primes prime to $p$. 

This property (P) is known (\cite{GRFAbook},\cite{shin}) in many cases, for instance if $\Pi$ is cuspidal, 
unless $n$ is even and for each archimedean place $v$ of $F$ the infinitesimal 
character $\Pi_v$ coincides with the infinitesimal character of a one dimensional algebraic representation 
of $\GL_n$. A weak version of (P) holds in all cases by \cite{chharris}.

\begin{thm}\label{thmetalenoncrit2} Assume that $U$, $E/F$ and $K$ are as
above. Let $\pi$ be a tempered automorphic representation of $U$ which is unramified above $p$ and such that $\pi^K \neq 0$. We assume 
that $\pi$ has multiplicity $1$ in the discrete spectrum of $U$, that (P) holds for $\pi$, 
and that for each $v \in S_p$ the crystalline Frobenius has $n$ distinct eigenvalues on $\Dc(\rho_{\pi,\widetilde{v}})$.

Let $\{\Ref_v\}$ be a collection of refinements of the $\rho_{\pi,\widetilde{v}}$ for $v \in S_p$ and 
let $x \in \mathcal{E}$ be the point associated to $(\pi,\{\Ref_v\})$. If $\Ref_v$ is non-critical and 
regular for each $v \in S_p$, then $\mathcal{E}$ is \'etale at $x$ over the weight space.
\end{thm}

By a result of Labesse~\cite[Thm. 5.4]{labesse}, such a $\pi$ has multiplicity $1$ if its base-change to $\GL_n(\AAA_E)$ is 
cuspidal, in which case it is necessariliy tempered. 
It follows from Langlands' conjectures that any $\pi$ as in the statement should have multiplicity one.

\begin{pf} We start with an observation. Let $\pi$ and $\pi'$ be any two
automorphic representations of $U$. Assume that $\pi$ and $\pi'$ are
unramified at each finite place $v$ of $F$ which is inert in $E$, and that
$\pi_v \simeq \pi'_v$ if $v$ is archimedean and for almost all finite places
$v$ of $F$. We claim that $\pi \simeq \pi'$. Indeed, let $\Pi$ and
$\Pi'$ be their respective base change to $\GL_n(\AAA_E)$. By
Labesse~\cite[Cor. 5.3]{labesse}, $\Pi$ and $\Pi'$ are induced from discrete
automorphic representations of a Levi subgroup of $\GL_n(\AAA_E)$, and
$\Pi_v \simeq \Pi'_v$ for almost all the finite places $v$ of $E$. By
Moeglin-Waldspurger classification the the discrete spectrum of $\GL_m$ and
by Jacquet-Shalika \cite[Thm. 4.4]{JS}, this implies that $\Pi \simeq \Pi'$
(at all places). By Labesse's theorem above again, this implies that $\pi \simeq
\pi'$. Indeed, the global base change is compatible with the obvious local base change
at all the finite places. At the split places it is the identity and at the inert places
it is the spherical base change, so in both cases it is injective, and we are done. As a
consequence, if $\psi$ is a system of eigenvalues of $\HH \otimes \ATL^-$ on
$F(W_\chi^\vee)$, then there is a unique automorphic representation
$\pi(\psi)$ of $U$ of weight $\chi$ such that $\HH$ acts on $\pi(\psi)_f^S$ as $\psi$.

We now go back to the proof of the theorem. The arguments of the proof of Thm.~\ref{thmetalenoncrit} 
apply verbatim until ``We have not used so far that $\pi$ has a cuspidal
base-change'' if for any 
$z \in Z(\Qpb)$ we set $\Pi_z=\{ \pi(\psi_z)\}$. If $x$ and $\pi$ are as in
the statement then $\pi(\psi_x) \simeq \pi$. As $m(\pi(\psi_x))=1$ by assumption,
we have $$\delta(x):={\dim_{\Qpb}}
 F(\Cc_{\chi_0})[\psi_x] = \dim_{\C}(\pi_{S\backslash S_p}^{K_{S\backslash S_p}}).$$
By Lemma~\ref{lemmehorspter} applied to each irreducible component of $V$, we obtain that for 
each $z \in Z(\Qpb)$ we have $\delta(z)=m(\pi(\psi_z))\delta(x) \geq \delta(x)$, which concludes the proof. 
This lemma applies as $\pi(\chi_z)$ is tempered and satisfies (P) for each
$z \in Z(\Qpb)$: at $x$ it is the assumption, and at $z \neq x$ it holds as $\kappa(z)$ is not in a 
Weyl wall of $X^*(T)$.
\end{pf}

\section{End of the proof of Theorem A and other global applications}

\subsection{Proof of Theorem~\ref{redthm}}\label{mainpf} Fix $x \in \cX(\rhob)$ a modular point. For 
each place $v$ above $p$, set $V_v:={\rho_x}_{|\G_{\tilde{v}}}$. Assume that for each
$v$, $\End_{\G_{\tilde{v}}(V_v)}=L$, and that the eigenvalues of the crystalline Frobenius of $D_{\rm cris}(V_v)$ are distinct
and in $k(x)$. 

Let $y \in \mathcal{E}(\rhob)$ be a refined modular point above $x$. Those $y$ are in natural bijection with 
the set of collections refinements $\Ref_v$ of $V_v$, $v \in S_p$, in the sense of \S~\ref{triref}. Indeed, 
they have the form $y=(x,\delta)$ and this bijection is $$\delta \mapsto (\mathcal{F}_v)_{v\in S_p}=(\delta_{1,v}(p)p^{k_{1,v}},\delta_{2,v}(p)p^{k_{2,v}},\delta_{3,v}(p)p^{k_{3,v}})_{v\in S_p},$$
where $k_{1,v} < k_{2,v} < k_{3,v}$ are the Hodge-Tate weights of $V_v$. For this reason we shall also denote 
those $y$ above $x$ by $(x,\{\mathcal{F}_v\})$. We use the notations of~\S\ref{triref}.

\begin{prop}\label{trisat} Let $(x,\{\mathcal{F}_v\}) \in \mathcal{E}(\rhob)$ be a refined modular point such 
that $\mathcal{F}_v$ is non-critical for each $v \in S_p$. For any artinian thickening 
${\rm Spec}(A) \hookrightarrow \cE(\rhob)$ at $x$, the associated Galois representation 
$\rho_A : \G_E \rightarrow \GL_3(A)$ satisfies $${\rho_A}_{|\G_{\tilde{v}}} \in \cX_{V_v,\mathcal{F}_v}(A).$$
\end{prop}

Note that this statement makes sense as $\End_{\G_{\tilde{v}}(V_v)}=L$ by Remark~\ref{forgetpi}.

\begin{pf} This is \cite[Thm. 4.4.1]{bch}. \end{pf}

Let us set $\tg_{v,\ast}:=\cX_{V_v,\ast}(L[\varepsilon])$. Fix $y=(x,\{\mathcal{F}_v\})$ a refined 
modular point and consider the natural map on tangent spaces: 
\begin{equation}\label{lastlasteq} T_y(\cE(\rhob)) \longrightarrow T_x(\cX(\rhob)) \overset{({\rm res_{\widetilde{v}}})_{v \in S_p}}{\longrightarrow} \prod_{v \in S_p} \tg_v.\end{equation}
By the proposition above, we know that if $\mathcal{F}_v$ is non critical
for each $v$ then the image of this 
map falls inside $\prod_{v \in S_p} \tg_{v,\mathcal{F}_v}$. However, this information is still 
too weak to localize this image: the space on the left, which is of global nature, tends to have 
dimension $3\cdot |S_p|$, whereas $\prod_{v \in S_p} \tg_{v,\mathcal{F}_v}$, which is purely locally 
defined, has dimension $7\cdot |S_p|$. A key idea is to neglect the crystalline deformations in the 
range of (\ref{lastlasteq}), as they should be conjecturally transversal to the global ones by the 
Bloch-Kato conjecture (see~\cite[Conj. 7.6.5]{bch}). This is confirmed by the following proposition.

\begin{prop}\label{propfin}Assume again that $\mathcal{F}_v$ is non-critical for each $v$, then the 
linear map (\ref{lastlasteq}) induces an isomorphism $$T_{(x,\{\Ref_v\})}(\cE(\rhob)) \isomo 
\prod_{v \in S_p} \tg_{v,\mathcal{F}_v}/\tg_{v,{\rm crys}}.$$
\end{prop}

\begin{pf} By Theorem~\ref{thmetalenoncrit} and Prop.~\ref{trisat}, the composite of the natural maps 
$$T_{(x,\{\Ref_v\})}(\cE(\rhob)) \longrightarrow \prod_{v\in S_p} \tg_{v,\mathcal{F}_v}/\tg_{v,{\rm crys}} \overset{ \pht}{\longrightarrow} \prod_{v\in S_p} L^{3}$$
is an isomorphism. But by Prop~\ref{maindefref} the second map is an isomorphism, hence so is the first one.
\end{pf}

\begin{cor} The image of 
$\bigoplus_{y=(x,\{\Ref_v\})} T_y(\cE(\rhob)) \longrightarrow T_x(\cX(\rhob))$
has $k(x)$-dimension at least $6[F:\Q]$.
\end{cor}

\begin{pf} Consider the natural map 
$$\bigoplus_{y=(x,\{\Ref_v\})} T_y(\cE(\rhob)) \longrightarrow T_x(\cX(\rhob)) \longrightarrow \prod_{v\in S_p} \tg_v/\tg_{v,{\rm crys}}$$
where we restrict in the sum on the left to those $\{\Ref_v\}$ such that $\Ref_v$ is non-critical for all $v$. By Theorem~\ref{mainlocalbis}, Example~\ref{dim3well} and Prop.~\ref{propfin}, the composite of the two maps above is surjective. The result follows as the space on the right has dimension ${\frac{n(n+1)}{2}}\cdot |S_p|=6 \cdot {|S_p|}$ by Prop.~\ref{maindefref}.
\end{pf}

\subsection{Generalisations and other settings}\label{genoset} Assume now that $n$ is any
odd (say) integer, and keep the notations of~\S\ref{prelimdefo}. Assume again
that $\rhob$ is modular. We expect that the strategy employed when $n=3$ will lead to a proof of the natural variant of Thm.\ref{mainthmobs} for any $n$ (with $3$ and $6$ replaced by $n$ and $\frac{n(n+1)}{2}$). It would actually follow by the same
proof if we could prove the following conjecture on the  genericity of
global Galois representations that we believe in:

\begin{conj}\label{generconj} ({\it Genericity conjecture}) For any modular point $x \in \cX(\rhob)$ and any affinoid
neighborhood $U$ of $x$ in $\cX(\rhob)$, there exists $y \in U$ such that
$\rho_{y,\widetilde{v}}$ is generic for each $v \in S_p$. 
\end{conj}

Of course, we expect to find such a point $y$ by starting form $x$ and moving in the piece of the infinite fern 
inside $U$. It would even be enough for our purpose to have a similar statement with
``is generic'' replaced by ``has a nested sequence of non-critical
refinements''. When $n=3$, Lemma~\ref{lemmazd} and Example~\ref{dim3well}
give a positive answer to this weak version, but even in this case the above
conjecture is still open as far as we see. Note that it would not be
too difficult to show that we may
find an $y$ in $U$ such that the Zariski-closure of the image of $\rho_{y,\widetilde{v}}$
is the whole of $\GL_n(\overline{k(x)})$ (extending Lemma~\ref{lemmazd} for any
$n$). However, this does not seem to imply that $\rho_{y,\widetilde{v}}$
is generic enough. In particular, there does not seem to be any variant of
Example~\ref{dim3well} when $n>3$. We hope to come back to this conjecture in the future. 
Here is another natural question. 

\begin{question} {\it How far geometric Galois representations can be non-generic ?
Are there some (even conjectural) global conditions implying the genericity at each place above
$p$ ?}
\end{question}

The example we have in mind is the classical conjecture that if a $p$-adic
Galois representation $\rho_f : G_{\Q} \longrightarrow \GL_2(\Qpb)$ attached
to a classical modular eigenform $f$ is split at $p$, then the form $f$ is CM.

In another direction, our infinitesimal approach also shed some lights back to the
standard case studied by Gouv\^ea and Mazur. The following result is a
simple consequence of our method and of Prop~\ref{transversality}.

\begin{prop}\label{retoursurgm} Let $\rho_f : G_{\Q,S} \rightarrow \GL_2(\Qpb)$ be a $p$-adic
Galois representation attached to a modular eigenform $f$ of weight $k>1$
and level prime to $p$. Assume that ${\rho_f}_{|G_{|\Q_p}}$ is
indecomposable and that the two eigenvalues of its crystalline Frobenius are
distinct. Then the two leaves of the infinite fern of Gouv\^ea-Mazur cross
transversally at $\rho_f$. 
\end{prop}

 Let us end this paper by considering the Hilbert modular analogue
of Thm~\ref{mainthmobs}, to which the arguments of Gouv\^ea-Mazur do not extend as well when $[F:\Q]>1$ (as far as we know).

Here $F$ is a totally real field which is {\it totally split at the odd prime $p$}, $S$ a finite set of places of $F$ containing the places above $p$ and $\infty$, 
and $\rhob : G_{F,S} \rightarrow \GL_2({\mathbb{F}_q})$ is absolutely irreducible
and totally odd. In this context, we say that $\rho$ is modular if it is isomorphic to the $p$-adic Galois
representation $\rho_\Pi$ attached to a cuspidal automorphic Galois
representation $\Pi$ of $\GL_2(\AAA_F)$ which is cohomological, unramified above $p$ and
outside $S$ (hence to a Hilbert modular eigenform), and we say that $\rhob$ is modular if it is isomorphic to the residual representation of a modular $\rho$. Let $\cX(\rhob)$ be the generic fiber
of the universal $G_{F,S}$-deformation of $\rhob$, as defined by Mazur. In the unobstructed case $H^2(G_{F,S},{\rm ad}(\rhob))=0$, then $\cX(\rhob)$ is the open unit ball over $\Q_q$ in $1+2[F:\Q]$ variables. This unobstructedness assumption also implies that $H^2(G_{F,S},\F_p)=0$, thus Leopold's defect $\got{d}$ of $F$ at $p$ vanishes. To avoid assuming Leopold's conjecture in general it is convenient to work with the notion of an {\it essentially modular $\rho$}, which we define as {\it the twist $\rho=\rho_\Pi \otimes \chi$ of a modular $\rho_\Pi$ by a continuous character $\chi: G_{F,S} \rightarrow \Qpb^*$ unramified above $p$}. We obtain this way notions of modular and essentially modular points in $\cX(\rhob)$. For simplicity we assume moreover that $[F:\Q]$ is even (see Rem.\ref{remhilb}).

\begin{thm}\label{cashbeigen} Assume that $\rhob$ is modular. Then the irreducible components of the Zariski-closure of the essentially modular points in $\cX(\rhob)$ all have dimension at least $1+\got{d}+2[F:\Q]$.

If $H^2(G_{F,S},{\rm ad}(\rhob))=0$ then the modular points are Zariski-dense in $\cX(\rhob)$.
\end{thm}
Of course, an important ingredient in the proof is the Hilbert modular eigenvarieties. These eigenvarieties have been studied by several authors: Hida, Kassaei, Buzzard, Kisin-Lai, Emerton, and Yamagami. We shall mostly rely on Buzzard's results.

Define ${\cXem}\subset \cX(\rhob)$ as the subset of essentially modular points, and ${\cXgem} \subset \cXem$ as the subset of $x$ such that for each $v|p$, $\rho_{x,v}$ is absolutely irreducible and with distinct crystalline Frobenius eigenvalues in $k(x)$ (in particular it is generic, note that each $\rho_{x,v}$ is crystalline by definition).

There is a finite index subgroup $\Gamma \subset \OO_F^*$ of the totally positive units such that the universal character $G_{F,S}^{\rm ab} \rightarrow R(\rhob)^*$ is trivial over $\Gamma$ via the reciprocity map. Let $\TT$ and $\WW$ be the rigid analytic spaces over $\Q_p$ parameterizing respectively the continuous $p$-adic characters of $(F_p^*)^2$ and of $(\OO_{F_p}^*)^2$ which are trivial on $\Gamma$ diagonally embedded. The space $\WW$ is called the {\it weight space}. It has equidimension $\dim(\WW)=1+\got{d}+[F:\Q]$. The $\rhob$-Hilbert modular eigenvariety $$\cE(\rhob) \subset \cX(\rhob) \times \TT$$ 
is the Zariski closure of the pairs $(x,\Ref)$ of refined essentially modular points (using the same translation as in \S\ref{mainpf} between refinements of an essentially modular $\rho$ and elements of $\TT$). 

$\cE(\rhob)$ has been studied by Buzzard in~\cite[\S III]{buzzard}, by switching to a totally definite quaternion algebra $D$ over $F$ which is split at all the finite places of $F$ and using the Jacquet-Langlands correspondence. Such a $D$ exists as $[F:\Q]$ is even. The important properties of $\cE(\rhob)$ are the following:

(i) {\it $\cE(\rhob)$ is equidimensional of dimension $\dim(\WW)$. The natural map to the weight space $\kappa: \cE(\rhob) \rightarrow \TT \rightarrow \WW$ is locally finite. The set of $(x,\Ref)$ with $x \in \cXem$ is a Zariski-dense accumulation subset of $\cE(\rhob)$. The same holds for the $(x,\Ref)$ with $x$ modular if (and only if) $\got{d}=0$.}

These properties are not stated this way {\it loc. cit}, but they follow simply from the construction there and an argument similar to the one of Thm.\ref{thmeigen}. There is a classicality criterion in this case as well. Arguing as in~Thm.\ref{thmetalenoncrit2} (note that the strong multiplicity one property holds for $D^*$), we obtain that: 

(ii) {\it If $x \in {\cXem}$ is such that $\rho_{x,v}$ has distinct crystalline Frobenius eigenvalues at each $v|p$, and if $\Ref_v$ is a non-critical refinement of $\rho_{x,v}$ for each $v|p$, then $\kappa$ is \'etale at $(x,\{\Ref_v\})$.}

Moreover, the following (easier) variant of Lemma~\ref{lemmazd} holds:

(iii) {\it Let $x \in {\cXem}$ and assume that a collection $\Ref=(\Ref_v)_{v\in S_p}$ of refinements of the $\rho_{x,v}$ has the property that for each $v|p$ such that $\rho_{x,v}$ is ordinary, then the first Frobenius eigenvalue of $\Ref_v$ is the one with the greatest $p$-adic valuation. Then the $(y,\Ref')$ with $y \in {\cXgem}$ accumulate at $(x,\Ref)$ in $\cE(\rhob)$}.

\begin{pf} (of Thm.\ref{cashbeigen}) The proof is only easier than in the $\U(3)$ case, so we shall be rather sketchy. 

Consider some $x \in {\cXgem}$. Let $(x,\Ref) \in \cE(\rhob)$ be an associated refined point and consider the composite of the natural maps  
\begin{equation}\label{cashilb1} T_{(x,\Ref)}(\cE(\rhob)) \longrightarrow T_x(\cX(\rhob)) \rightarrow \prod_{v\in S_p} \tg_v/\tg_{v,\rm crys}\end{equation}
as in \eqref{lastlasteq}. Using the \'etaleness of $\cE(\rhob) \rightarrow \WW$ at $(x,\Ref)$, Prop.\ref{maindefref}, and arguing as in~\S\ref{mainpf}, we see that the image of \eqref{cashilb1} is exactly the subspace $\tg_\Ref\cap \tg^0$ where $\tg_\Ref=\prod_{v\in S_p} \tg_{v,\Ref_v}/\tg_{v,\rm crys}$ and $\tg^0 \subset \prod_{v\in S_p} \tg_v/\tg_{v,\rm crys}$ is the subspace parameterizing the deformations of $(\rho_{x,v})_{v\in S_p}$ whose determinant, viewed as a morphism $F_p^* \rightarrow k(x)^*$, factors through $\Gamma$. Remark that 
$$\dim_{k(x)} t^0 = 3[F:\Q]-([F:\Q]-1-\got{d})=1+\got{d}+2[F:\Q].$$
Moreover $\sum_{\Ref} t_\Ref\cap t^0 = t^0$ by Thm.\ref{mainlocal} (remark that $\cap_\Ref t_\Ref \rightarrow (\prod_{v \in S_p}t_v)/t^0$ is surjective as  
``families of twists are $\Ref$-trianguline for each $\Ref$''). It follows that the image of the natural map  $$\bigoplus_{\Ref}T_{(x,\Ref)}(\cE(\rhob)) \rightarrow T_x(\cX(\rhob))$$
has $k(x)$-dimension at least $1+\got{d}+2[F:\Q]$. The first statement follows by the same argument as in \S\ref{stratpar}. The second statement (unobstructed case) follows as well by property (i) of $\cE(\rhob)$ and by the local-global numerical coincidence $$\dim \cX(\rhob) = \dim  \tg^0.$$
\end{pf}

\begin{remark}\label{remhilb} {\rm When $[F:\Q]$ is odd, the same result holds, with the same proof, if we restrict to essentially modular points $\rho_\Pi\otimes \chi$ such that $\Pi_w$ is ess. square integrable at some finite place $w \in S\backslash S_p$, as we may switch to a suitable definite quaternion algebra in this case as well. However there may sometimes be no such modular point (e.g. when $S=S_p \cup S_\infty$). We could actually still conclude in most cases in general using some well chosen quadratic base-change and results of Kisin-Lai, but this would lead us too far away from our purposes here.}
\end{remark}

\section{An application to adjoint Selmer groups}\label{adjointprime}

\subsection{Bloch-Kato Selmer groups} Let $p$ be a prime, $L$ a finite
extension of $\Q_p$. In what follows, an $L$-representation of a topological group $G$
will always mean a continuous, finite-dimensional, $L$-linear
representation. Let $\ell$ be a prime, $M$ a finite extension 
extension of $\Q_\ell$, and let $V$ be an $L$-representation of $G_M$. Recall that Bloch and Kato
defined in \cite[\S 3]{bkato} a subspace $H^1_f(M,V) \subset H^1(M,V)$. Set
$H^1_s(M,V)=H^1(M,V)/H^1_f(M,V)$. Let $F$ be a number field, 
let $S$ be a finite set of finite places of $F$ containing the set $S_p$ of
places dividing $p$, the set $S_\infty$ of archimedean places, and let $V$ be an $L$-representation of $G_{F}$ unramified outside $S$. Let 
$H^1_{f'}(F,V)$ be the kernel of the natural map 
$H^1(G_{F,S},V) \rightarrow \prod_{v \in S\backslash S_p} H^1_s(F_v,V)$, it
does not depend on $S$ as above. In
what follows, we shall be mostly interested in the surjectivity of the natural map:

\begin{equation}\label{glolo} H^1_{f'}(F,V) \rightarrow \prod_{v| p} H^1_s(F_v,V).
\end{equation}

The kernel of this map is usually denoted by $H^1_f(F,V)$. We may like viewing the surjectivity 
of \eqref{glolo} as some splitting of the Poitou-Tate exact 
sequence. By the global duality theorem, the image of \eqref{glolo} is Tate's orthogonal complement 
of the image of $H^1_f(F,V^\vee(1)) \rightarrow \prod_{v|p} H^1_f(F_v,V^\vee(1))$. In particular, \eqref{glolo} 
is surjective if and only if the natural map
$H^1_{f}(F,V^\vee(1)) \rightarrow \prod_{v | p} H^1(F_v,V^\vee(1))$ is the zero map. 

\begin{remark}\label{rembk}{\rm We shall be mainly concerned with cases where $V$ is geometric
(say in the strongest sense) of pure
weight $0$. In this case, the conjectures of Bloch-Kato
(\cite[\S 3.4]{FP},\cite{bkato}) assert that $H^1_f(F,V)=0$, so \eqref{glolo} is
expected to be injective. Assume furthermore that $F$ is totally
real, let $c_v$ denote a complex conjugation at the archimedean place
$v$, and denote by $h_v^{(0,0)}$ the dimension of the $(0,0)$-part of the
$F_v$-Hodge structure associated to $V$. If $${\rm trace}({c_v}_{|V})=-h_v^{(0,0)}$$ 
for each real place $v$, then those conjectures predict as well that
$H^1_f(F,V^\vee(1))=0$ (compare
them with~\cite[\S 3]{serre}). In particular, \eqref{glolo} 
should be an isomorphism in this case ! As far as we know, very few is known about this conjecture, possibly nothing when $\dim(V) \geq 3$ before the results of this paper. }
\end{remark}

\subsection{The {\sc Adjoint'} Selmer group of a Galois representation of type
$\U(n)$}\label{adselmer}

Assume now that $F$ is a totally real field. Let $E$ be a totally imaginary quadratic 
extension of $F$ and assume that $p$ is totally split in $E$ (hence in $F$). For each 
$v | p$ fix some place $\widetilde{v}$ of $E$ above $v$, so $F_v=E_{\widetilde{v}}=\Q_p$. 
Let $c \in G_F$ be a complex conjugation.

Let $\Pi$ be a cuspidal automorphic representation of $\GL_n(\AAA_E)$ such that 
$\Pi^\vee \simeq \Pi^c$ and such that $\Pi_v$ is algebraic regular at each archimedean 
place $v$ of $E$. Let $\rho_\Pi$ be the $p$-adic Galois representation of $G_E$ associated 
to $\Pi$ and $\iota$, say with coefficients in $L$ and normalized so that 
$\rho_\Pi^\vee \simeq \rho_\Pi(n-1)$. We shall assume that:

\begin{itemize}
\item[(i)] For each $v | p$, $\Pi_v$ is unramified, its Langlands class has $n$ distinct eigenvalues, and
$\End_{G_{\widetilde{v}}}({\rho_{\Pi}}_{|G_{\widetilde{v}}})=L$,\ps

\item[(ii)] $\Pi$ is tempered and $\Pi \mapsto \rho_\Pi$ satisfies (P) of~\S\ref{noncritetalepar}.
\end{itemize}

Recall that (ii) is actually known to be automatically satisfied except in the rather special case where: $n$ is even, each $\Pi_v$ for $v$ archimedean has the same infinitesimal character as a one dimensional algebraic representation of $\GL_n$, and there is no finite place $v$ such that $\Pi_v$ is square-integrable. Note moreover that (i) implies that $\rho_\Pi$ is absolutely irreducible. In what follows, we shall assume furthermore that $L$ is big enough so as to contain (via $\iota$) the eigenvalues alluded in (i), which is harmless.

We are going to define below an $L$-representation $\adr(\rho_\Pi)$ of $G_F$
of dimension $n^2$ that extends the standard representation ${\rm ad}(\rho_\Pi)$ of $G_E$. From the self-duality like
condition, there is a $P \in \GL_n(L)$, unique up to $L^*$, such that 
$${}^t\rho_{\Pi}(cgc^{-1})^{-1}=P\rho_\Pi(g) P^{-1}\chi(g)^{n-1}, \forall g \in G_E$$ ($\chi : G_F \rightarrow \Q_p^*$ is the cyclotomic character). It follows that ${}^t P = {\pm} P$, and the main theorem of \cite{bchsign} ensures that $P$ is a {\it symmetric} matrix. Following~\cite[\S 2]{cht}, let us consider the linear algebraic $\Q$-group
$$\mathcal{G}_n:=(\GL_n \times \GL_1)\rtimes \Gal(E/F),$$
where $c$ acts on $\GL_n \times \GL_1$ by $(g,\lambda) \mapsto
({}^{t}g^{-1}\lambda,\lambda)$. We check at once that there exists a unique
morphism $\widetilde{\rho_\Pi}: G_F \rightarrow \mathcal{G}_n(L)$ such that
$\widetilde{\rho_\Pi}(g)=(\rho_\Pi(g),\chi(g)^{1-n})$ for $g \in G_E$ and
$\widetilde{\rho_\Pi}(c)=({}^t P^{-1},1)c$. We denote by $$\adr(\rho_\Pi)$$ 
the representation of $G_F$ defined by the adjoint representation of 
$\widetilde{\rho_\Pi}$ on $M_n(L)={\rm Lie}(\GL_n \times \{1\}) \otimes L \subset {\rm Lie}(\mathcal{G}_n)\otimes L$. The $G_F$-equivalence class of $\adr(\rho_\Pi)$ only depends on the $G_E$-equivalence class of $\rho_\Pi$. The map $(X,Y)\mapsto {\rm trace}(XY)$ defines a $\mathcal{G}_n(L)$-equivariant pairing $M_n(L) \otimes M_n(L) \rightarrow L$, thus $\adr(\rho_\Pi)$ is selfdual. Note that the complex conjugation $c$ acts on $M_n(L)=\adr(\rho_\Pi)$ as $X \mapsto - P{}^t\!X P^{-1}$. In particular, $\adr(\rho_\Pi)^{G_F}=0$, and the diagonal matrices induce an embedding $\varepsilon_{E/F} \rightarrow \adr(\rho_\Pi)$, where $\varepsilon_{E/F}$ is the character of order $2$ of $G_F$ associated to $E/F$.

In what follows we shall be mainly interested in the Selmer group
$$H^1_{f'}(F,\adr(\rho_\Pi)).$$
The purity assumption (ii) on $\rho_\Pi$ ensures that
$H^0(F_v,\adr(\rho_\Pi)^\vee(1))$ vanishes, hence so do
$H^2(F_v,\adr(\rho_\Pi))$ and $H^1_s(F_v,\adr(\rho_\Pi))$, for each finite place
$v$. In particular, for each finite set $S$ containing $S_p$, $S_\infty$ and
the ramification of $\adr(\rho_\Pi)$, we have 
\begin{equation}\label{setS}
H^1(G_{F,S},\adr(\rho_\Pi))=H^1_{f'}(F,\adr(\rho_\Pi)).
\end{equation} 

\medskip

\begin{conj}\label{conjadr} If $V=\adr(\rho_\Pi)$, then \eqref{glolo} is an isomorphism,
and $H^1_f(F,V)=H^1_f(F,V^\vee(1))=0$.
\end{conj}

Indeed, we are in the case of Remark~\ref{rembk}. Note that for each real
place $v$, $h_v^{(0,0)}=n$ by the description of the Hodge-Tate numbers of
$\rho_\Pi$, and ${\rm trace}({c_v}_{|V})=-n$ by the result of \cite{bchsign}
recalled above.

%\footnote{Indeed, if $v \notin S_p$ is a finite place of $F$ unramified in $E$, and if $w$ is a place of $E$ above $v$ such that $\Pi_w$ is unramified, then $H^1_f(F_v,\rho_\Pi)=H^1(F_v,\rho_\Pi)$ (each class is unramified) as the cardinal of the residue field of $E_w$ is not the quotient of two eigenvalues of $\rho_\Pi({\rm Frob}_w)$ by (ii). The other cases are similar}  

The {\it adjoint'} \,Selmer group of $\rho_\Pi$ has a natural description as the tangent space of a suitable deformation functor. Fix an $S$ as above and let $\mathcal{C}$ be the category defined in~\S\ref{setdef}. For an object $A$ of $\mathcal{C}$ let $\cX_{\rho_\Pi}(A)$ be the set of $A$-isomorphism classes of continuous representations $\rho_A : G_{E,S} \rightarrow \GL_n(A)$ such that $\rho_A \otimes_A L \simeq \rho_\Pi$ and such that ${\rm trace}(\rho_A^\vee)={\rm trace}(\rho_A^c)\chi^{n-1}$.
This defines a functor $\cX_{\rho_\Pi}: \mathcal{C} \rightarrow {\rm Sets}$, equipped with a natural morphism 
$$\cX_{\rho_\Pi} \rightarrow \prod_{v\in S_p} \cX_{{\rho_{\Pi,\widetilde{v}}}}=:\cX_p,$$
where we have set $\rho_{\Pi,\widetilde{v}}:={\rho_{\Pi}}_{|G_{\widetilde{v}}}$. This allows to define a collection of subfunctors of the left-hand side by pulling back the subfunctors of the right-hand side studied in~\S\ref{tridef}. Let $\Ref=\{\Ref_v\}$ be a collection of refinements of the $\rho_{\Pi_{\widetilde{v}}}$ for $v \in S_p$. We set

$$\cX_{\rho_\Pi,\Ref}=\cX_{\rho_\Pi}\times_{\cX_p} \prod_{v|p} \cX_{\rho_{\Pi_{\widetilde{v}}},\Ref_v} \, \, \, {\rm and}\, \, \, \cX_{\rho_\Pi,f}=\cX_{\rho_\Pi}\times_{\cX_p} \prod_{v|p} \cX_{\rho_{\Pi_{\widetilde{v}}},{\rm crys}}.$$
Note that conditions (i), (ii), (iii) and (iv) of~\S\ref{triref} on $\rho_{\Pi,\widetilde{v}}$ are satisfied by conditions (i) and (ii) above (and by the known description of its Hodge-Tate weights for (iv)). Similar functors were previously introduced in~\cite[\S 7.6]{bch}.

\begin{prop}\label{proprepglob} $\cX_{\rho_\Pi} : \mathcal{C} \rightarrow
{\rm Sets}$ is pro-representable and its tangent space is canonically isomorphic to $H^1_{f'}(F,\adr(\rho_\Pi))$.

Moreover, for $\ast \in \{ f, \Ref\}$, $\cX_{\rho_\Pi,\ast}$ is a pro-representable subfunctor of $\cX_{\rho_\Pi}$. The tangent space of $\cX_{\rho_\Pi,f}$ coincides with $H^1_f(F,\adr(\rho_\Pi))$. When $\Ref_v$ is non-critical for each $v\in S_p$, then $\cX_{\rho_\Pi,f}\subset \cX_{\rho_\Pi,\Ref}$. 
\end{prop}

\begin{pf} The functor $\cX_{\rho_\Pi}$ is pro-representable by Mazur as $\rho_\Pi$ is absolutely irreducible (the condition ${\rm trace}(\rho_\ast^\vee)={\rm trace}(\rho_\ast^c)\chi^{n-1}$ being obviously relatively representable). The identification of its tangent space with $H^1(G_{F,S},\adr(\rho_\Pi))$ is similar to Lemma~\ref{calctanun} (using the $\mathcal{G}_n$ here) and is checked in~\cite[\S 2]{cht}. The other assertions follow from Prop.~\ref{maindefref}.
\end{pf}

\begin{definition}\label{defhunF} For each collection of refinements $\Ref=\{\Ref_v\}$ for $v \in S_p$, define $H^1_\Ref(F,\adr(\rho_\Pi))$ (resp. $H^1_{\Ref_v}(F_v,\adr(\Pi))$ as the tangent space of  $\cX_{\rho_\Pi,\Ref}$ (resp. $\cX_{\rho_{\Pi,\widetilde{v}},\Ref_v}$). 
\end{definition}

\subsection{The main results}\label{mainthmgeneraln} ${}^{}$ \ps  Fix $\Pi$ as in~\S\ref{adselmer}. The following statement (in some special cases) was conjectured in~\cite[\S 7.6]{bch}. 

\begin{thm}\label{thmadselmer1} Assume that $\Ref_v$ is non-critical and regular for each $v \in S_p$. Then the natural map 
$$H^1_\Ref(F,\adr(\rho_\Pi)) \longrightarrow \prod_{v|p} H^1_{\Ref_v}(F_v,\adr(\rho_\Pi))/H^1_f(F_v,\adr(\rho_\Pi))$$
is surjective. In particular, $\dim_L H^1_\Ref(F,\adr(\rho_\Pi))=\dim_L H^1_f(F,\adr(\rho_\Pi))+n[F:\Q]$.
\end{thm}

\begin{pf} Let us choose $F'/F$ a totally real quadratic extension of $F$. Assume that $F'/F$ is split above $p$, that $F'_v=E_v$ for each finite place $v$ of $F$
which is either ramified in $E$ or which is inert and such that $\Pi_v$ is
ramified. Let $\Pi'$ be Arthur-Clozel's quadratic base change of $\Pi$ to $E \cdot F'$, it is cuspidal by (i). As the $f$ and $\Ref$ conditions can be checked after any finite base-change which is split above $p$, the inflation-restriction sequence induces isomorphisms
$$H^1_{\ast}(F,V)=H^1_{\ast}(F',V_{|G_{F'}})^{\Gal(F'/F)}$$
for $\ast \in \{ f', f, \Ref\}$. As a consequence, it is
enough to show the surjectivity of the map of the statement when $F,E$ and $\Pi$ are replaced by $F',E\cdot F'$ and $\Pi'$. In particular, we may assume that $[F:\Q]$ is even, that $E/F$ is unramified everywhere, and that the places $w$ of $E$ such that $\Pi_w$ is ramified are split over $F$.  (This kind of trick has been used in another context by
Blasius-Rogawski and Harris-Taylor.)

	As $[F:\Q]$ is even, there exists a unitary group $U/F$ as in~\S\ref{infclasscrit}
which is furthermore quasi-split at all finite places. By Labesse~\cite[Thm. 5.4]{labesse}, $\Pi$ admits a strong descent to an automorphic representation
$\pi$ of $U$ which has furthermore multiplicity $1$. Let $S$ be a finite set of places of $F$ which split in $E$, containing $S_p$, and such that $\pi_f$ is unramified outside $S$. Let $K=\prod_v K_v$ be a compact open subgroup of $\U(\AAA_{F,f})$ as in~\S\ref{infclasscrit} such that $K_v$ is hyperspecial for each $v \notin S$ and small enough so that $\pi_f^K \neq 0$. Let $X$ be the eigenvariety of $U$ of level $K$ associated to $\iota$ and the spherical Hecke
algebra outside $S$. Let $\{\Ref_v\}_{v \in S_p}$ be as in the statement and let $x \in X$ be the point associated to $(\pi,\{\Ref_v\})$. Up to extending the scalars in $X$, me may assume that $X$ is defined over the field $L$ of~\S\ref{adselmer}. The natural pseudo-character $G_{E,S} \rightarrow \OO(X)$ (see e.g. \cite[Prop. 7.5.4]{bch}) and the absolute irreducibility of $\rho_\Pi$ define a canonical $L$-linear map
on tangent spaces
\begin{equation}\label{maptaneigen} T_x(X) \longrightarrow \cX_{\rho_\Pi}(L[\varepsilon])=H^1(G_{F,S},\arb(\rho_\Pi))=H^1_{f'}(F,\adr(\rho_\Pi)). \end{equation}
By~\cite[Thm. 4.4.1]{bch} (using that $\Ref_v$ is non-critical and regular for each $v$), the image of \eqref{maptaneigen} falls inside $H^1_\Ref(F,\adr(\rho_\Pi))$. Consider now the induced map:   
\begin{equation}\label{maptaneigen2} T_x(X) \longrightarrow \prod_{v|p} H^1_{\Ref_v}(F_v,\arb(\rho_\Pi))/H^1_f(F_v,\arb(\rho_\Pi)). \end{equation}
\noindent We claim that this map is an isomorphism, which will conclude the proof. Note that the subspace $T_x(X) \subset H^1_\Ref(F,\adr(\rho_\Pi))$ even furnishes a canonical section of the map of the statement. Consider the natural map $$\prod_{v|p} H^1_{\Ref_v}(F_v,\arb(\rho_\Pi))/H^1_f(F_v,\arb(\rho_\Pi))\overset{(\pht)_v}{\longrightarrow} (L^n)^{S_p}.$$
By Prop.~\ref{maindefref} this is an isomorphism, so it is enough to check that the composite of~\eqref{maptaneigen2} with this map is an isomorphism. But this composite is the derivative of the map from $X$ to the weight-space at $x$, which is an isomorphism by Thm.~\ref{thmetalenoncrit2}, and we are done.
\end{pf}

\begin{thm}\label{mainthmselmer} Assume that $\forall \, \,v | p$, $\rho_{\Pi,\widetilde{v}}$ is weakly generic and regular (\S~\ref{examplelocalmain}). Then $$H^1_{f'}(F,\adr(\rho_\Pi)) = \sum_{\Ref} H^1_{\Ref}(F,\adr(\rho_\Pi))$$
and the natural map 
$$H^1_{f'}(F,\adr(\rho_\Pi)) \longrightarrow \prod_{v|p} H^1_s(F_v,\adr(\rho_\Pi))$$
is surjective. 
\end{thm}

For each $v |p$, let $\Ref_{i,v}$  be a weakly nested sequence of non-critical regular refinements of $\rho_{\Pi,\widetilde{v}}$. The equality $H^1_{f'}(F,\adr(\rho_\Pi)) = \sum_{\Ref} H^1_{\Ref}(F,\adr(\rho_\Pi))$ even holds if we restrict the sum to the $n^{[F:\Q]}$ refinements $\Ref$ such that $\Ref_v$ is one of the $\Ref_{i,v}$ for each $v$.

\begin{pf} It is an immediate consequence of Thm.~\ref{mainlocalbis} and Thm.~\ref{thmadselmer1}.
\end{pf}

\begin{cor}\label{hdeuxhunf} Under the assumptions of Theorem~\ref{mainthmselmer}, and for each $S$ as in~\eqref{setS}, 
$$\dim_L H^1_f(F,\adr(\rho_\Pi))=\dim_L H^2(G_{F,S},\adr(\rho_\Pi)).$$
In particular, if $H^2(G_{F,S},\overline{\adr(\rho_\Pi)})=0$ for some
residual $\overline{\F}_p$-representation $\overline{\adr(\rho_\Pi)}$ associated to $\adr(\rho_\Pi)$,
then $H^1_f(F,\adr(\rho_\Pi))=0$.
\end{cor}

\begin{pf} By Thm.~\ref{mainthmselmer} and Prop.~\ref{maindefref}, we have 
\begin{equation}\label{dimintermediaire}\dim_L H^1_{f'}(F,\adr(\rho_\Pi))-\dim_L H^1_f(F,\adr(\rho_\Pi))=\frac{n(n+1)}{2}[F:\Q].\end{equation}
As we already saw, $H^0(G_{F,S},\adr(\rho_\Pi))=0$ and $\adr(\rho_\Pi)(c)(X)=-P^{-1}{}^t\!XP$ for some symmetric invertible matrix $P$ by~\cite{bchsign}, so $$\dim_L \adr(\rho_\Pi)^{c=-1} = \frac{n(n+1)}{2}.$$
This numerical coincidence with the right-hand side of~\eqref{dimintermediaire}, as well as \eqref{setS} and the global Euler characteristic formula conclude the proof.
\end{pf}

Theorem~\ref{mainthmselmer} is very general as a lot of Selmer groups are direct summands of Adjoint' Selmer groups. Here is an important special case.

Let $F$ be a totally real field in which $p$ is totally split and let $\Pi$ be a regular algebraic cuspidal 
automorphic representation of $\GL_n(\AAA_F)$ such that $\Pi^\vee \isomo \Pi \eta$ for some Hecke character 
$\eta$ of $F$ such that the sign $\eta_\infty(-1):=\eta_v(-1)$ is
independent\footnote{Actually, Langland's conjectures imply that this last
condition should follow from the others.} of the archimedean place $v$ 
of $F$. Denote by $q \in \Z$ the integer such that $\eta |.|^{-q}$ is an
Artin character of $F$. Let $\rho_\Pi : G_F \rightarrow \GL_n(L)$ be the $p$-adic Galois representation associated to $\Pi$ 
and $\iota$ (\cite{chharris}), and denote by $\eta_\iota : G_F \rightarrow
L^*$ the character associated to $\eta$ and $\iota$. Recall that $\rho_\Pi^\vee \simeq \rho_\Pi
\eta_\iota \chi^{n-1}$ 
and that $\rho_\Pi$ is symplectic if $n$ is even and
$\eta_\infty(-1)(-1)^q=1$, and orthogonal otherwise by \cite[Cor. 1.3]{bchsign}. Assume that $\Pi$ satisfies (i) and (ii) 
of~\S\ref{adselmer} (with $\widetilde{v}$ replaced by $v$ in (i), again (ii) is known in most cases). 
Let $E$ be any totally imaginary quadratic extension of $F$ in which $p$ is totally split.

\begin{thm}\label{thmselmertotreel} Let $\Pi$ be as above and assume that for each $v|p$, 
${\rho_\Pi}_{|G_{F_v}}$ is weakly generic regular. 

(Symplectic case) If $n$ is even and $\eta_\infty(-1)(-1)^q=1$, then \eqref{glolo} is surjective 
for $V=\Lambda^2 \rho_\Pi \otimes \eta_\iota\varepsilon_{E/\Q}\chi^{n-1}$ and $V=\Sym^2 \rho_\Pi \otimes
\eta_\iota \chi^{n-1}$. 

(Orthogonal case) Otherwise, \eqref{glolo} is surjective for $V=\Lambda^2 \rho_\Pi \otimes
\eta_\iota\chi^{n-1}$ and $V=\Sym^2 \rho_\Pi \otimes \eta_\iota \varepsilon_{E/\Q}\chi^{n-1}$. 

\end{thm}

\begin{pf} Indeed, let $\Pi'$ be Arthur-Clozel's base change $\Pi'$ of $\Pi$ to $E$, it is 
cuspidal by property (i) of $\Pi$. By~\cite[Lemma 4.1.4]{cht}, there is an algebraic Hecke 
character $\nu$ of $E$ which is unramified at $p$ and such that $\Pi' \otimes \nu$ is conjugate 
self-dual. This latter representation clearly satisfies properties (i) and (ii) of~\S\ref{adselmer} 
by the corresponding assumptions on $\Pi$. A simple computation left to the reader shows that 
$$\adr(\rho_{\Pi'\otimes \nu}) =  \Lambda^2 \rho_\Pi \otimes \eta_\iota\chi^{n-1} \oplus \Sym^2 \rho_\Pi \otimes
\eta_\iota \varepsilon_{E/\Q}\chi^{n-1}$$
in the first case, and 
$$\adr(\rho_{\Pi'\otimes \nu}) =  \Lambda^2 \rho_\Pi \otimes \eta_\iota \varepsilon_{E/\Q}\chi^{n-1} \oplus \Sym^2 \rho_\Pi \otimes
\eta_\iota\chi^{n-1},$$
in the second.
\end{pf}

Let $f$ be a classical modular eigenform of weight $k>1$ and let $\rho_f :
G_{\Q,S} \rightarrow \GL_2(\Qpb)$ be the $p$-adic Galois representation associated to $f$ and $\iota$. Assume that the level of $f$ is prime to $p$ and let $E$ be a quadratic imaginary field split at $p$.

\begin{thm}\label{corQthm} Assume that ${\rho_f}_{|G_{\Qp}}$ is absolutely indecomposable and that the ratio of its two crystalline Frobenius eigenvalues is not a root of unity.
Then \eqref{glolo} is an isomorphism for:\begin{itemize}
\item[(i)] $V={\Sym^n}(\rho_f)\det(\rho_f)^{-n/2}$ if $n=2,6$. \ps
\item[(ii)] $V={\Sym^n}(\rho_f)\det(\rho_f)^{-n/2}\varepsilon_{E/\Q}$ if
$n=0,4$.
\end{itemize}
Moreover, for each such $V$ we have $H^1_f(\Q,V)=H^2(G_{\Q,S},V)$.
\end{thm}

\begin{pf} Let $\Pi_0$ be the cuspidal automorphic representation of 
$\GL_2(\AAA_\Q)$ associated to $\rho_f$ and let 
$\eta$ be its central character. We normalize $\Pi_0$ so that 
$\eta = |.|^{2-k}$ times a Dirichlet character, so $\rho_f=\rho_{\Pi_0}$. As is well known, we 
have $\eta_\infty(-1)(-1)^{k-2}=1$. 

By Kim-Shahidi~\cite{KSh} the representation 
$\Pi:={\rm Sym}^3 \Pi_0$ of $\GL_{4}(\AAA_\Q)$ defined at all places 
by the local Langlands correspondence is automorphic and cuspidal (note that
$f$ is not CM by the assumption at $p$). By 
assumption on $\Pi_0$, $\Pi$ satisfies conditions (i) and (ii) 
of~\S\ref{adselmer}, and we have $\Pi^\vee \simeq \Pi \eta^{-3}$. Note that
$\Pi$ is symplectic as $(\eta_\infty(-1)(-1)^{(k-2)})^3=1$. The first part of the theorem is then 
a consequence of Thm.\ref{thmselmertotreel} applied to $\Pi$. Its
assumptions are satisfied by Examples~\ref{casexsym} and~\ref{symtrois}. The last sentence follows
from the global Euler characteristic formula as in the proof of Cor.~\ref{hdeuxhunf}.
\end{pf}

\section*{Appendix:  Some unobstructed Galois representations of type $\U(3)$}\label{appendix}

Let $A$ be an elliptic curve over $\Q$. Let $p\geq 5$ be a prime of good reduction of $A$ and assume that the representation 
$$\G_\Q \rightarrow {\rm
Aut}A[p] \simeq \GL_2(\F_p)$$ is surjective. Let $E$ be a quadratic imaginary field and $S$ the set of places containing $\infty$, $p$ and the primes dividing
$\disc(E){\rm cond}(A)$. Assume to simplify that $p$ splits in $E$ and that $\gcd(\disc(E),{\rm cond}(A))=1$. By 
Example~\ref{exmodf}, the representation $$\rhob={\rm Symm^2}A[p](-1)_{|\G_{E,S}}$$ 
is modular of type $\U(3)$. Let $\arb$ be the $\G_{\Q,S}$-module associated to $\rhob$ as in~\S\ref{prelimdefo}. \ps

\begin{lemma} There is an isomorphism of $\F_p[\G_{\Q,S}]$-modules :
$$\arb \simeq \rhob \oplus \epsilon_{E/\Q} \oplus ({\rm
Symm}^4A[p])(-2)\otimes \epsilon_{E/\Q},$$ 
where $\epsilon_{E/\Q}$ is the non-trivial character of $\Gal(E/\Q)$. 
\end{lemma}

%\begin{pf} Indeed, consider the natural representation $V=\F_p^2$ of $G=\GL_2(\F_p)$. As
%$V^* = V \det(V)^{-1}$, the representation $W=\Symm^2(V) \det(V)^{-1}$ is orthogonal. Thus we may
%define a representation $\tau : G \times \Gal(E/\Q) \longrightarrow
%\mathcal{G}_3(\F_p)$ extending $W$ and such that $\tau(\Gal(E/\Q)) \subsetneq \GL_3(\F_p)$. An immediate computation shows that
%$\mathfrak{g}_3 \simeq \Lambda^2(W)\oplus \epsilon_{E/F} \Symm^2(W)$. But
%$\Lambda^2(W) \simeq W$ and $\Symm^2 \Symm^2 V = \Symm^4 V \oplus
%\det(V)^2$ as $p\geq 5$, so
%%for $p=3$, it is a sum of copies of the Steinberg representation. 
%$\mathfrak{g}_3 = \Symm^2(V)\det(V)^{-1} \oplus \epsilon_{E/F}(1\oplus
%\Symm^4 V \det(V)^{-2})$.\end{pf}

We are looking for a set of conditions ensuring that $H^2(\G_{\Q,S},\arb)=0$.  Let $K$ be the field of definition of some $\F_p$-line in $A[p]$ (so
$[K:\Q]=p+1$), $\chi : G_K \rightarrow \F_p^*$ the natural character on the quotient of $A[p]$ by this
line, and $K'$ the extension of $K$ cut out by the character $\chi^4\epsilon_{E/\Q}(-1)$. If $\ell$ is a prime of good reduction of $A$, 
we set $a_\ell=\ell+1-|A(\F_\ell)|$.

\begin{prop} If the following conditions are satisfied, then $\rhob$
is unobstructed. \ps

(i) $\forall \ell \in S$, $H^0(\Q_\ell, \arb(1))=0$.

(ii) $p$ does not divide the degree of a modular parameterization of $A$, 

(iii) The class numbers of $E$ and $K'$ and prime to $p$.

\noindent Moreover, $H^0(\Q_\ell, \arb(1))=0$ if : 

(a) $\ell=p$, unless $A$ is ordinary at $p$, $A[p]_{|\G_{\Q_p}}$ is split, and $a_p = \pm 1 \bmod p$.

(b) $\ell \, | \,{\rm disc}(E)$, $\ell \neq 1 \bmod p$ and $a_\ell \neq \pm (\ell+1) \bmod p$.

(c) $\ell \, | \,{\rm cond}(A)$, $A$ is semistable at $\ell$, $A[p]_{|\G_{\Q_\ell}}$ is non-split, $\ell^2$ and  $\ell^3 \epsilon_{E/\Q}(\ell) \neq  1 \bmod p$.

\end{prop}

\begin{pf} As $\OO_E^*\otimes \F_p={\rm Cl}(E)\otimes \F_p=0$ we have $\Sha^1_S(\G_{\Q,S},\epsilon_{E/\Q}(1))=0$, thus 
assumption (i) ensures first that $H^2(G_{\Q,S},\epsilon_{E/\Q})=0$. 
By (i) again, $H^0(\Q_\ell,\Symm^2A[p])=0$ for each $\ell \in S$, thus assumption (ii) implies the vanishing of 
$H^2(G_{\Q,S},\rhob)$ by a result of Flach~\cite[\S 3]{flach}. We now deal with the last term 
$U:=\Symm^4A[p](-2)\otimes \epsilon_{E/\Q}=U^*$. Note that $H^0(G_{\Q,S},U^*(1))=0$ by assumption (i), so 
Tate's global duality theorem shows that $H^2(G_{\Q,S},U)=0$ if and only if $\Sha^1_S(G_{\Q,S},U(1))=0$. 

Let $B \subset \GL_2(\F_p)$ denote the upper triangular subgroup
and $\chi : B \rightarrow \F_p^*$ the unique quotient of $V_{|B}$ (in the
notations of the proof of the above lemma). We have an $\F_p[G]$-equivariant injection
$\Symm^4 V\rightarrow {\rm Ind}_B^G (\chi^4)$ and the quotient $Q$ is well known to be irreducible of dimension $p-4$, isomorphic to 
$\det(V)^2$ when $p=5$. In particular, in all cases $H^0(\G_{\Q,S},Q(-1)\otimes \epsilon_{E/\Q})=0$ and we have natural injections
$$H^1(\G_{\Q,S},U(1)) \longrightarrow H^1(G_{K,S},\chi^4(-1)\epsilon_{E/\Q}) \longrightarrow  H^1(G_{K',S},\F_p).$$
Via this injection, $\Sha^1_S(\G_{\Q,S},U(1))$ is mapped into ${\rm Cl}(K')\otimes \F_p$, hence (iii). 

\noindent The second part of the statement is an immediate computation and is left to the reader.\end{pf}

We use the notations of Cremona's tables of elliptic curves~\cite{cremona}.

\begin{prop} (Under GRH) Let $p=5$, $E=\Q(i)$ and assume that $A$ belongs to one of the isogeny classes 
$$17A, \, \, 21A, \, \, 37B,\, \,  39A,\, \,  51A,\, \,  53A,\, \,  69A, \, \, 73A, \, \, 83A,\,\mbox{and }\,\, 91B.$$
Then $H^2(\G_{\Q,S},\arb)=0$.
\end{prop}
%28 curves
These curves are exactly the elliptic curves $A$ of odd, square-free, conductor $N < 100$ such that each prime divisor of $N$ is congruent to $\pm 2 \bmod 5$, 
and such that $a_2(A) \neq \pm 2$ (this rules out the classes $37A$, $43A$, $67A$ and $91A$). These last two conditions are the ones in (b) and part of (c) above, and 
are actually necessary local conditions for the unobstruction of the deformation functor of $A[5]$ itself. 

The third row of the following table gives the degrees of the prime isogenies in a class, note that there is no isogeny of degree $5$. The sixth row gives the $j$ 
invariant of a certain element in the class, namely of : $17A4$, \,$21A4$, \,$37B3$, 
\,$39A4$, \,$51A1$, \,$53A1$, \,$69A1$, \,$73A2$, \,$83A1$, and $91B1$. Their valuation at any $\ell$ dividing 
${\rm cond}(A)$ is $<5$, so $A[p]_{|\G_{\Q_\ell}}$ is not split for any such $\ell$. This also shows that 
$\G_\Q \rightarrow {\rm Aut}(A[5])$ is surjective for each $A$. Using the table below, we see that the criteria (a), (b) and (c) of the proposition apply, hence (i) holds.

The seventh row gives the modular degree of the strong Weil curve in each class, (ii) follows (recall that there is no $5$-isogeny within a class). 

For each class we computed the field $K$ and found $K=\Q(x)$ with : {\tiny
$x^6 + 4x^5 + 17x - 17=0$,  $x^6 - 2x^5 + 9x^4 - 12x^3 + 21x^2 - 15x + 3=0$, $x^6 - 2x^5 - 5x^4 + 10x^3 - 15x^2 + 8x - 36=0$, 
$x^6 - 2x^5 - 5x^4 + 20x^3 - 25x^2 - 25x + 75=0$, $x^6 - 2x^5 - 5x^4 + 20x^3 - 25x^2 - 25x + 75=0$, $x^6 - x^5 - 5x^4 + 20x^3 - 15x^2 - 46x + 101=0$, 
$x^6 - x^5 - 5x^4 - 80x^3 + 35x^2 + 184x + 1521=0$, $x^6 - 2x^5 - 5x^4 + 20x^3 - 25x^2 - 59x + 143=0$, $x^6 - 3x^5 - 5x^4
 - 10x^3 + 10x^2 + 75x + 245=0$}, and {\tiny $x^6 - 2x^5 - 5x^4 + 30x^3 - 35x^2 + 20x + 80=0$.} The last row below gives the computation by Pari \cite{pari} of the 
class number of $K'=K(\cos(\pi/10))$, which concludes the proof. This last computation
(and only this one) depends on GRH; it would be interesting to make it
unconditional!

\begin{center}
{\small
\begin{tabular}{|c||c|c|c|c|c|c|c|c|c|c|c|} 
\hline
class &17A & 21A & 37B & 39A & 51A & 53A & 69A & 73A & 83A & 91B \\
\hline
cardinal & $4$ & $6$ & $3$ & $4$ & $2$ & $1$ & $2$ & $2$ & $1$ & $3$ \\
\hline
isogenies & $2$ & $2$ & $3$ & $2$ & $3$ & & $2$ & $2$ & & $3$\\
\hline
$a_2$ & $-1$ & $-1$ & $0$ & $1$ & $0$ & $-1$ & $1$ & $1$ & $-1$ & $0$ \\
\hline
$a_5$ & $-2$ & $-2$ & $0$ & $2$ & $3$ & $0$ & $0$ & $2$ & $-2$ & $-3$ \\
\hline
$j$ & {\tiny $\frac{3^3 11^3}{17}$} & {\tiny $\frac{47^3}{3.7}$} & {\tiny $\frac{2^{15} 5^3}{37}$} & {\tiny $\frac{23^3}{3.13}$} & 
{\tiny$\frac{2^{15}}{3^317}$} & {\tiny $\frac{3^5 5^3}{53}$} & {\tiny $\frac{-5^6}{3^2 23}$} 
& {\tiny $\frac{3^3 19^3}{73}$} & {\tiny $\frac{47^3}{83}$} & {\tiny $\frac{-2^{15}11^3}{7.13}$} \\
\hline
degree & 1 & 1 & 2 & 2 & 2 & 2 & 2 & 3 & 2 & 4 \\
\hline
$h_{K'}$ & $2$ & $2$ & $8$ & $2$ & $8$ & $2$ & $2$ & $32$ & $2$ & $2$ \\
\hline
\end{tabular}
}
\end{center}

\end{document}